\documentclass[12pt,twoside]{article}
\usepackage{amsmath,amsfonts,amssymb,a4,enumitem,epsfig,array,psfrag}
\usepackage{url}
\usepackage{amsthm}
\usepackage{etoolbox} 
\usepackage{color}
\usepackage{calc}
\usepackage{tikz, graphics} 
\usepackage{stmaryrd} 

\numberwithin{equation}{section}

\theoremstyle{plain}
\newtheorem{theo}{Theorem}
\newtheorem{prop}{Proposition}[section]

\newtheorem{coro}[prop]{Corollary}

\newtheorem{lemma}[prop]{Lemma}
\newtheorem{question}[prop]{Question}

\theoremstyle{definition}

\AtEndEnvironment{definition}{\null\hfill$\diamond$}%
\newtheorem{remark}[prop]{Remark}
\AtEndEnvironment{remark}{\null\hfill$\diamond$}%
\newtheorem{example}[prop]{Example}
\AtEndEnvironment{example}{\null\hfill$\diamond$}%

\newcommand{\NN}{{\mathbb{N}}}
\newcommand{\ZZ}{{\mathbb{Z}}}

\newcommand{\RR}{{\mathbb{R}}}

\newcommand{\Ss}{{\mathbb{S}}}
\newcommand{\BB}{{\mathbb{B}}}
\newcommand{\DD}{{\mathbb{D}}}

\newcommand{\FF}{{\mathbb{F}}}
\newcommand{\KK}{{\mathbb{K}}}

\newcommand{\cE}{{\cal E}}

\newcommand{\fa}{{\mathfrak a}}

\newcommand{\bF}{{\mathbf F}}
\newcommand{\bL}{{\mathbf L}}
\newcommand{\bN}{{\mathbf N}}
\newcommand{\bQ}{{\mathbf Q}}
\newcommand{\bS}{{\mathbf S}}
\newcommand{\bp}{{\mathbf p}}
\newcommand{\cU}{{\mathcal U}}
\newcommand{\cO}{{\mathcal O}}

\newcommand\thin{\operatorname{thin}}
\newcommand\thick{\operatorname{thick}}
\newcommand\anti{\operatorname{anti}}

\newcommand\Cliff{\operatorname{Cl}}
\newcommand\NL{\operatorname{NL}}
\newcommand\nc{\operatorname{nc}}
\newcommand\Flag{\operatorname{Flag}}

\newcommand\GL{\operatorname{GL}}
\newcommand\SL{\operatorname{SL}}
\newcommand\SO{\operatorname{SO}}
\newcommand\Spin{\operatorname{Spin}}
\newcommand\so{\operatorname{\mathfrak{so}}}
\newcommand\spin{\operatorname{\mathfrak{spin}}}
\newcommand\Trace{\operatorname{Trace}}
\newcommand\diag{\operatorname{diag}}
\newcommand\Diag{\operatorname{Diag}}
\newcommand\HQuat{\operatorname{HQuat}}
\newcommand\Quat{\operatorname{Quat}}
\newcommand\Pos{\operatorname{Pos}}
\newcommand\Lo{\operatorname{Lo}}
\newcommand\lo{{\mathfrak{lo}}}
\newcommand\Up{\operatorname{Up}}

\newcommand{\BL}{\operatorname{BL}}
\newcommand{\BLS}{\operatorname{BLS}}
\newcommand{\BLC}{\operatorname{BLC}}
\newcommand{\Bru}{\operatorname{Bru}}
\newcommand\inv{\operatorname{inv}}
\newcommand\Inv{\operatorname{Inv}}
\newcommand\sign{\operatorname{sign}}

\newcommand{\Block}{\operatorname{Block}}
\newcommand{\block}{\operatorname{block}}
\newcommand{\nmesmo}{\llbracket n \rrbracket}
\newcommand{\nmaisum}{\llbracket n+1 \rrbracket}

\newcommand\fl{{\mathfrak{l}}}

\begin{document}
\title{On the homotopy type of \\
intersections of two real Bruhat cells}

\author{Em{\'\i}lia Alves
\footnote{
Em{\'\i}lia Alves,
emiliaalves@id.uff.br,
Departamento de Matem\'atica Aplicada,
Instituto de Matem\'atica e Estat{\'\i}stica,
Universidade Federal Fluminense,
Rua Professor Marcos Waldemar de Freitas Reis s/n,
Niterói, RJ 24210-201, Brazil}
\and Nicolau C. Saldanha
\footnote{Nicolau C. Saldanha, saldanha@puc-rio.br,
Departamento de Matem\'atica, PUC-Rio,
Rua Marqu\^es de S\~ao Vicente 225,
Rio de Janeiro, RJ 22451-900, Brazil}
}



\maketitle

\begin{abstract}
Real Bruhat cells give an important and well studied stratification
of such spaces as $\GL_{n+1}$, $\Flag_{n+1} = \SL_{n+1}/B$,
$\SO_{n+1}$ and $\Spin_{n+1}$.
We study the intersections of a top dimensional cell
with another cell (for another basis).
Such an intersection is naturally identified with a subset
of the lower nilpotent group $\Lo_{n+1}^{1}$.
We are particularly interested in the homotopy type of such intersections.
In this paper we define a stratification of such intersections.
As a consequence, we obtain a finite CW complex
which is homotopically equivalent to the intersection.

We compute the homotopy type for several examples.
It turns out that for $n \le 4$ all connected components 
of such subsets of $\Lo_{n+1}^1$ are contractible:
we prove this by explicitly constructing the corresponding CW complexes.
Conversely, for $n \ge 5$ and the top permutation,
there is always a connected component with even Euler characteristic,
and therefore not contractible.
This follows from formulas for the number of cells per dimension
of the corresponding CW complex.
For instance, for the top permutation $S_6$,
there exists a connected component with Euler characteristic equal to $2$.
We also give an example of a permutation in $S_6$ for which 
there exists a connected component which is homotopically equivalent
to the circle $\Ss^1$.
\end{abstract}


\section{Introduction}

For a permutation $\sigma \in S_{n+1}$,
let $P_\sigma$ be the corresponding permutation matrix,
with entries $(P_\sigma)_{i,i^\sigma} = 1$ and $0$ otherwise.
The top permutation is denoted by $\eta$,
so that $i^\eta = n+2-i$
(for all $i \in \nmaisum = \{1, 2, \ldots, n, n+1\}$).
Let $\Lo_{n+1}^1$ be the nilpotent group 
of real lower triangular matrices
with diagonal entries equal to $1$.
Following the Bruhat decomposition,
partition $\Lo_{n+1}^1$ into subsets $\BL_\sigma$ for $\sigma \in S_{n+1}$:
\begin{equation}
\label{equation:BL}
\BL_\sigma = \{ L \in \Lo_{n+1}^1 \;|\;
\exists U_0, U_1 \in \Up_{n+1}, L = U_0 P_\sigma U_1 \}.
\end{equation}
Here $\Up_{n+1}$ is the group of real upper triangular matrices
with non zero diagonal entries,
a Borel subgroup of $\GL_{n+1}$.
Notice that $\BL_\eta$ is open and dense.
The set $\BL_\sigma$ is naturally homeomorphic to
the intersection of two Bruhat cells with different basis,
one of the of top dimension: see Section \ref{section:history}.
The aim of this paper is to study the homotopy type 
of the sets $\BL_\sigma$.
The following theorem sums up some of our main results.

\begin{theo}
\label{theo:one}
Consider $\sigma \in S_{n+1}$ and $\BL_\sigma \subset \Lo_{n+1}^1$
as in Equation \eqref{equation:BL}.
\begin{enumerate}
\item{For $n \le 4$, every connected component of every set
$\BL_\sigma$ is contractible.}
\item{For $n = 5$ and $\sigma = 563412 \in S_6$,
there exist connected components of $\BL_\sigma$
which are homotopically equivalent to the circle $\Ss^1$.}
\item{For $n \ge 5$, there exist connected components 
of $\BL_\eta$ which have even Euler characteristic.}
\end{enumerate}
\end{theo}

We also list the number of connected components per permutation.
Here $563412$ denotes the permutation $\sigma$ with
$1^\sigma = 5$, $2^\sigma = 6$, $3^\sigma = 3$ and so on until $6^\sigma = 2$. 

\begin{figure}[ht]
\begin{center}
\includegraphics[scale=0.25]{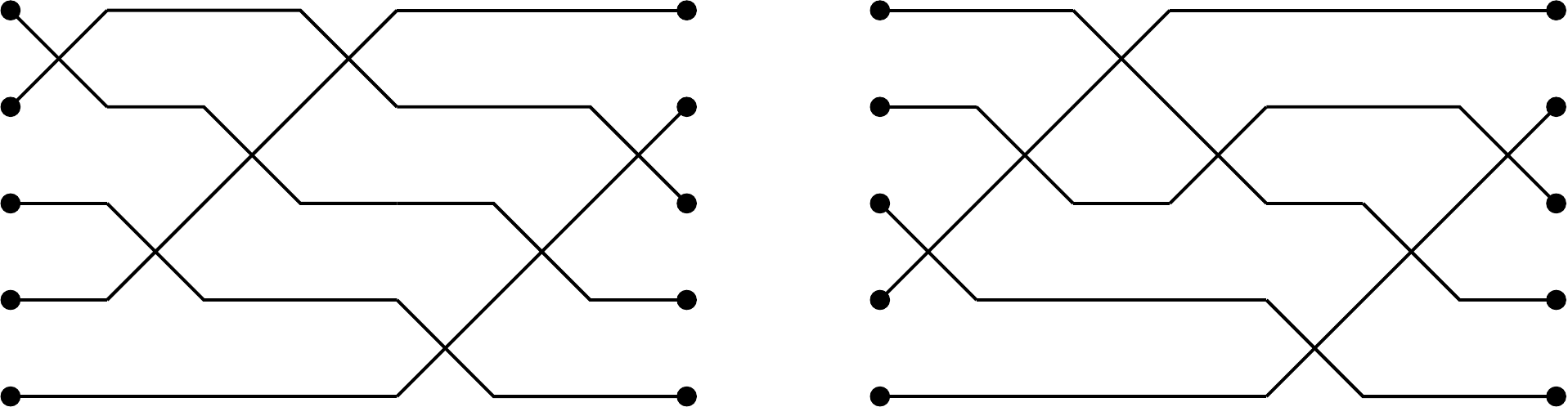}
\end{center}
\caption{The words
$a_1a_3a_2a_1a_4a_3a_2$ and $a_3a_2a_1a_2a_4a_3a_2$
are both reduced and represent the same permutation
$\sigma = 43512 \in S_5$, $\inv(\sigma) = 7$.}
\label{fig:reduced}
\end{figure}

In order to compute the homotopy type of $\BL_\sigma$,
we first choose a \textit{reduced word}
(also called \textit{reduced decomposition})
for $\sigma \in S_{n+1}$:
\begin{equation}
\label{equation:reduced}
\sigma = a_{i_1} \cdots a_{i_\ell}, \qquad \ell = \inv(\sigma).
\end{equation}
We extensively use the standard Coxeter-Weyl  generators
$a_1, \ldots, a_n$ of the symmetric group $S_{n+1}$;
$a_i$ is the simple transposition $(i,i+1)$.
Also, $\inv(\sigma)$ is the number of inversions of $\sigma$,
so that a reduced word is a word of shortest length.
A reduced word for a permutation is represented by a wiring diagram,
as in Figure \ref{fig:reduced}.
(Notice that there are different conventions in the literature;
in our system, each crossing is a generator,
read left-to-right.) 

\begin{figure}[ht]
\begin{center}
\includegraphics[scale=0.25]{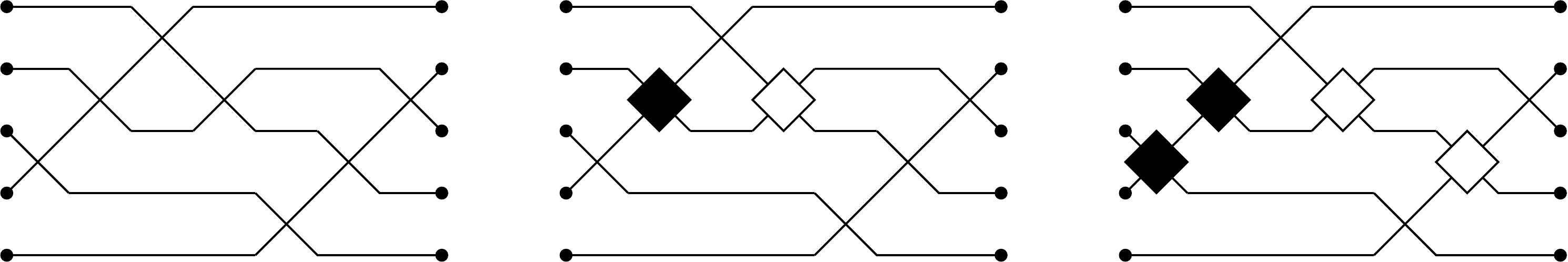}
\end{center}
\caption{Preancestries for the word $a_3a_2a_1a_2a_4a_3a_2$
with dimensions $0, 1, 2$.}
\label{fig:preancestry}
\end{figure}

There are usually several reduced words for a given permutation $\sigma$
but we shall keep our word fixed in the construction
of the stratification of $\BL_\sigma$.
A \textit{preancestry} for a reduced word
is obtained by marking certain crossings with black and white squares,
as in Figure \ref{fig:preancestry}.
The first example is the empty preancestry.
We always mark the same number $d$ of black and white squares:
$d$ is called the \textit{dimension} of the preancestry.
There are other conditions in the definition of a preancestry,
to be made precise in Section \ref{section:preancestry}
(see Equations \eqref{equation:preancestry0} and
\eqref{equation:preancestry}).

\begin{figure}[ht]
\begin{center}
\includegraphics[scale=0.25]{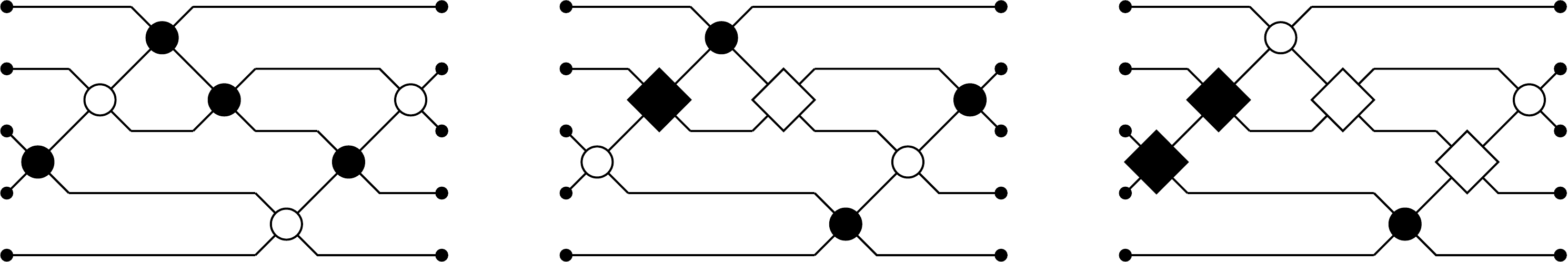}
\end{center}
\caption{Ancestries corresponding to the preancestries
in Figure \ref{fig:preancestry}.
These examples are
$(-1,+1,-1,-1,+1,-1,+1)$,
$(+1,-2,-1,+2,-1,+1,-1)$ and
$(-2,-2,+1,+2,-1,+2,+1)$. }
\label{fig:ancestry}
\end{figure}

An \textit{ancestry} is obtained from a preancestry
by marking the remaining crossings with black and white disks,
as in Figure \ref{fig:ancestry}
(see also Section \ref{section:ancestry}).
Given a preancestry of dimension $d$,
there are $2^{\ell - 2d}$ corresponding ancestries.
Figure \ref{fig:ancestry} shows examples of ancestries.
More formally, an ancestry is a sequence $\varepsilon$ of length $\ell$
assuming values in $\{\pm 1, \pm 2\}$:
square stands for $\pm 2$, disk stands for $\pm 1$,
black stands for negative and white for positive.

For each ancestry $\varepsilon$,
we define a smooth contractible submanifold
$\BLS_\varepsilon \subset \BL_\sigma$
of codimension $d = \dim(\varepsilon)$.
For a fixed reduced word, such submanifolds are disjoint
and their union is $\BL_\sigma$.
This decomposition is not as nice as might be desired.
In particular, except in the simplest cases,
it does not satisfy Whitney's condition:
Example \ref{example:bcbabdcb} exhibits
ancestries $\varepsilon_0$ and $\varepsilon_1$ such that
\[ \BLS_{\varepsilon_0} \cap \overline{\BLS_{\varepsilon_1}} \ne \emptyset,
\qquad
\BLS_{\varepsilon_0} \not\subseteq \overline{\BLS_{\varepsilon_1}}. \]
Nevertheless, a dual construction is possible.

\begin{theo}
\label{theo:CWcomplex}
For $\sigma \in S_{n+1}$, there exist a finite CW complex $\BLC_\sigma$
and a continuous map $i_\sigma: \BLC_\sigma \to \BL_\sigma$
with the following properties:
\begin{enumerate}
\item{The map $i_\sigma$ is a homotopy equivalence.}
\item{The cells $\BLC_\varepsilon$ of $\BLC_\sigma$
are labeled by ancestries $\varepsilon$.
For each ancestry $\varepsilon$ of dimension $d$,
the cell $\BLC_\varepsilon$ has dimension $d$.}
\end{enumerate}
\end{theo}

The desired CW complex $\BLC_\sigma$
and the homotopy equivalence are constructed in Section \ref{section:CW}.
The recursive step is provided by Lemma \ref{lemma:topolemma},
a topological result similar to many well known
results but for which we could not locate a published proof.
The CW complex  $\BLC_\sigma$ is, 
roughly speaking, a dual cell structure to the stratification
(see \cite{Hatcher}, Section 3.3,
particularly the figure in page 232).
The construction is reasonably explicit,
but the description of the glueing maps
is not as direct as might be desired.
For a small value of $n$ and an explicit $\sigma$,
the CW complex $\BLC_\sigma$ can be constructed by hand.
We shall construct $\BLC_\sigma$ explicitly for several examples
in the text and
this is how we prove the first two items of Theorem~\ref{theo:one}.

\begin{figure}[ht]
\begin{center}
\includegraphics[scale=0.25]{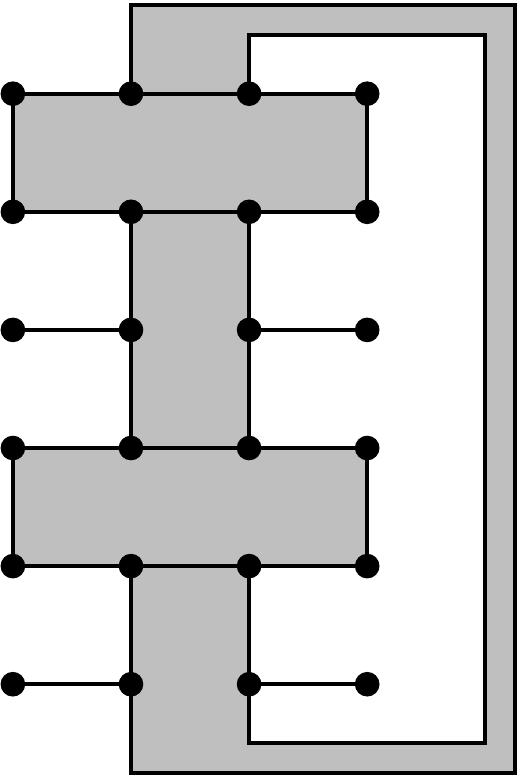}
\end{center}
\caption{A connected component of $\BLC_\sigma$
for $\sigma = 563412$;
see also Figure \ref{fig:563412}. }
\label{fig:563412CW}
\end{figure}

\begin{figure}[p]
\begin{center}
\includegraphics[scale=0.35]{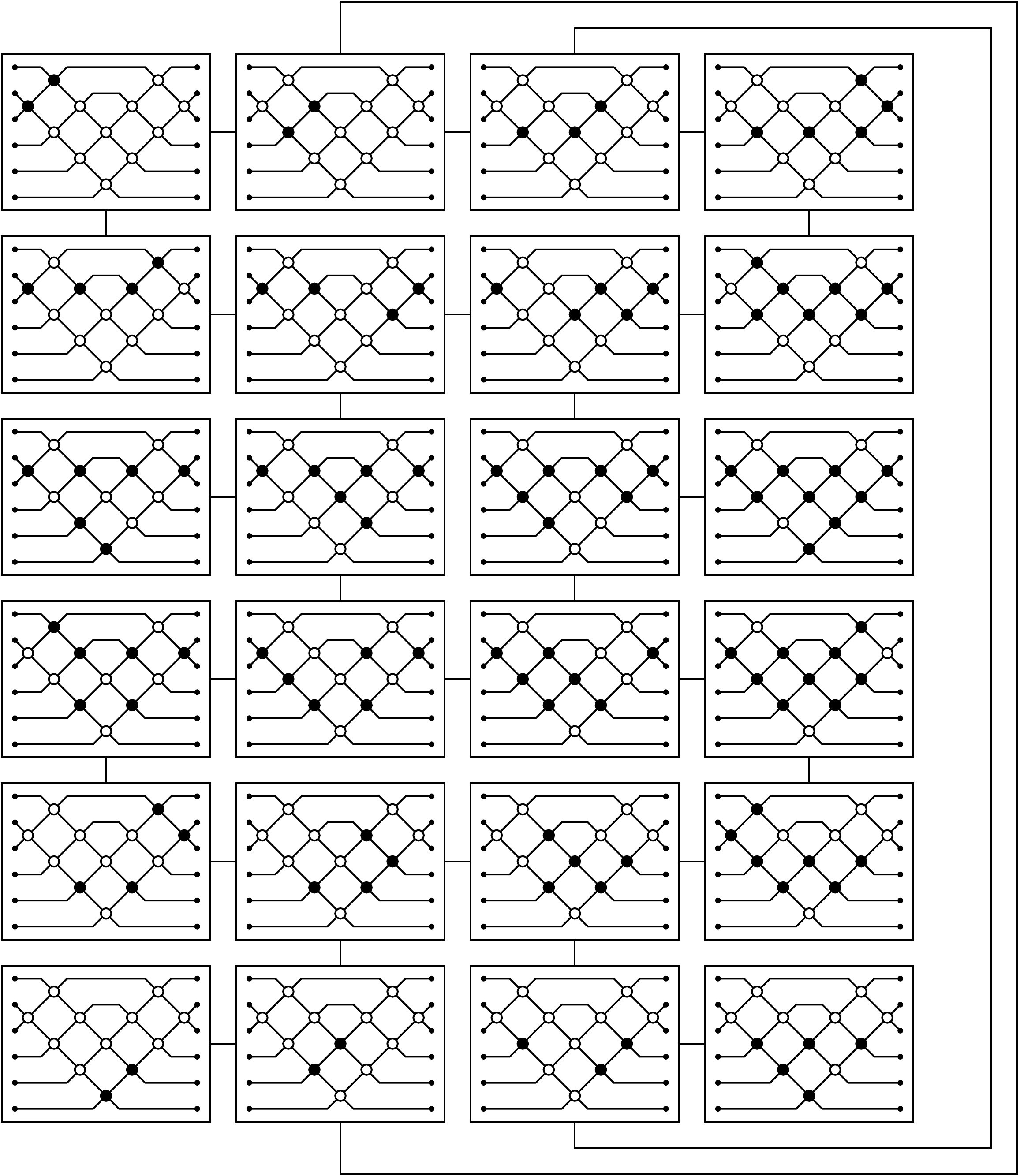}
\end{center}
\caption{A connected component of $\BLC_\sigma$
for $\sigma = 563412$.}
\label{fig:563412}
\end{figure}

Figures \ref{fig:563412CW} and \ref{fig:563412}, for instance,
show a connected component of $\BLC_\sigma$,
$\sigma = 563412 = a_2a_1a_3a_2a_4a_3a_5a_4a_2a_1a_3a_2 \in S_6$.
In Figure \ref{fig:563412},
cells of dimension $0$ are indicated by boxed ancestries
and cells of dimension $1$ by edges.
Cells of dimension $2$ are shown in Figure \ref{fig:563412CW}:
there are two octagons and two hexagons.
There are no cells of dimension $3$ or higher.
It should be clear from Figure \ref{fig:563412CW} that this is homotopically
equivalent to $\Ss^1$, as claimed in the second item of Theorem~\ref{theo:one}.

In order to prove the third item of Theorem \ref{theo:one}
we need a better understanding of the connected components of $\BL_\sigma$:
we need to know which subsets $\BLS_\varepsilon$
are contained in which connected component.
We shall discuss the connected components further in this paper.
For the moment let us introduce a partition
of $\BL_\sigma$ into $2^{n+1}$ subsets which are both open and closed
(but which may be empty or disconnected):
\begin{equation}
\label{equation:BLz}
\BL_\sigma = \bigsqcup_{z \in \acute\sigma \Quat_{n+1}} \BL_z.
\end{equation}
The definition of the sets $\BL_z$ is given in
Equation \eqref{equation:BLz2}, Section \ref{section:gsobis}.
Following the notation of \cite{Goulart-Saldanha0},
we have a finite group $\tilde B_{n+1}^{+} \subset \Spin_{n+1}$
and a surjective homomorphism $\Pi: \tilde B_{n+1}^{+} \to S_{n+1}$
with kernel $\Quat_{n+1}$,
a normal subgroup of order $2^{n+1}$.
For $\sigma \in S_{n+1}$, we have $\acute\sigma \in \tilde B_{n+1}^{+}$
with $\Pi(\acute\sigma) = \sigma$ so that
$\acute\sigma \Quat_{n+1} = \Pi^{-1}[\{\sigma\}] \subset  \tilde B_{n+1}^{+}$
is a coset (both left and right). 
We shall go over this in Sections \ref{section:tildeB} and \ref{section:gsobis}.

\begin{figure}[ht]
\begin{center}
\includegraphics[scale=0.25]{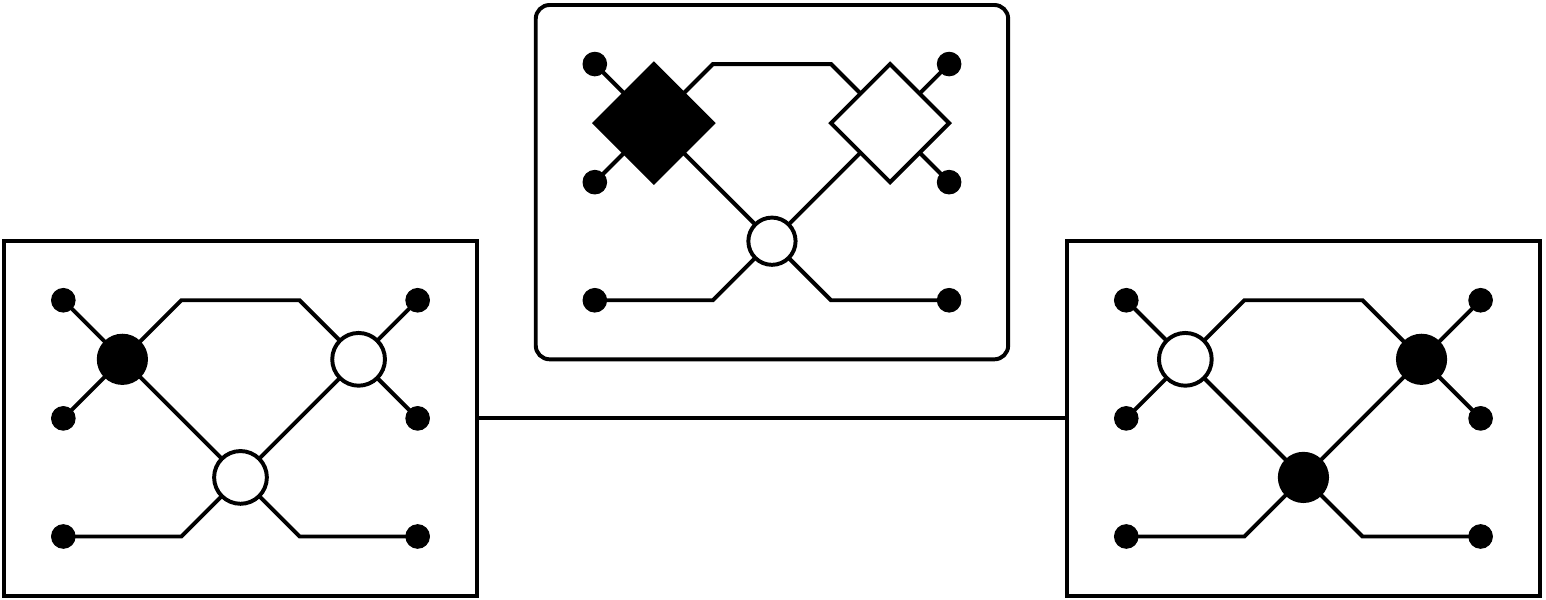}
\end{center}
\caption{The CW complex $\BLC_{z_0} \subset \BLC_{\sigma}$
is homotopically equivalent to $\BL_{z_0}$: both are contractible.
The box on the left and right represent the vertices
$\BLC_{(-1,+1,+1)}$ and $\BLC_{(+1,-1,-1)}$, respectively.
The only edge is $\BLC_{(-2,+1,+2)}$.}
\label{fig:aba-CW}
\end{figure}

\begin{example}
\label{example:abaintro}
Take $n = 2$ and $\sigma = \eta = a_1a_2a_1$.
The subset $\BL_{\eta} \subset \Lo_{3}^{1}$ is open
with $6$ connected components, all contractible:
\begin{equation}
\label{equation:abaintro}
\BL_{\eta} = \{L\;|\;z \ne 0, z \ne xy\}, \qquad
L = \begin{pmatrix}
1 & 0 & 0 \\ x & 1 & 0 \\ z & y & 1 
\end{pmatrix}.
\end{equation}
The connected component $\BL_{z_0} = \{L \;|\; z > \max\{0,xy\}\}$
is decomposed into three strata,
two of them open (corresponding to $y > 0$ and $y < 0$).
The third stratum is the half plane $y = 0$, $z > 0$.
Figure \ref{fig:aba-CW} shows $\BLC_{z_0} = i_{\eta}^{-1}[\BL_{z_0}]$, 
the connected component of $\BLC_{\eta}$ corresponding to $\BL_{z_0}$.
This is the simplest non trivial example
and shall be discussed in Examples
\ref{example:aba}, \ref{example:abacE}, 
\ref{example:abapos}, \ref{example:abazero} and \ref{example:abaL}. 
\end{example}

In order to prove the first item of Theorem \ref{theo:one},
we construct the CW complexes $\BLC_z$ for several $z \in \tilde B_{n+1}^{+}$,
as in the previous example,
and verify that they are, indeed, contractible.
It turns out that we have a somewhat stronger result.

\begin{theo}
\label{theo:collapse}
For $n \le 4$ and $z \in \tilde B_{n+1}^{+}$,
each connected component $X \subseteq \BLC_{z}$
collapses to a point: $X \searrow \{\bullet\}$.
\end{theo}

Theorem \ref{theo:collapse} implies the first item of Theorem \ref{theo:one}.

Given a preancestry $\varepsilon_0$ and $z \in \acute\sigma \Quat_{n+1}$,
our next main result gives a formula for
the number $N_{\varepsilon_0}(z)$ of ancestries $\varepsilon$
for which $\BLS_\varepsilon \subseteq \BL_{z}$.
Such a formula yields information about
the connected components of $\BL_\sigma$,
implying in particular the last item of Theorem \ref{theo:one}.
The formula is also very helpful as a check
when working out concrete examples, such as in Figure \ref{fig:563412}.

Before stating our next main result,
we need a few remarks about notation.
We interpret $\Spin_{n+1}$ to be contained
in the Clifford algebra $\Cliff_{n+1}^{0}$.
For $z \in \Cliff_{n+1}^{0}$,
let $\Re(z) = \langle 1,z\rangle \in \RR$
be its real part,
so that $2^n \Re: \Spin_{n+1} \to \RR$ is a character.
We shall see that if
$z \in \acute\sigma \Quat_{n+1} \subset \tilde B_{n+1}^{+}$
then $\Re(z)$ equals either $0$ or $\pm 2^{-(n+1-c)/2}$,
where $c = \nc(\sigma)$ is the number of cycles of $\sigma$.
There is always an element $z_0 \in \acute\sigma \Quat_{n+1}$
with $\Re(z_0) = 2^{-(n+1-c)/2} > 0$.
Also, given a preancestry $\varepsilon_0$ we define
a subgroup $H_{\varepsilon_0} \le \Quat_{n+1}^{+}$
in Equation \eqref{equation:Hvarepsilon}, Section \ref{section:ancestry}.

\begin{theo}
\label{theo:two}
Consider a permutation $\sigma \in S_{n+1}$,
a reduced word and a preancestry $\varepsilon_0$.
Let $z_0 \in \acute\sigma \Quat_{n+1}$
be such that $\Re(z_0) > 0$.

For any $z = qz_0 \in \acute\sigma\Quat_{n+1}$, we have:
\begin{align*}
N_{\varepsilon_0}(z) - N_{\varepsilon_0}(-z) &= 2^{\frac{\ell-2d}{2}} \Re(z);
\\
N_{\varepsilon_0}(z) + N_{\varepsilon_0}(-z) &= 
\begin{cases}
2^{\ell-2d+1}/|H_{\varepsilon_0}|, & q \in H_{\varepsilon_0}, \\
0, & q \notin H_{\varepsilon_0}.
\end{cases} 
\end{align*}
\end{theo}

Theorem \ref{theo:two} above can be used
(together with previous results by several authors,
detailed in Section \ref{section:history})
to give an effective enumeration of the $3 \cdot 2^n$
connected components of $\BL_\eta$ for $n \ge 5$ (see Proposition \ref{prop:effective}).
The theorem also implies that,
if $z \in \acute\eta\Quat_{n+1}$ and $\Re(z) < 0$
then the Euler characteristic $\chi(\BL_z)$ is even.
For $n \le 4$, $\BL_z$ has an even number
of contractible connected components.
For $n \ge 5$, however, the number of connected components
is odd and all except one are contractible:
this last connected component has therefore even Euler characteristic.

\bigskip

Section \ref{section:history} outlines the context
and history of this problem.
In Section \ref{section:permutations} we review some notation
and useful facts about permutations;
we are particularly interested in 
the \textit{Bruhat order} and in \textit{reduced words}.
In Section \ref{section:preancestry} we define \textit{preancestries},
one of the main combinatorial tools in our work.
In Section \ref{section:clifford} we present a quick construction
of Clifford algebras and prove a few results.
In Section \ref{section:tildeB} we construct the finite group
$\tilde B_{n+1}^{+} \subset \Spin_{n+1}$.
\textit{Ancestries} are defined based on preancestries
and on the group $\tilde B_{n+1}^{+}$
in Section \ref{section:ancestry}.
The concept of a \textit{thin} ancestry is discussed
in Section \ref{section:thin}:
this will tie in with the concept of total positivity.
Section \ref{section:gsobis} is largely a review of \cite{Goulart-Saldanha0}:
we discuss the Bruhat decomposition of $\Spin_{n+1}$
and define the \textit{lifted Bruhat order} in $\tilde B_{n+1}^{+}$.
This induces a partial order in the set of ancestries
(for a fixed reduced word).
We already mentioned in Equation~\eqref{equation:BLz}
that the set $\BL_{\sigma}$ is the disjoint union
of closed and open disjoint subsets $\BL_z$:
this is explained in Section \ref{section:BLz}.

The second part of the paper begins with Section \ref{section:firstexamples}
where we present first examples of the stratification of $\BL_\sigma$
(or of $\BL_z$) into subsets $\BLS_\varepsilon$
where $\varepsilon$ is an ancestry.
The stratification is presented in full generality
in Section \ref{section:stratification}:
this is where we complete the proof of Theorem \ref{theo:two}.
Basic properties of the strata $\BLS_\varepsilon$
are proved in Section \ref{section:goodstrata}.
These are smooth contractible submanifolds.
Whitney's property does not hold (Example \ref{example:bcbabdcb}).
Nevertheless, the partial order among ancestries
defined in Section \ref{section:gsobis}
(based on the lifted Bruhat order)
allows us to construct a dual CW complex and
prove Theorem \ref{theo:CWcomplex} is Section \ref{section:CW}.
Lemma \ref{lemma:topolemma} gives us the recursive step:
this is a topological result similar to many well known results,
such as Poincaré duality,
but we did not find a reference with the required generality.
In Section \ref{section:euler} we apply Theorem \ref{theo:CWcomplex}
to deduce information about the Euler characteristic
of $\BL_z$ and of its connected components:
we complete the proof of the third item of Theorem \ref{theo:one}
in Corollary \ref{coro:eveneuler} and Remark \ref{remark:eveneuler}.
As with any CW complex, the glueing maps for $\BLC_\sigma$
are extremely important:
they are studied in Section \ref{section:glue}.
We define the concept of a \textit{tame} ancestry:
ancestries in low dimensional examples are tame and
for tame ancestries, the glueing map is easy.
In Section \ref{section:moreexamples}
we apply our results to compute the homotopy type of $\BL_z$
in several examples;
we also complete the proof of Theorems \ref{theo:one} and \ref{theo:collapse}.
The proof of Theorem \ref{theo:collapse} is by then a long computation.
Some examples are given in
Examples \ref{example:4231}, \ref{example:4321},
\ref{example:54321} and \ref{example:54231}.
Other examples are discussed in greater detail
in \cite{Alves-Saldanha-2}.
In particular, the case $n = 4$ (corresponding to $S_5$)
and $\sigma = \eta$, the top permutation,
is a larger example than the ones discussed here,
with cells of dimension up to $4$.
This example is carefully described 
in \cite{Alves-Saldanha-2}.

\smallskip
\noindent \emph {Acknowledgements.}
The authors warmly thank Dan Burghelea,
Victor Goulart, Giovanna C.~Leal, Boris Shapiro, Michael Shapiro
and Cong Zhou for helpful conversations and the referee for a very careful report.
The first author wants to acknowledge the
hospitality of the Department of Mathematics, Stockholm University
in February 2019 when Boris Shapiro proposed the problem
which would eventually lead to this paper.
The first author was working at the Mathematics Department of PUC-Rio
during part of the time when this work was produced.
The second author gratefully acknowledges
the support of CNPq, CAPES and Faperj (Brazil).



\section{Context and history of the problem}
\label{section:history}

Bruhat/Schubert cell decompositions of Grassmannians and various spaces of
flags have been used in mathematics for more than a century and are standard
objects/tools in e.g. topology, enumerative geometry and representation theory.
Intersections of pairs and more general collections of Bruhat cells appear
naturally in several areas such as singularity theory,  Kazhdan-Lusztig
theory, matroid theory, to mention a few.
In spite of their importance, to the best of our
knowledge, there is hardly any topological information available about such
intersections, see e.g. \cite{SSV3} and references therein. 

One  exception from this general situation is
the problem of counting connected
components in pairwise intersections of big (i.e., top-dimensional)
Bruhat cells over the reals.
Substantial progress was obtained in the late 90's,
see \cite {SSV1, SSV2, SSVZ, Ri1, Ri2, GeShVa, Z}. 
In short,  this problem can be reduced to  counting the orbits 
of a certain finite group of symplectic transvections 
acting on a  finite-dimensional vector space over the
finite field $\FF_2$ (also denoted by $\ZZ/(2)$). 
Both the group and the vector space are uniquely
determined by   the pair of Bruhat cells under consideration, see \cite{SSVZ}.
Further information about counting such orbits can be found in \cite{Se}.

Consider the real flag space
$\Flag_{n+1}=\SL_{{n+1}}/\Up_{n+1} = \Spin_{n+1}/\Quat_{n+1}$:
notice that $\Up_{n+1}$ is the  Borel subgroup
(usually called $B$).
The intersection of two opposite big Bruhat cells in $\Flag_{n+1}$
is homeomorphic to $\BL_\eta$
(as in Equation \eqref{equation:BL}).
The number of connected components of $\BL_\eta$ equals
$2, 6, 20, 52$ for $n=1,2,3,4$ respectively.
Starting from $n=5$, the number of connected components
stabilizes and is given by $3\cdot 2^{n}$.
This is explained by the possibility to embed, for $n\ge 5$,
the lattice $E_6$  in a certain lattice arising in this problem,
see \cite{SSV2}.  

Observe that the relative positions of two big Bruhat cells
in $\Flag_{n+1}$ are in bijective correspondence with $S_{n+1}$.
In particular, 
opposite big Bruhat cells correspond to the longest permutation $\eta$.
The study of  the number of connected components
in the intersection of two big cells in a given relative
position $\sigma$ was initiated in \S~7 of \cite{SSV2}. 
For each concrete $\sigma$,
the number can, in principle,
be deduced from the results of \cite{Se} obtained about two decades ago.
However,  to the best of our knowledge, 
there is no closed formula.

The main goal of the present paper
is to introduce the appropriate tools
which allow us to study the homotopy type the latter  intersections. 
As a corollary, we obtain in Proposition \ref{prop:effective}
an effective labeling
of the connected components of $\BL_\eta$.
We make extensive use of notation and results
from \cite{Goulart-Saldanha0},
where the Bruhat cell decomposition of $\Spin_{n+1}$ is studied:
we review the needed facts in Section \ref{section:gsobis}.
We use Clifford algebras.
In order to keep this paper more self-contained
we review the basic constructions in Section \ref{section:clifford};
see \cite{Atiyah-Bott-Shapiro, Lawson-Michelsohn} 
for far more complete discussions.
We also need some basic facts from the theory of
totally positive matrices;
\cite{Berenstein-Fomin-Zelevinsky} is an excellent reference
and contains much more than we need.

The stratification we construct appears to be
different but closely related to the one constructed by V.~V.~Deodhar.
As of this writing, the precise relationship
between the two stratifications is unfortunately not clear to us.

\section{Permutations}
\label{section:permutations}

Recall that
if $\sigma \in S_{n+1}$ is a permutation,
we write $i^\sigma$ (and not $\sigma(i)$)
and we compose left to right so that
$i^{\sigma_0\sigma_1} = (i^{\sigma_0})^{\sigma_1}$;
we use several notations for permutations $\sigma \in S_{n+1}$.
The \textit{complete} notation just lists the values
of $1^\sigma, 2^\sigma, \ldots, (n+1)^\sigma$:
we usually enclose the list in square brackets
but sometimes omit them for brevity.
Another common notation is to write $\sigma$
as a product of disjoint cycles:
Figure \ref{fig:reduced}, for instance,
shows $\sigma = [43512] = (14)(235)$.
Perhaps the most important notation for us
is to write $\sigma$ as a \textit{reduced word},
a word of minimal length $\ell = \inv(\sigma)$
in the standard Coxeter-Weyl  generators
$a_1, \ldots, a_n$ (see Equation \eqref{equation:reduced}).
Here $a_i$ is the transposition $(i,i+1)$.
Figure \ref{fig:reduced} shows a \textit{wiring diagram},
a visual representation of a reduced word.
In that example, $\sigma = a_1a_3a_2a_1a_4a_3a_2$;
in our system,
each crossing is a generator, read left-to-right,
and the row indicates the value of $i$.
The identity permutation is denoted by $1 \in S_{n+1}$.

Given $\sigma \in S_{n+1}$,
an \textit{inversion} is a pair $(i,j) \in \nmaisum^2$
with $i < j$ and $i^\sigma > j^\sigma$.
Let $\Inv(\sigma)$ be the set of inversions of $\sigma$;
let $\inv(\sigma) = |\Inv(\sigma)|$.
Crossings in a wiring diagram for a reduced word for $\sigma$
naturally correspond to inversions of $\sigma$.
Indeed, for $i \in \nmaisum$,
the \textit{wire} $i$ is the curve
joining the point $i$ on the left to the point $i^\sigma$ on the right.
Given an inversion $(i,j)$, the wires $i$ and $j$ cross exactly once:
that is the crossing $(i,j)$;
sometimes we omit parenthesis and comma for brevity.
For the first reduced word in Figure \ref{fig:reduced},
the inversions are, in order,
$12$, $34$, $14$, $24$, $35$, $15$, $25$.
For the second reduced word the inversions
are the same but the order is different:
$34$, $24$, $14$, $12$, $35$, $15$, $25$.
More algebraically, if $\ell = \inv(\sigma)$
and $k$ is the crossing $(i,j)$ then
\begin{equation}
\label{equation:crossing}
\begin{aligned}
\sigma =
a_{i_1}\cdots a_{i_k} \cdots a_{i_\ell}
&= (i j)\;a_{i_1}\cdots a_{i_{k-1}} a_{i_{k+1}} \cdots a_{i_\ell} \\
&= a_{i_1}\cdots a_{i_{k-1}} a_{i_{k+1}} \cdots a_{i_\ell}\;
(j^\sigma i^\sigma).
\end{aligned}
\end{equation}
Notice that the subword
$a_{i_1}\cdots a_{i_{k-1}} a_{i_{k+1}} \cdots a_{i_\ell}$
is usually not reduced.

We use the strong Bruhat order throughout the paper,
and write $\sigma_0 \vartriangleleft \sigma_1$
if $\sigma_0 < \sigma_1$ (in the strong Bruhat order)
and $\inv(\sigma_1) = 1+\inv(\sigma_0)$.
Equivalently, $\sigma_0 \vartriangleleft \sigma_1$
if there exists a reduced word 
$\sigma_1 = a_{i_1} \cdots a_{i_\ell}$ such that
$\sigma_0 = a_{i_1}\cdots a_{i_{k-1}} a_{i_{k+1}} \cdots a_{i_\ell}$
is also a reduced word.

A permutation $\sigma \in S_{n+1}$ {\em blocks at $j$}
(where $1 \le j \le n$)
if and only if $i \le j$ implies $i^\sigma \le j$.
Equivalently, $\sigma$ blocks at $j$
if and only if $a_j$ does not appear in a reduced word for $\sigma$.
Let $\Block(\sigma)$ be the set of $j$ such that
$\sigma$ blocks at $j$ and $b = \block(\sigma) = |\Block(\sigma)|$.
A permutation $\sigma$ \textit{does not block} if $\block(\sigma) = 0$.

Given $\sigma_0 \in S_j$ and $\sigma_1 \in S_{k}$
define
\[
\sigma = \sigma_0 \oplus \sigma_1 \in S_{j+k}, \qquad
\sigma^i =
\begin{cases} \sigma_0^i, & i \le j, \\ \sigma_1^{i-j}+j, & i > j. 
\end{cases} \]
If $\sigma \in S_{n+1}$ blocks at $j$ then there exist permutations
$\sigma_0 \in S_j$ and $\sigma_1 \in S_{n+1-j}$
such that $\sigma = \sigma_0 \oplus \sigma_1$.
In this case, the permutation matrices satisfy
$P_\sigma = P_{\sigma_0} \oplus P_{\sigma_1}$,
i.e., the matrix $P_\sigma$
has two diagonal blocks $P_{\sigma_0}$ and $P_{\sigma_1}$
(and is zero elsewhere).

\begin{remark}
\label{remark:blockfree}
If $\sigma = \sigma_0 \oplus \sigma_1$ then
$\BL_{\sigma} = \BL_{\sigma_0} \oplus \BL_{\sigma_1}$,
meaning that $L \in \BL_{\sigma}$ if and only if 
there exist $L_0 \in \BL_{\sigma_0}$ and $L_1 \in \BL_{\sigma_1}$
with $L = L_0 \oplus L_1$.
As a manifold, we have that 
$\BL_{\sigma}$ is diffeomorphic to $\BL_{\sigma_0} \times \BL_{\sigma_1}$.
Since our aim is to study the homotopy type of the sets $\BL_\sigma$,
we may focus on permutations $\sigma$ which do not block.
\end{remark}


\section{Preancestries}
\label{section:preancestry}

A \textit{preancestry} for a reduced word $\sigma = a_{i_1}\cdots a_{i_\ell}$
is a sequence $(\rho_k)_{0 \le k \le \ell}$ of permutations
with the following properties:
\begin{enumerate}
\item{$\rho_0 = \rho_\ell = \eta$;}
\item{for all $k \in \llbracket \ell \rrbracket$,
we either have
$\rho_k = \rho_{k-1}$ or $\rho_k = \rho_{k-1} a_{i_k}$;}
\item{for all $k \in \llbracket \ell \rrbracket$,
if $\rho_{k-1} a_{i_k} > \rho_{k-1}$
then
$\rho_k = \rho_{k-1} a_{i_k}$.}
\end{enumerate}

\goodbreak

It is usually more convenient to represent a preancestry $(\rho_k)$
by the sequence
$\varepsilon_0: \llbracket \ell \rrbracket \to \{0,\pm 2\}$
with
\begin{equation}
\label{equation:preancestry0}
\varepsilon_0(k) = 
\begin{cases}
0, & \rho_k = \rho_{k-1}, \\
-2, & \rho_k = \rho_{k-1} a_{i_k} < \rho_{k-1}, \\
+2, & \rho_k = \rho_{k-1} a_{i_k} > \rho_{k-1}. 
\end{cases}
\end{equation}
Thus, a sequence $\varepsilon_0$ is a \textit{preancestry}
if the sequence of permutations
$(\rho_k)_{0 \le k \le \ell}$
defined by the recursion below is a preancestry: 
\begin{equation}
\label{equation:preancestry}
\rho_0 = \eta, \quad
\rho_k = 
\begin{cases}
\rho_{k-1} a_{i_k}, & \varepsilon_0(k) \ne 0, \\
\rho_{k-1}, &\varepsilon_0(k) = 0.
\end{cases}
\end{equation}
The preancestry $\varepsilon_0$ is represented
in the wiring diagram for $\sigma$:
a black (resp. white) square indicates $-2$ (resp. $+2$).

\begin{figure}[ht]
\begin{center}
\includegraphics[scale=0.15]{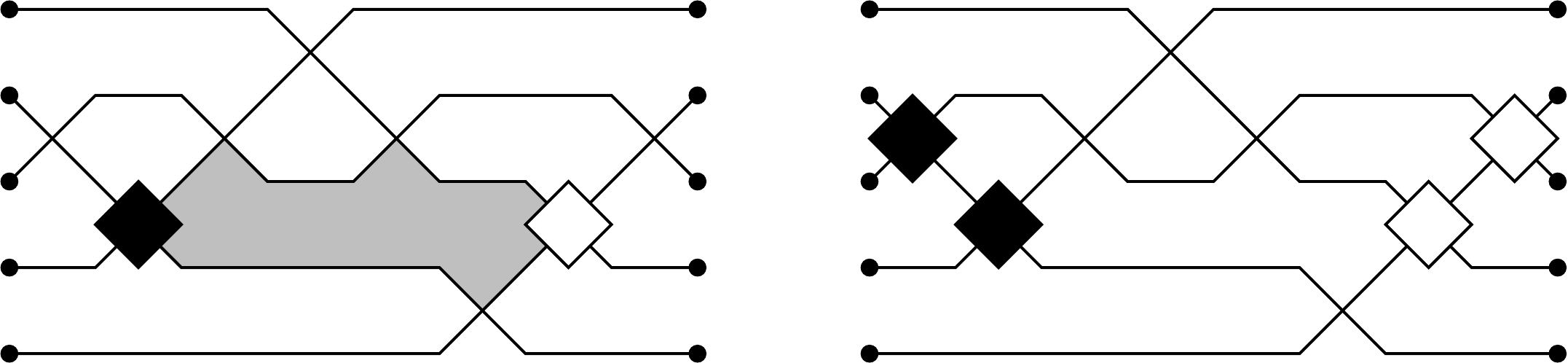}
\end{center}
\caption{Two preancestries
for $\sigma = a_2a_3a_2a_1a_2a_4a_3a_2 \in S_5$.
The shaded region corresponds to the preancestry of dimension $1$.
The second preancestry has dimension $2$.}
\label{fig:bcbabdcb}
\end{figure}

\begin{figure}[ht]
\begin{center}
\includegraphics[scale=0.15]{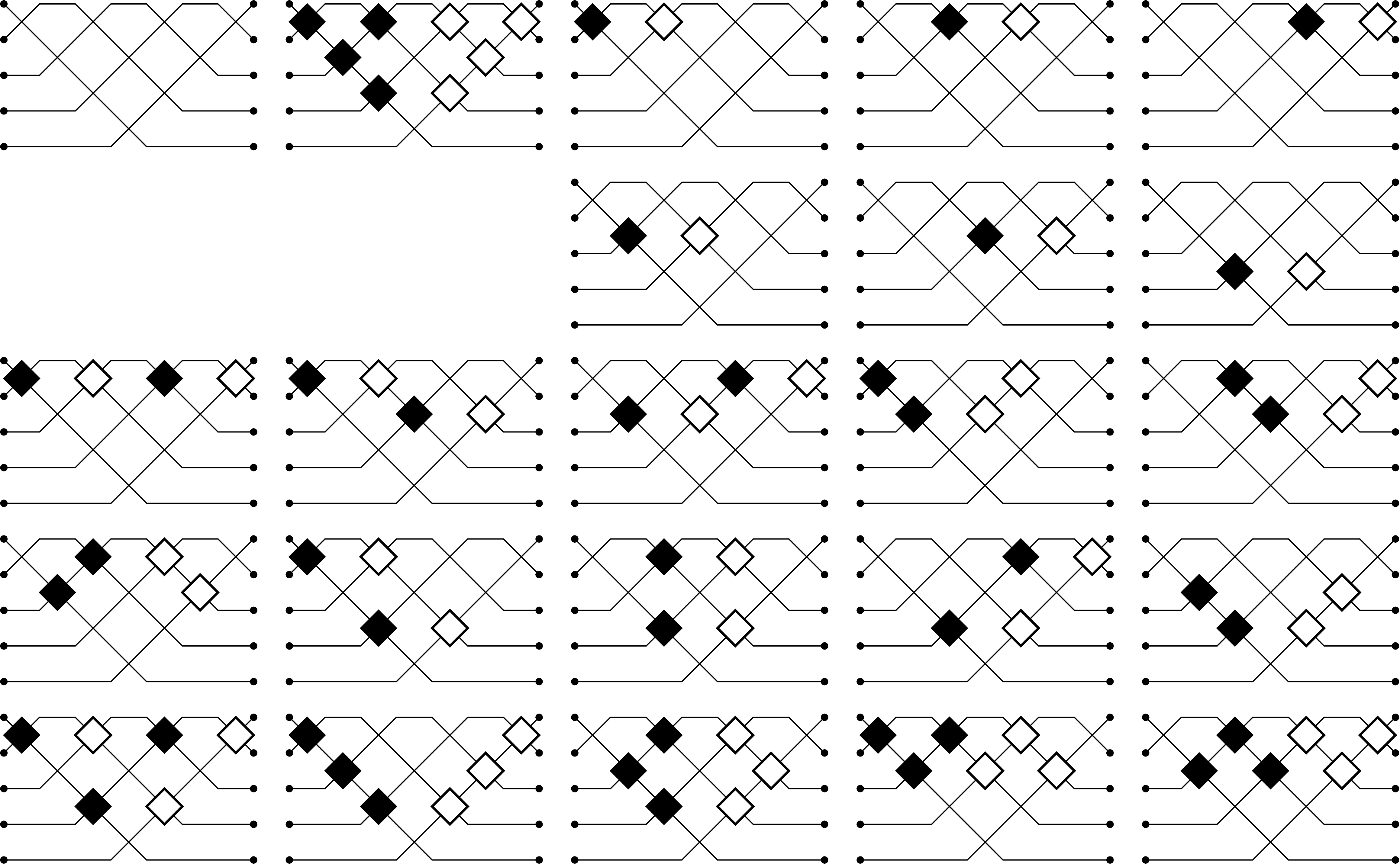}
\end{center}
\caption{The reduced word $\eta = a_1a_2a_1a_3a_2a_1a_4a_3a_2a_1 \in S_5$
admits respectively $1, 6, 10, 5, 1$ preancestries of dimensions
$0, 1, 2, 3, 4$.}
\label{fig:preancestry54321}
\end{figure}

Clearly, in any preancestry $\varepsilon_0$,
the number of $k \in \llbracket\ell\rrbracket$
for which $\varepsilon_0(k) = -2$ equals
the number of $k \in \llbracket\ell\rrbracket$
for which $\varepsilon_0(k) = +2$:
this number is denoted by $d = \dim(\varepsilon_0)$,
the \textit{dimension} of the preancestry.
Figure \ref{fig:bcbabdcb} shows two preancestries 
for $\sigma = a_2a_3a_2a_1a_2a_4a_3a_2 \in S_5$,
of dimensions $1$ and $2$.
Figure \ref{fig:preancestry54321} shows all preancestries for
$\eta = a_1a_2a_1a_3a_2a_1a_4a_3a_2a_1 \in S_5$.

\begin{lemma}
\label{lemma:countpreancestry}
Consider a permutation $\sigma \in S_{n+1}$.
The number of preancestries per dimension 
does not depend on the choice of the reduced word.
\end{lemma}

\begin{proof}
We know that two reduced words for $\sigma$
can be joined by a finite number of steps 
of the two following moves:
\begin{gather*}
a_{i_1}\cdots a_{i_k}a_{i_{k+1}} \cdots a_{i_\ell} =
a_{i_1}\cdots a_{i_{k+1}}a_{i_{k}} \cdots a_{i_\ell}, \quad
|i_k - i_{k+1}| \ne 1; \\
a_{i_1}\cdots a_{i_k}a_{i_{k+1}} a_{i_k} \cdots a_{i_\ell} =
a_{i_1}\cdots a_{i_{k+1}}a_{i_{k}} a_{i_{k+1}} \cdots a_{i_\ell}, \quad
|i_k - i_{k+1}| = 1.
\end{gather*}
The effect of the first move is trivial.
In fact, in Figures \ref{fig:563412} and \ref{fig:preancestry54321}
diagrams have been drawn with several disjoint crossings
on the same vertical line,
thus blurring the distinction between reduced words
differing by the first move.

We must therefore prove that if two reduced words differ by the second move then the number of preancentries of each dimension is the same.
This can be verified by a finite computation, considering which of the three intersections involved in the move are marked by a preancestry.
We omit the rather long case-by-case verification. 
Notes that there are three wires involved in the second move: they do not cross to the left or right of the move.  
\end{proof} 

\begin{remark}
\label{remark:lowdimpreancestry}
There is always a unique preancestry of dimension $0$,
obtained by marking no vertices.
Preancestries of dimension $1$ are also easy to classify:
we must mark two consecutive intersections on the same row.
In other words, a preancestry of dimension $1$
corresponds to a bounded component of the complement
of the wiring diagram,
and the two marked intersections are the left and right extremes.
The component is shown shaded in Figure \ref{fig:bcbabdcb}.
The number of preancestries of dimension $1$
is $\ell - n + b$
where $\ell = \inv(\sigma)$ and $b = \block(\sigma)$.

For dimension $2$ the situation is slightly more complicated.
For a preancestry $\varepsilon_0$ of dimension $2$,
let $k_1 < k_2 < k_3 < k_4$ be such that $|\varepsilon_0(k_i)| = 2$.
We always have $\varepsilon_0(k_1) = -2$ and $\varepsilon_0(k_4) = +2$.
If $\varepsilon_0(k_2) = +2$, we have $\varepsilon_0(k_3) = -2$,
$i_{k_1} = i_{k_2}$, $i_{k_3} = i_{k_4}$.
In this case, intersections $k_1$ and $k_2$ are consecutive on row $i_{k_1}$
and intersections $k_3$ and $k_4$ are consecutive on row $i_{k_3}$.
The first three preancestries on the third row of
Figure \ref{fig:preancestry54321} are examples.
If $\varepsilon_0(k_2) = -2$ and $|i_{k_1} - i_{k_2}| > 1$
we also have two pairs of consecutive intersections on two rows. 
The preancestry on the fourth row, third column of
Figure \ref{fig:preancestry54321} is an example.
In both cases, the preancestry $\varepsilon_0$ is said to be of type I.

If $\varepsilon_0(k_2) = -2$
(which implies $\varepsilon_0(k_3) = +2$) and $|i_{k_1} - i_{k_2}| = 1$
then $\varepsilon_0$ is said to be of type II.
In this case, $i_{k_1} = i_{k_4}$, $i_{k_2} = i_{k_3}$
and intersections $k_2$ and $k_3$ are consecutive on row $i_{k_2}$.
Intersection $k_1$ is the last on row $i_{k_1}$ before $k_2$
and intersection $k_4$ is the first on row $i_{k_1}$ after $k_3$.
There is no bound, however, on the number of intersections 
on row $i_{k_1}$ between $k_2$ and $k_3$.
Figure \ref{fig:bcbabdcb} shows a preancestry of type II
for which there are two such intersections.
\end{remark}

The subword formed by all marked letters has value $1$.
The subword formed by unmarked letters is also informative,
as the following lemma shows.
It also gives an estimate
for the possible dimensions of preancestries.

\begin{lemma}
\label{lemma:transpositions}
Consider $\sigma \in S_{n+1}$ and a fixed reduced word
of length $\ell = \inv(\sigma)$.
Consider a preancestry $\varepsilon_0$ of dimension $d = \dim(\varepsilon_0)$.
There are $\delta = \ell - 2d$ unmarked crossings $k_1, \ldots, k_\delta$.
Assume that the unmarked crossing $k_j$ is
$(\iota_{j,0},\iota_{j,1}) \in \Inv(\sigma)$.
We then have
\[ \sigma =
(\iota_{\delta,0} \iota_{\delta,1}) \cdots (\iota_{1,0} \iota_{1,1}). \]
If $c = \nc(\sigma)$ is the number of cycles then
$2d \le l+c-n-1$.
\end{lemma}

\begin{proof}
Join marked crossings into subwords $\sigma_j$ to write
\[ \sigma = \sigma_0 a_{i_{k_1}} \sigma_1 \cdots
\sigma_{\delta-1} a_{i_{k_\delta}} \sigma_\delta, \qquad
1 = \sigma_0 \sigma_1 \cdots
\sigma_{\delta-1} \sigma_\delta. \]
Apply Equation \eqref{equation:crossing} to the crossing $k_\delta$
to obtain
\[ \sigma = (\iota_{\delta,0} \iota_{\delta,1})
\sigma_0 a_{i_{k_1}} \sigma_1 \cdots
\sigma_{\delta-1} \sigma_\delta. \]
Repeat to obtain
\[ \sigma = (\iota_{\delta,0} \iota_{\delta,1}) \cdots 
(\iota_{1,0} \iota_{1,1})
\sigma_0 \sigma_1 \cdots \sigma_{\delta-1} \sigma_\delta
= (\iota_{\delta,0} \iota_{\delta,1}) \cdots 
(\iota_{1,0} \iota_{1,1}), \]
proving the first claim.
The second claim follows from observing 
that any cycle of length $k$ can be written 
as a product of $k-1$ transpositions,
but not fewer.
\end{proof}

\begin{figure}[ht]
\begin{center}
\includegraphics[scale=0.2]{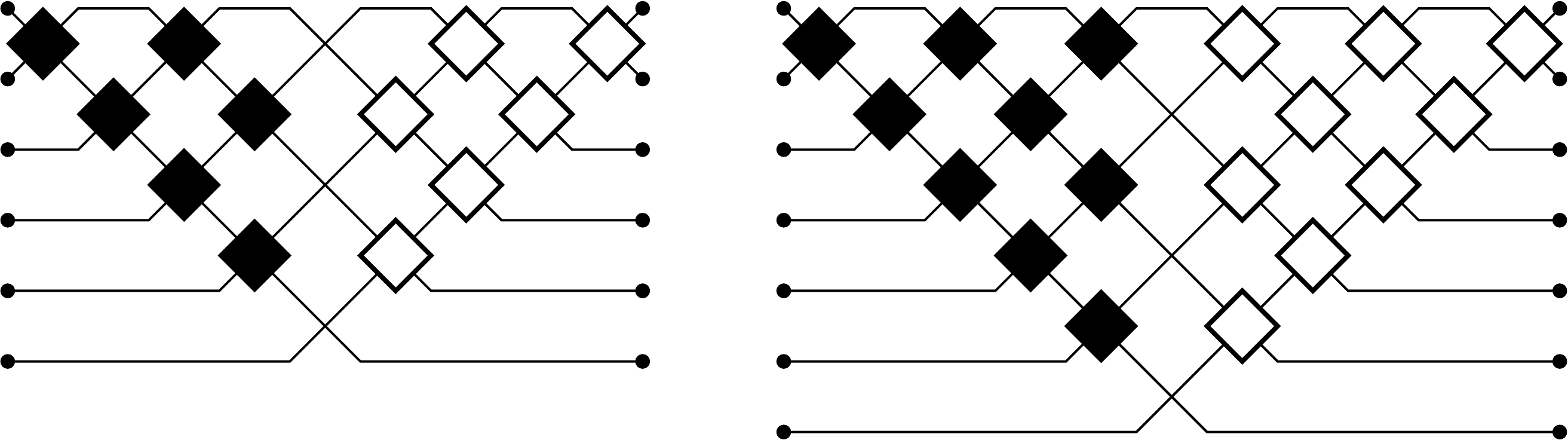}
\end{center}
\caption{Preancestries of top dimension $d_{\max}$ for $\eta \in S_{n+1}$, $n=5$ and $n=6$.}
\label{fig:preancestryeta}
\end{figure}

\begin{lemma}
\label{lemma:eta}
For $n \ge 2$, let $\eta \in S_{n+1}$ be the top permutation.
The largest possible dimension among all preancestries is
\[ d_{\max} = \left\lfloor \frac{n^2}{4} \right\rfloor. \]
Furthermore, there exists a unique preancestry of dimension $d_{\max}$.
\end{lemma}

\begin{proof}
From Lemma \ref{lemma:countpreancestry} we may use the reduced word
\[ \eta = a_1a_2a_1a_3a_2a_1\cdots a_na_{n-1} \cdots a_2a_1. \]
The case $n = 4$ is addressed by Figure \ref{fig:preancestry54321}.
Figure \ref{fig:preancestryeta} shows a preancestry of dimension $d_{\max}$
for cases $n = 5$ and $n = 6$. 
The pattern should be clear.
The value of $d_{\max}$ follows from these
examples and from Lemma \ref{lemma:transpositions}.
Since $\eta$ can be written as product of
$\lceil n/2 \rceil$ transpositions in a unique way,
the result follows.
\end{proof}





\section{Clifford algebras}
\label{section:clifford}

For each $k \in \NN^\ast$,
let $J_k$ be the $2^k \times 2^k$ real matrix defined by:
\[ J_1 = \begin{pmatrix} 0 & -1 \\ 1 & 0 \end{pmatrix};
\qquad
J_k = \begin{pmatrix} 
0 & 0 & -I & 0 \\
0 & 0 & 0 & I \\
I & 0 & 0 & 0 \\
0 & -I & 0 & 0 \end{pmatrix}, \;
I \in \RR^{2^{k-2} \times 2^{k-2}}, \;k > 1. \]
For $k \in \nmesmo = \{1,2,\ldots,n\}$, let
$\hat a_k$ be the real $2^n \times 2^n$ matrix
with $2^{n-k}$ blocks $J_k$ along the diagonal:
\[ \hat a_k = \begin{pmatrix}
J_k & & \\ & \ddots & \\ & & J_k \end{pmatrix}. \]
The following relations are easily verified:
\[ (\hat a_i)^2 = -I, \quad
\hat a_j \hat a_i = (-1)^{[|i-j| = 1]} \hat a_i \hat a_j. \]
We use Iverson's bracket notation:
thus, for instance,
$[|i-j| = 1] = 1$ if $|i-j| = 1$ and
$[|i-j| = 1] = 0$ otherwise.
Let $\Cliff_{n+1}^{0}$ be the associative algebra
generated by $\hat a_1, \ldots, \hat a_n$.
In the notation of \cite{Atiyah-Bott-Shapiro, Lawson-Michelsohn},
this is the subalgebra $\Cliff_{n+1}^{0} \subset \Cliff_{n+1}$
of even elements of the Clifford algebra $\Cliff_{n+1}$;
we shall abuse notation and refer to  $\Cliff_{n+1}^{0}$
as the Clifford algebra. 
The set $\HQuat_{n+1}$ below is a real basis of $\Cliff_{n+1}^{0}$:
\[ \HQuat_{n+1} = \{ [\hat a_1] [\hat a_2] \cdots [\hat a_n] \} 
=
\{1, \hat a_1, \hat a_2, \hat a_1 \hat a_2,
\hat a_3, \hat a_1\hat a_3,
\dots, \hat a_1\hat a_2\cdots \hat a_n \}; \]
here the square brackets
denote that the presence of the corresponding element
in the product is optional.
Notice that from now on we write $1$ for $I \in \Cliff_{n+1}^{0}$.


The set $\Quat_{n+1} = \HQuat_{n+1} \sqcup (-\HQuat_{n+1})
\subset \Cliff_{n+1}^{0}$
is a group
of cardinality $2^{n+1}$ generated by $\hat a_i$, $i \in \nmesmo$. 
The subgroup $\{ \pm 1\} \vartriangleleft \Quat_{n+1}$
is normal and contained in the center:
we have the exact sequence
\begin{equation}
\label{equation:QuatDiag}
1 \to \{ \pm 1 \} \to \Quat_{n+1} \to \Diag_{n+1}^{+} \to 1;
\end{equation}
the quotient $\Diag^{+}_{n+1}$ is naturally identified
with the subgroup of diagonal matrices in $\SO_{n+1}$.
Group-theoretic commutators  
in $\Quat_{n+1}$ will be occasionally important:
set 
\begin{equation}
\label{equation:commutator}
[q_0,q_1] = q_0^{-1}q_1^{-1}q_0q_1 \in \{\pm 1\} =
[\Quat_{n+1},\Quat_{n+1}] \vartriangleleft \Quat_{n+1}. 
\end{equation}
For instance,
if $q = \pm \hat a_1^{\varepsilon(1)} \cdots \hat a_n^{\varepsilon(n)}$
(with $\varepsilon(k) \in \ZZ$) then
$[\hat a_i,q] = (-1)^{\varepsilon(i-1)+\varepsilon(i+1)}$
(we follow the convention that 
$\varepsilon(0) = \varepsilon(n+1) = 0$).

Define in $\Cliff_{n+1}^{0}$ a positive definite inner product
so that $\HQuat_{n+1}$ is an orthonormal basis:
$\langle z_1, z_2 \rangle = 2^{-n} \Trace(z_1z_2^{\top})$.
Define the \textit{real part} of $z \in \Cliff_{n+1}^{0}$
by $\Re(z) = 2^{-n} \Trace(z) = \langle z,1 \rangle$:
equivalently, $\Re(z) = c_1$ if $z = \sum_{q \in \HQuat} c_q q$.
From now on we almost never consider
the elements of $\Cliff_{n+1}^{0}$ as matrices,
rather as elements of an explicitly given associative algebra.

Define $\fa_i^{\Spin} = \frac12 \hat a_i \in \Cliff_{n+1}^{0}$
and $\fa_i^{\SO}$ to be the $(n+1) \times (n+1)$
real skew symmetric matrix whose only non zero entries
are $(\fa_i^{\SO})_{i+1,i} = 1$ and $(\fa_i^{\SO})_{i,i+1} = -1$.
Let $\spin_{n+1} \subset \Cliff_{n+1}^{0}$ be the Lie algebra generated by 
$\fa_i^{\Spin}$, $i \in \nmesmo$.
It can easily be verified that there exists
a unique isomorphism of Lie algebras
$\Pi: \spin_{n+1} \to \so_{n+1}$ satisfying $\Pi(\fa_i^{\Spin}) = \fa_i^{\SO}$.
From now on we often identify the two isomorphic Lie algebras 
$\spin_{n+1}$ and $\so_{n+1}$.
We write $\alpha_i(\theta) = \exp(\theta \fa_i)$,
defining one parameter subgroups
$\alpha_i^{\Spin}: \RR \to \Spin_{n+1}$ and
$\alpha_i^{\SO}: \RR \to \SO_{n+1}$;
we often omit the superscripts.
More explicitly, we have
\[
\alpha_i^{\SO}(\theta) =
\begin{pmatrix} I & & & \\
& \cos(\theta) & -\sin(\theta) & \\
& \sin(\theta) & \cos(\theta) & \\
& & & I \end{pmatrix}, 
\quad
\alpha_i^{\Spin}(\theta) = \cos\left(\frac{\theta}{2}\right) +
\sin\left(\frac{\theta}{2}\right) \hat a_i; \]
in the first formula, the central block ocupies rows and columns
$i$ and $i+1$.
Let $\Spin_{n+1} \subset \Cliff_{n+1}^{0}$
be the group generated by $\alpha_i^{\Spin}(\theta)$
for $i \in \nmesmo$ and $\theta \in \RR$.
Taking exponentials, we have the familiar double cover
$\Pi: \Spin_{n+1} \to \SO_{n+1}$.
The function $\Re: \Spin_{n+1} \to \RR$ is
$2^{-n}$ times a character,
the Clifford algebra being the corresponding representation.

Let $\Diag_{n+1}$ be the group of diagonal matrices
with diagonal entries in $\{\pm 1\}$:
this group acts by conjugations on $\SO_{n+1}$. 
The quotient $\cE_n = \Diag_{n+1}/\{\pm I\}$
is naturally isomorphic to $\{\pm 1\}^{\nmesmo}$:
the matrix $D \in \Diag_{n+1}$ is taken to $E \in \cE_n =  \{\pm 1\}^{\nmesmo}$,
$E_i = D_{i,i} D_{i+1,i+1}$.
The group $\cE_n$ also acts by automorphisms on $\SO_{n+1}$.

This action can be lifted to $\Spin_{n+1}$ and 
then extended to $\Cliff_{n+1}^{0}$.
Indeed, each $E \in \cE_n$ defines automorphisms of
$\Spin_{n+1}$ and $\Cliff_{n+1}^{0}$ given by 
\begin{equation}
\label{equation:cE}
(\hat a_i)^E = E_i \hat a_i, 
\qquad
(\alpha_i(\theta))^E = \alpha_i(E_i \theta).
\end{equation}

\begin{remark}
\label{remark:cE}
If $n$ is even, conjugation by an element of $\Quat_{n+1}$
defines the same set of automorphisms as the action above.
If $n$ is odd, however, 
the element $\hat a_1\hat a_3 \cdots \hat a_{n-2} \hat a_n$
is in the center of $\Cliff_{n+1}^{0}$.
In this case, conjugation accounts for only half
of the automorphisms above.
\end{remark}

The following result gives a formula for $\Re(z)$, $z \in \Spin_{n+1}$,
in terms of the eigenvalues of $Q = \Pi(z) \in \SO_{n+1}$.

\begin{lemma}
\label{lemma:realpartQ}
For $z \in \Spin_{n+1} \subset \Cliff_{n+1}^0$,
let $Q = \Pi(z) \in \SO_{n+1}$ have eigenvalues
$\exp(\pm \theta_1 i), \ldots, \exp(\pm \theta_k i), 1, \ldots, 1$.
Then 
\[ \Re(z) = \pm
\cos\left(\frac{\theta_1}{2}\right) \cdots
\cos\left(\frac{\theta_k}{2}\right). \]
\end{lemma}

Notice, in particular, that $\Re(z) = 0$ if and only if
$-1$ is an eigenvalue of $Q$.

\begin{proof}
The function $\Re: \Spin_{n+1} \to \RR$ is
invariant under conjugation.
We may therefore assume that 
$Q = \Pi(z) \in \SO_{n+1}$ is of the form:
\begin{align*}
Q &= \diag\left(
\begin{pmatrix} \cos \theta_1 & -\sin \theta_1 \\
\sin \theta_1 & \cos \theta_1 \end{pmatrix}, \cdots ,
\begin{pmatrix} \cos \theta_k & -\sin \theta_k \\
\sin \theta_k & \cos \theta_k \end{pmatrix},
1, \ldots, 1 \right), \\
z &= \alpha_1\left(\theta_1\right)
\alpha_3\left(\theta_2\right)
\cdots
\alpha_{2k-1}\left(\theta_k\right) \\
&=
\left( \cos\left(\frac{\theta_1}{2}\right) +
\sin\left(\frac{\theta_1}{2}\right) \hat a_1 \right)
\cdots
\left( \cos\left(\frac{\theta_k}{2}\right) +
\sin\left(\frac{\theta_k}{2}\right) \hat a_{2k-1} \right).
\end{align*}
The result now follows from an explicit computation. 
\end{proof}


\goodbreak
\section{The group $\tilde B_{n+1}^{+}$}
\label{section:tildeB}

Let $B_{n+1}$ be the Coxeter-Weyl group
of signed permutation matrices;
let $B_{n+1}^{+} = B_{n+1} \cap \SO_{n+1}$.
For the double cover 
$\Pi: \Spin_{n+1} \to \SO_{n+1}$,
let
$\tilde B_{n+1}^{+} = \Pi^{-1}[B_{n+1}^{+}]$.
The finite group $\tilde B_{n+1}^{+} \subset \Spin_{n+1}$
is generated by 
\[ \acute a_i = \alpha_i\left(\frac{\pi}{2}\right)
= \frac{1+\hat a_i}{\sqrt{2}}, \quad
\grave a_i = (\acute a_i)^{-1} = \alpha_i\left(-\frac{\pi}{2}\right)
= \frac{1-\hat a_i}{\sqrt{2}},
 \quad i \in \nmesmo \] 
(of course, $\alpha_i = \alpha_i^{\Spin}$). 
Clearly, we have $\hat a_i = (\acute a_i)^2 = -(\grave a_i)^2$.
There is a natural short exact sequence
\begin{equation}
\label{equation:exactsequence}
1 \to \Quat_{n+1} \to \tilde B_{n+1}^{+} \to S_{n+1} \to 1. 
\end{equation}
If $\sigma = a_{i_1} \cdots a_{i_\ell}$ is a reduced word,
define
$\acute\sigma = \acute a_{i_1} \cdots \acute a_{i_\ell}
\in \tilde B_{n+1}^{+}$:
it turns out that different reduced words yield the same answer
(\cite{Goulart-Saldanha0}, Lemma 3.2).
If $\Pi = \Pi_{\tilde B_{n+1}^{+}, S_{n+1}}:
\tilde B_{n+1}^{+} \to S_{n+1}$ is as in
Equation \eqref{equation:exactsequence},
the set $\Pi^{-1}[\{\sigma\}] = \acute\sigma \Quat_{n+1} \subset
\tilde B_{n+1}^{+}$ is a coset (both left and right). 
Recall that the group $\cE_n = \{\pm 1\}^{\nmesmo}$
acts by automorphisms on $\Spin_{n+1}$
(see Equation \eqref{equation:cE}):
these are also automorphisms of $\tilde B_{n+1}^{+}$.
Furthermore, $\Re(z^E) = \Re(z)$
for all $z \in \tilde B_{n+1}^{+}$ and $E \in \cE_n$.

We now compute $\Re(z)$ for $z \in \tilde B_{n+1}^{+}$
using information about $\sigma = \Pi(z) \in S_{n+1}$.

\begin{remark}
\label{remark:Chebyshev}
The Chebyshev polynomials $T_k$ are recursively defined by
\[ T_0(x) = 1, \quad T_1(x) = x, \quad T_{k+1}(x) = 2x T_k(x) - T_{k-1}(x) \]
and satisfy $\cos(kt) = T_k(\cos(t))$.
Thus, for $k > 0$, the leading term of $T_k(x)$ is $2^{k-1} x^k$
and the independent term is $0$ if $k$ is odd and
$(-1)^{k/2}$ if $k$ is even.
\end{remark}

\begin{example}
\label{example:oddcycle}
Consider $z_0 \in \tilde  B_{n+1}^{+}$ such that
$Q_0 \in B_{n+1}^{+}$ is a cycle of odd length $k$.
More precisely, there exist $i_1, \ldots, i_k$ with
$(e_{i_j})^\top Q_0 = (e_{i_{j+1}})^\top$
for $1 \le j < k$, 
$(e_{i_k})^\top Q_0 = (e_{i_{1}})^\top$
and
$(e_{i})^\top Q_0 = (e_{i})^\top$
otherwise.
The eigenvalues of $Q_0$ are $1$ (with some large multiplicity)
and simple eigenvalues
\[
\exp\left(\pm\frac{2\pi i}{k}\right),
\exp\left(\pm\frac{4\pi i}{k}\right), \ldots, 
\exp\left(\pm\frac{(k-1)\pi i}{k}\right).
\]
It follows from Lemma \ref{lemma:realpartQ} that
$\Re(z) = \pm P$ where $P$ is the product
\[ P = \cos\left(\frac{\pi}{k}\right)  \cos\left(\frac{2\pi}{k}\right) 
\cdots
\cos\left(\frac{(k-3)\pi}{2k}\right)
\cos\left(\frac{(k-1)\pi}{2k}\right). \]
But $-P^4$ is the product of the roots of $T_{2k}(x) - 1 = 0$:
since the leading coefficient of $T_{2k}$ is $2^{2k-1}$
and the independent coefficient is $-1$
we have $P^4 = 2^{-2k+2}$.
We therefore have
$\Re(z_0) = \pm 2^{(-k+1)/2}$.
\end{example}

\begin{example}
\label{example:evencycle}
Consider $z_0 \in \tilde  B_{n+1}^{+}$ such that
$Q_0 \in B_{n+1}^{+}$ is a cycle of even length $k$.
More precisely, there exist $i_1, \ldots, i_k$ with
$(e_{i_j})^\top Q_0 = (e_{i_{j+1}})^\top$
for $1 \le j < k$, 
$(e_{i_k})^\top Q_0 = - (e_{i_{1}})^\top$
and
$(e_{i})^\top Q_0 = (e_{i})^\top$
otherwise.
By adapting the computations in Example \ref{example:oddcycle}
we have
$\Re(z_0) = \pm 2^{(-k+1)/2}$.
\end{example}

Let $\Diag_{n+1}^{+} = \Diag_{n+1} \cap \SO_{n+1}$.
Given a partition $X$ of the set $\llbracket n+1\rrbracket = \{1,\ldots,n+1\}$
we define a subgroup $H_{\Diag,X} \le \Diag_{n+1}^{+}$
of index $2^{|X|-1}$.
Given $X$, $H_{\Diag,X}$ is the set
of matrices $E \in \Diag_{n+1}^{+} \subset \SO_{n+1}$ such that,
if $A = \{i_1, i_2, \ldots, i_k\} \in X$ then
the product $E_{i_1i_1}E_{i_2i_2} \cdots E_{i_ki_k}$ equals $1$.
Let $H_X = \Pi^{-1}[H_{\Diag,X}] \le \Quat_{n+1}$,
where $\Pi: \Quat_{n+1} \to \Diag_{n+1}^{+}$
is the restriction of $\Pi: \Spin_{n+1} \to \SO_{n+1}$.

For a permutation $\sigma$, consider 
the partition $X_\sigma$ of  $\llbracket n+1\rrbracket$ into cycles.
Let $H_{\sigma} = H_{X_\sigma} \le \Quat_{n+1}$.
We have $|H_\sigma| = 2^{n+2-c}$.


\begin{example}
\label{example:Hsigma}
Take $n = 4$, $\sigma = 53421 = (15)(234) = a_1a_2a_3a_2a_1a_4a_3a_2a_1$.
The subgroup $H_{\Diag,\sigma}$ is spanned by
\[ \diag(-1,1,1,1,-1), \diag(1,-1,-1,1,1), \diag(1,1,-1,-1,1)
\in \Diag_5^{+}: \]
these diagonal matrices change an even number of signs in a single cycle.
Lifting to $H_\sigma$ we have the generators
\begin{gather*}
\hat a_1\hat a_2\hat a_3\hat a_4, \hat a_2, \hat a_3 \in \Quat_5, \\
H_\sigma = \{ \pm 1, \pm \hat a_2, \pm \hat a_3, \pm \hat a_2\hat a_3,
\pm \hat a_1\hat a_4, \pm \hat a_1\hat a_2\hat a_4, 
\pm \hat a_1\hat a_3\hat a_4, \pm \hat a_1\hat a_2\hat a_3\hat a_4 \}. 
\end{gather*}
Notice that have $c = 2$ and therefore $|H_\sigma| = 2^4 = 16$.
\end{example}

\begin{lemma}
\label{lemma:realpartsigma}
Consider $\sigma \in S_{n+1}$:
assume that $\sigma$ has $c = \nc(\sigma)$ disjoint cycles.
Choose $z_0 \in \acute\sigma\Quat_{n+1}$ such that $\Re(z_0) > 0$.
For $q \in \Quat_{n+1}$, we have 
\[ |\Re(qz_0)| = |\Re(z_0q)| =
\begin{cases}
2^{-(n+1-c)/2}, & q \in H_\sigma, \\
0, & q \notin H_\sigma.
\end{cases} \]
There are $2^{n+1-c}$ values of $q \in \Quat_{n+1}$
such that $\Re(qz_0) > 0$
(and similarly for $\Re(z_0 q)$).
Also, if $z_0$ is expanded in
the canonical basis as $z_0 = \sum_{q \in \HQuat_{n+1}} c_q q$
then $c_q \ne 0$ if and only if $q \in H_\sigma$.
\end{lemma}

\begin{proof}
For $z \in \acute\sigma\Quat_{n+1}$, take $Q = \Pi(z)$
and consider its cycles and submatrices $Q_1, \ldots, Q_c$
of dimensions $k_1, \ldots, k_c$ with $k_1 + \cdots + k_c = n+1$.
If $\det(Q_i) = -1$ for some $i$ then $-1$ is an eigenvalue of $Q_i$
and therefore of $Q$ which implies that $\Re(z) = 0$. 
Otherwise, apply Lemma \ref{lemma:realpartQ} to compute $\Re(z)$
as product, imitating Examples \ref{example:oddcycle} and \ref{example:evencycle}, 
to obtain $\Re(z) = \pm 2^{-(n+1-c)/2}$.
It is easy to relate these properties with the definition
of the subgroup $H_\sigma$.
\end{proof}

\begin{example}
\label{example:Hsigmaz0}
Take $n = 4$ and $\sigma = a_1a_2a_3a_2a_1a_4a_3a_2a_1$
as in Example \ref{example:Hsigma}.
A computation gives 
\[ \acute\sigma =
\frac{-\hat a_1 - \hat a_1\hat a_2
+ \hat a_1\hat a_3 - \hat a_1\hat a_2\hat a_3
- \hat a_4 + \hat a_2\hat a_4 
- \hat a_3\hat a_4 - \hat a_2\hat a_3\hat a_4}{2\sqrt{2}}: \]
notice that $1/\sqrt{8} = 2^{-(n+1-c)/2}$.
We may choose
\[ z_0 = \hat a_1 \acute\sigma = 
\frac{1 + \hat a_2 - \hat a_3 + \hat a_2\hat a_3 
- \hat a_1\hat a_4 + \hat a_1\hat a_2\hat a_4
- \hat a_1\hat a_3\hat a_4 - \hat a_1\hat a_2\hat a_3\hat a_4}{2\sqrt{2}}; \]
notice that the nonzero coefficients match the elements of $H_\sigma$.
\end{example}


\section{Ancestries}
\label{section:ancestry}

An \textit{ancestry} is a sequence $(\varrho_k)_{0 \le k \le \ell}$
of elements of $\tilde B_{n+1}^{+}$ such that:
\begin{enumerate}
\item{$\varrho_0 = \acute\eta$,
$\varrho_\ell \in \acute\eta \Quat_{n+1}$;}
\item{for all $k$, we have 
$\varrho_k = \varrho_{k-1}$ or
$\varrho_k = \varrho_{k-1} \acute a_{i_k}$ or
$\varrho_k = \varrho_{k-1} \hat a_{i_k}$;}
\item{the sequence $(\rho_k)$ defined by $\rho_k =
\Pi_{\tilde B_{n+1}^{+}, S_{n+1}}(\varrho_k)$
is a preancestry.}
\end{enumerate}
In this case, we say that the ancestry $(\varrho_k)$
corresponds to the preancestry $(\rho_k)$.
The third condition is equivalent to saying that 
$\Pi(\varrho_{k-1}) < \Pi(\varrho_{k-1}) a_{i_k}$ implies 
$\varrho_k = \varrho_{k-1} \acute a_{i_k}$ (for all $k$).

There are two sequences of integers
and a sequence of elements of $\Quat_{n+1}$
which conveniently describe an ancestry. 
The sequence $\xi: \llbracket \ell \rrbracket \to \{0,1,2\}$
is recursively defined by
$\varrho_k = \varrho_{k-1} (\acute a_{i_k})^{\xi(k)}$.
Equivalently, we have
\[ \varrho_k = (\acute a_{i_1})^{\xi(1)} \cdots (\acute a_{i_k})^{\xi(k)}. \]
The sequence $(q_k)_{0 \le k \le \ell}$, $q_k \in \Quat_{n+1}$,
is defined by $\acute\rho_k q_k = \varrho_k$,
so that in particular $q_\ell = \grave\eta \varrho_\ell$.

The third sequence may seem at first artificial,
but we shall later see that it appears naturally in our problem. 
Given an ancestry (and therefore $\xi$, $(\rho_k)$ and $(q_k)$),
define a sequence
$\varepsilon:  \llbracket \ell \rrbracket \to \{\pm 1, \pm 2\}$ by
\begin{equation}
\label{equation:xi2varepsilon}
\varepsilon(k) = \begin{cases}
-2, & \xi(k) = 1, \; \rho_k < \rho_{k-1}, \\
+2, & \xi(k) = 1, \; \rho_k > \rho_{k-1}, \\
(1-\xi(k))[\hat a_{i_k},q_{k-1}], & \xi(k) \ne 1.
\end{cases}
\end{equation}
Conversely, $\xi$ and $(\varrho_k)$ can be recovered
from $\varepsilon$ by:
\begin{equation}
\label{equation:varepsilon2xi}
\xi(k) = \begin{cases}
0, & 
\varepsilon(k) = [\hat a_{i_k},q_{k-1}], \\
2, & 
\varepsilon(k) = - [\hat a_{i_k},q_{k-1}], \\
1, & 
|\varepsilon(k)| = 2; \end{cases}
\qquad
\varrho_k = \varrho_{k-1} (\acute a_{i_k})^{\xi(k)}.
\end{equation}
Here $[\hat a_{i_k},q_{k-1}] =
(\hat a_{i_k})^{-1}q_{k-1}^{-1} \hat a_{i_k} q_{k-1} \in \{\pm 1\}$
is the commutator, as in Equation \eqref{equation:commutator}.
We also have
$q_k = (\acute\rho_k)^{-1} \varrho_k =
(\acute a_{i_k})^{-\sign(\varepsilon(k))} q_{k-1} \acute a_{i_k}$.
Thus, given the reduced word, each of the sequences
$(\varrho_k)$, $\xi$ and $\varepsilon$
allows us to obtain $(q_k)$ and the other two sequences;
given the preancestry and $(q_k)$,
we can also obtain the three sequences above.
We thus consider these three sequences to be
alternate descriptions of an ancestry.
An ancestry is represented in a wiring diagram by
the sequence $\varepsilon$:
$-2$ and $+2$ are represented by black and white squares, 
respectively, as for preancestries;
$-1$ and $+1$ are represented by black and white disks.
For an ancestry $\varepsilon$, set
$P(\varepsilon) = \acute\sigma q_\ell^{-1} \in \acute\sigma \Quat_{n+1}^{+}$;
we therefore have
$\varrho_\ell = \acute\eta (P(\varepsilon))^{-1} \acute\sigma$.

\begin{lemma}
\label{lemma:P}
For any ancestry $\varepsilon$, we have
\[ P(\varepsilon) =
(\acute a_{i_1})^{\sign(\varepsilon(1))} \cdots
(\acute a_{i_\ell})^{\sign(\varepsilon(\ell))}. \]
\end{lemma}

\begin{proof}
Let $\sigma_k = a_{i_1}\cdots a_{i_k}$;
we prove by induction that
\[ \acute\sigma_k q_k^{-1} = 
(\acute a_{i_1})^{\sign(\varepsilon(1))} \cdots
(\acute a_{i_k})^{\sign(\varepsilon(k))}. \]
The result follows by taking $k = \ell$.
\end{proof}

\begin{example}
\label{example:ancestry}
Let $\varepsilon_1$, $\varepsilon_2$ and $\varepsilon_3$
be the three ancestries shown in Figure \ref{fig:ancestry}.

We have $\varepsilon_1 = (-1,+1,-1,-1,+1,-1,+1)$:
for this ancestry we have $\rho_k = \eta$ for all $k$.
We apply the recursion in Equation \eqref{equation:varepsilon2xi}.
We have $q_0 = 1$ and therefore $[\hat a_3, q_0] = 1$:
thus, $\xi(1) = 2$, $\varrho_1 = \acute\eta\hat a_3$
and $q_1 = \hat a_3$.
We have $[\hat a_2, q_1] = -1$ and therefore $\xi(2) = 2$,
$\varrho_2 = \acute\eta\hat a_3\hat a_2$
and $q_2 = -\hat a_2\hat a_3$.
We have $[\hat a_1, q_2] = -1$ and therefore $\xi(3) = 0$,
$\varrho_3 = -\acute\eta\hat a_2\hat a_3$
and $q_3 = -\hat a_2\hat a_3$.
Proceed in this manner to obtain
$\xi = (2,2,0,0,2,2,0)$, $\varrho_7 = -\acute\eta\hat a_2\hat a_4$,
$q_7 = -\hat a_2\hat a_4$ 
and therefore 
\[
P(\varepsilon_1) = \acute\sigma (-\hat a_2\hat a_4) = 
\grave a_3 \acute a_2 \grave a_1 \grave a_2 \acute a_4 \grave a_3 \acute a_2,
\]
consistently with Lemma \ref{lemma:P}.
\end{example}




Given a preancestry $\varepsilon_0$
and $z \in \acute\sigma \Quat_{n+1}$,
let $L_{\varepsilon_0}(z)$ be the set of ancestries
$\varepsilon$ corresponding to $\varepsilon_0$ and
satisfying $P(\varepsilon) = z$;
let $\NL_{\varepsilon_0}(z) = |L_{\varepsilon_0}(z)|$.

\begin{remark}
\label{remark:NLN}
Recall that $N_{\varepsilon_0}(z)$ appears in the statement
of Theorem \ref{theo:two}.
We shall see in Section \ref{section:stratification} that
$N_{\varepsilon_0}(z) = \NL_{\varepsilon_0}(z)$.
The following result is therefore one important step
towards the proof of Theorem \ref{theo:two}. 
\end{remark}

\begin{lemma}
\label{lemma:preminus}
For any $z \in \acute\sigma\Quat_{n+1}$, we have 
$\NL_{\varepsilon_0}(z) - \NL_{\varepsilon_0}(-z) =
2^{\frac{\ell-2d}{2}} \Re(z)$.
\end{lemma}

\begin{proof}
Set $\delta = \ell - 2d$. 
We define the element of the Clifford algebra
\[ S = \sum_\varepsilon P(\varepsilon) =
\sum_z \NL_{\varepsilon_0}(z) z \in \Cliff_{n+1}^{0}; \]
the first summation is understood to be over all ancestries $\varepsilon$
associated with $\varepsilon_0$.
We claim that $S = 2^{\frac{\delta}{2}} \cdot 1$
(where $1 \in \Cliff_{n+1}^{0}$ is the unit of the Clifford algebra).
Indeed, $S$ can be written as a product
$S = 2^{\frac{\delta}{2}} w_1w_2\cdots w_\ell$
where $w_k \in \Cliff^0_{n+1}$ is defined as follows.
If $\varepsilon_0(k) = 2$ set $w_k = \acute a_{i_k}$;
if $\varepsilon_0(k) = -2$ set $w_k = \grave a_{i_k}$;
if $\varepsilon_0(k) = 0$
set $w_k = (\acute a_{i_k} + \grave a_{i_k})/\sqrt{2} = 1 \in \Cliff^0_{n+1}$.
Set $s_k = w_1w_2\cdots w_k$:
we have $S = 2^{\frac{\delta}{2}} s_\ell$.
From the recursive definition of $\rho_k$
in Equation \eqref{equation:preancestry}
we have $s_k = \grave\eta\acute\rho_k$.
It follows that $s_\ell = \grave\eta \acute\eta = 1$,
completing the proof of the claim.

On the other hand, since two distinct elements of $\acute\sigma\Quat_{n+1}$
are either antipodal or orthogonal, we have
$2^{\frac{\delta}{2}} \Re(z) = \langle z, S \rangle = 
\NL_{\varepsilon_0}(z) - \NL_{\varepsilon_0}(-z)$,
completing the proof of the lemma.
\end{proof}

Given a preancestry $\varepsilon_0$,
we define a partition $X_{\varepsilon_0}$ 
which is coarser than $X_\sigma$.
Recall that $X_\sigma$ is the partition of $\llbracket n+1\rrbracket$
into cycles of $\sigma$.
If $k \in \llbracket\ell\rrbracket$ is such that $\varepsilon_0(k) = 0$
and the $k$-th crossing is $(i_0,i_1)$ then
$\{i_0, i_1\}$ is contained in some $A \in X_{\varepsilon_0}$.
The partition $X_{\varepsilon_0}$ is the finest partition
respecting the conditions above.
Let 
\begin{equation}
\label{equation:Hvarepsilon}
H_{\varepsilon_0} = H_{X_{\varepsilon_0}} \le \Quat_{n+1}:
\end{equation}
we clearly have $H_{\sigma} \le H_{\varepsilon_0}$.

\begin{example}
\label{example:Hvarepsilon0}
Let $n = 4$ and
$\sigma = \eta = a_1a_2a_1a_3a_2a_1a_4a_3a_2a_1 = (15)(24)(3)$.
We have $X_\sigma = \{\{1,5\},\{2,4\},\{3\}\}$ and
$H_\sigma = \{\pm 1, \pm \hat a_2\hat a_3, 
\pm \hat a_1\hat a_4, \pm \hat a_1\hat a_2\hat a_3\hat a_4\}$.

As seen in Figure \ref{fig:preancestry54321},
there is only one preancestry of top codimension $d_{\max} = 4$,
which is $(-2,-2,-2,-2,0,2,0,2,2,2)$:
we then have $X = X_\sigma$ and therefore $H = H_\sigma$.

There are $5$ preancestries of codimension $3$.
For $(-2,-2,0,-2,0,0,0,2,2,2)$, we have
$X = \{\{1,5\},\{2,3,4\}\}$ and
$H$ is generated by $H_\sigma$ and $\hat a_2$ so that
\[ H = \{\pm 1, \pm \hat a_2, \pm \hat a_3, \pm \hat a_2\hat a_3,
\pm \hat a_1\hat a_4, \pm \hat a_1\hat a_2\hat a_4,
\pm \hat a_1\hat a_3\hat a_4, \pm \hat a_1\hat a_2\hat a_3\hat a_4 \}. \]
For $(-2,0,2,-2,0,-2,0,2,0,2)$, we have
$X = \{\{1,3,5\},\{2,4\}\}$ and
$H$ is generated by $H_\sigma$ and $\hat a_1\hat a_2$.
For the other three, we have
$X = \{\{1,2,4,5\},\{3\}\}$ and
$H$ is generated by $H_\sigma$ and $\hat a_1$.

There are $10$ preancestries of codimension $2$.
We have $X = \{\{1,2,3,4,5\}\}$
and therefore $H = \Quat_5$
for $9$ of them.
The exception is $(0,-2,-2,0,0,2,0,0,2,0)$ 
for which $X = \{\{1,2,4,5\},\{3\}\}$.
For all $6$ preancestries of codimension $1$
and for the preancestry of codimension $0$,
we have $X = \{\{1,2,3,4,5\}\}$ and $H = \Quat_5$.
\end{example}

The following result is another important step
towards Theorem \ref{theo:two}.

\begin{lemma}
\label{lemma:preplus}
Consider a preancestry $\varepsilon_0$ and 
the subgroup $H_{\varepsilon_0} \le \Quat_{n+1}$.
Choose $z_0 \in \acute\sigma\Quat_{n+1}$ with $\Re(z) > 0$.
For $z = q z_0$, we have
\[ \NL_{\varepsilon_0}(z) + \NL_{\varepsilon_0}(-z) = 
\begin{cases}
2^{\ell-2d+1}/|H_{\varepsilon_0}|, & q \in H_{\varepsilon_0}, \\
0, & q \notin H_{\varepsilon_0}.
\end{cases} \]
\end{lemma}

\begin{proof}
For this proof it is convenient to introduce
the (associative but not commutative) group ring $\RR[B_{n+1}^{+}]$.
Elements of $\RR[B_{n+1}^{+}]$ are
formal linear combinations $\sum_{Q \in B_{n+1}^{+}} c_Q Q$, $c_Q \in \RR$.
Thus, the set $B_{n+1}^{+}$ is a real basis for $\RR[B_{n+1}^{+}]$;
multiplication is inherited from $B_{n+1}^{+}$.
We use names of elements of $\tilde B_{n+1}^{+}$,
such as $\acute a_i$,
as abbreviations of $1\cdot Q \in \RR[B_{n+1}^{+}]$, $Q = \Pi(\acute a_i)$. 
We thus have (in $\RR[B_{n+1}^{+}]$) $(\hat a_i)^2 = 1$,
$\acute a_i + \grave a_i =
\acute a_i (1+\hat a_i) = \grave a_i (1+\hat a_i)$
and $(1+\hat a_i)^2 = 2(1+\hat a_i)$.

The commutative subring $\RR[\Diag^{+}_{n+1}] \subset \RR[B_{n+1}^{+}]$
will also interest us.
Given a subgroup $H \le \Diag^{+}_{n+1}$,
define 
\[ \bp_H = \frac{1}{|H|} \sum_{E \in H} 1\cdot E \in \RR[\Diag^{+}_{n+1}]. \]
We have $\bp_H^2 = \bp_H$ so that multiplication by $\bp_H$
on $\RR[\Diag^{+}_{n+1}]$ is the projection onto the subring $R_H$:
\[ \sum_{E \in \Diag^{+}_{n+1}} c_E \cdot E \in R_H \quad\iff\quad
\forall E_0, E_1 \in \Diag^{+}_{n+1}, \;
(E_0E_1^{-1} \in H \;\to\; c_{E_0} = c_{E_1}). \]
If $H_0, H_1 \le \Diag^{+}_{n+1}$ we have
$\bp_{H_0}\bp_{H_1} = \bp_{H_0H_1}$.

For a transposition $\tau = (i j) \in S_{n+1}$ (with $i < j$), write
\[ E_\tau =
\Pi(\hat a_i\hat a_{i+1} \cdots \hat a_{j-1}) \in \Diag^{+}_{n+1}, \]
the diagonal matrix with entries $(E_\tau)_{i,i} = (E_\tau)_{j,j} = -1$
and $(E_\tau)_{k,k} = +1$ for $k \notin \{i,j\}$.
Write $H_\tau = \{I, E_\tau\}$ and $\bp_\tau = \bp_{H_\tau}$.
A simple computation verifies the following identities in $\RR[B_{n+1}^{+}]$:
\begin{equation}
\label{equation:movep}
\acute a_i \bp_\tau = \bp_{\tau^{a_i}} \acute a_i, \qquad
\grave a_i \bp_\tau = \bp_{\tau^{a_i}} \grave a_i. 
\end{equation}
Notice that $\tau^{a_i} = a_i\tau a_i$ is also a transposition.

We need to compute $N_z = \NL_{\varepsilon_0}(z) + \NL_{\varepsilon_0}(-z)$,
which is the coefficient $c_Q$ of $Q = \Pi(z) \in B_{n+1}^{+}$ in 
$S = s_\ell \in \RR[B_{n+1}^{+}]$,
defined recursively by
\[ s_0 = 1, \qquad
s_k = \begin{cases}
s_{k-1}\cdot (\acute a_{i_k} + \grave a_{i_k}), & \varepsilon_0(k) = 0, \\
s_{k-1}\cdot (\acute a_{i_k})^{\sign(\varepsilon_0(k))}, &
\varepsilon_0(k) \ne 0.
\end{cases} \]
Let $\delta = \ell - 2d$.
Choose $\varepsilon$ with $P(\varepsilon) = z_0$
and rewite $S$ as $S = 2^{\delta} \tilde S = 2^{\delta} \tilde s_\ell$,
\[ \tilde s_0 = 1, \qquad
\tilde s_k = \begin{cases}
\tilde s_{k-1}\cdot
(\acute a_{i_k})^{\sign(\varepsilon(k))}\cdot \bp_{a_{i_k}}, &
\varepsilon_0(k) = 0, \\
\tilde s_{k-1}\cdot (\acute a_{i_k})^{\sign(\varepsilon_0(k))}, &
\varepsilon_0(k) \ne 0.
\end{cases} \]
As in the proof of Lemma \ref{lemma:transpositions},
let $k_1, \ldots, k_\delta$ be the unmarked crossings.
Again assume that the unmarked crossing $k_j$ is 
$(\iota_{j,0}, \iota_{j,1}) \in \Inv(\sigma)$.
Let $\tau_j = (\iota_{j,0} \iota_{j,1})$.
Applying Equation \eqref{equation:movep}
to bring the terms $\bp_\ast$ to the left,
we have
\[ \tilde S = \bp_{\tau_1}\cdots \bp_{\tau_\delta} \Pi(z_0)
\in \RR[B_{n+1}^{+}]. \]
On the other hand, we have
$H_{\tau_1}\cdots H_{\tau_\ell} = H_{\Diag,X_{\varepsilon_0}}$
and therefore
\[ \bp_{\tau_1}\cdots \bp_{\tau_\delta} =
\bp_{H_{\Diag,X_{\varepsilon_0}}} =
\frac{2}{|H_{\varepsilon_0}|}
\sum_{\Pi(E) \in H_{\varepsilon_0}}
1\cdot E \in \RR[\Diag_{n+1}^{+}] \]
and $S = 2^\delta \bp_{H_{\Diag,X_{\varepsilon_0}}} \Pi(z_0)
\in \RR[B_{n+1}^{+}]$,
completing the proof.
\end{proof}


\section{Thin ancestries}
\label{section:thin}

Assume given a permutation $\sigma \in S_{n+1}$ and a reduced word for $\sigma$. 
An ancestry $\varepsilon$ of dimension $0$
is \textit{thin} if $i_{k_0} = i_{k_1}$ implies
$\varepsilon(k_0) = \varepsilon(k_1)$.
Otherwise, an ancestry is \textit{thick}.
There are therefore $2^{n-b}$ thin ancestries, where $b = \block(\sigma)$.
Consistently with Remark \ref{remark:blockfree},
we assume from now on that $\sigma$ does not block, i.e., that $b = 0$.
We shall see in Sections \ref{section:firstexamples} and later
that thin ancestries correspond to certain
contractible connected components of $\BL_\sigma$.

For now, we are interested in counting thin solutions
for equations of the form $P(\varepsilon) = z$.
In other words: let $\varepsilon_0$ be the empty preancestry
and consider a fixed element $z \in \acute\sigma\Quat_{n+1}$.
By definition, there are $\NL_{\varepsilon_0}(z)$ ancestries
$\varepsilon$ corresponding to $\varepsilon_0$
and satisfying $P(\varepsilon) = z$.
Lemmas \ref{lemma:preminus} and \ref{lemma:preplus}
allow us to compute $\NL_{\varepsilon_0}(z)$.
Among these, let $\NL_{\thin}(z)$
be the number of such ancestries $\varepsilon$ which are thin.
We want to compute  $\NL_{\thin}(z)$.

Recall that
the multiplicative abelian  group $\cE_n$
acts as an automorphism group of $\Cliff^0_{n+1}$
(see Equation \eqref{equation:cE} and Remark \ref{remark:cE}).
For $z \in \Spin_{n+1}$, define $\cE_z \subseteq \cE_n$
as the isotropy group of $z$, i.e., 
\begin{equation}
\label{equation:isotropy}
\cE_z = \{ E \in \cE_n \;|\; z^{E} = z \}.
\end{equation}
We want to study the group $\cE_{\acute\sigma} \le \cE_n$.

\begin{lemma}
\label{lemma:cEsigma}
Given $\sigma \in S_{n+1}$,
let $c = \nc(\sigma)$ be the number of cycles of $\sigma$.
We have $|\cE_{\acute\sigma}| = 2^{\tilde c}$
where $\tilde c \in \ZZ$, $c-2 \le \tilde c \le c$.
\end{lemma}

\begin{remark}
\label{remark:cEsigma}
The value of $c - \tilde c \in \{0,1,2\}$
can be deduced by following the proof.
It does not, however, appear to have a simple formula.
Further investigation is desireable.
\end{remark}

\begin{proof}[Proof of Lemma \ref{lemma:cEsigma}]
Let $\tilde\cE_z = \{ E \in \cE_n \;|\; z^{E} = z \}$
so that $\cE_z \le \tilde\cE_z$.
If $Q = \Pi(z)$ (and $\Pi: \Spin_{n+1} \to \SO_{n+1}$)
we have $\tilde\cE_z = \cE_Q = \{ E \in \cE_n \;|\; Q^{E} = Q \}$.

Take $z = \acute\sigma$ so that $Q \in B_{n+1}^{+}$
is a signed permutation matrix.
A diagonal matrix $E \in \Diag_{n+1}$ commutes with $Q$
if and only if $E_{ii} = E_{jj}$ whenever $i$ and $j$
are in the same cycle of $\sigma$.
We may therefore choose one sign per cycle,
for a total of $2^c$ elements of $\Diag_{n+1}$.
We have $2^{c_1}$ such elements of $\Diag_{n+1}^{+}$
where $c-1 \le c_1 \le c$.
Thus, $|\tilde\cE_z| = |\cE_Q| = 2^{c_1}$.
Finally, $|\cE_{\acute\sigma}| = 2^{\tilde c}$, as claimed. 
\end{proof}

Let $\acute\sigma^{\cE_n}$ be the \textit{orbit}
of $\acute\sigma$ under $\cE_n$:
\[ \acute\sigma^{\cE_n} = \{ \acute\sigma^E, E \in \cE_n \}. \]
The following result is now easy.

\begin{lemma}
\label{lemma:Nthin}
Let $\tilde c$ be as in Lemma \ref{lemma:cEsigma}.
For $z \in \acute\sigma\Quat_{n+1}$, we have
\[ \NL_{\thin}(z) = \begin{cases}
2^{n-\tilde c}, & z \in \acute\sigma^{\cE_n}, \\
0, & z \notin \acute\sigma^{\cE_n}. \end{cases} \]
Furthermore, $|\acute\sigma^{\cE_n}| = 2^{\tilde c}$.
\end{lemma}

\begin{proof}
Each $E\in\cE_n$ gives us a thin ancestry $\varepsilon$
with $P(\varepsilon) = \acute\sigma^E$.
The result follows by the usual combinatorics
of actions over finite sets.
\end{proof}

\begin{coro}
\label{coro:thick}
Consider $\sigma \in S_n$ which does not block and
consider the empty preancestry $\varepsilon_0$. 
If $\ell = \inv(\sigma) > 2n+2$ then
for all $z \in \acute\sigma\Quat_{n+1}$ we have
$\NL_{\varepsilon_0}(z) > \NL_{\thin}(z)$.
\end{coro}

\begin{remark}
\label{remark:thicketa}
In particular, for $\sigma = \eta$ and $n > 4$ 
the above condition holds.
\end{remark}

\begin{proof}[Proof of Corollary \ref{coro:thick}]
From Lemma \ref{lemma:preplus},
$\NL_{\varepsilon_0}(z) + \NL_{\varepsilon_0}(-z) \ge 2^{\ell-n}$.
From Lemma \ref{lemma:preminus},
$\NL_{\varepsilon_0}(z) - \NL_{\varepsilon_0}(-z) \ge -2^{\ell/2}$.
From Lemma \ref{lemma:Nthin}, $\NL_{\thin} \le 2^n$.
Thus $\NL_{\varepsilon_0}(z) - \NL_{\thin} \ge
2^{\ell-n-1} - 2^{\ell/2-1} - 2^n \ge
2^{\ell-n-1} - 2^{\ell/2} > 0$.
\end{proof}


\section{Bruhat cells in the spin group}
\label{section:gsobis}

This section contains further review of notation and results
from \cite{Goulart-Saldanha0}.
We define the Bruhat cells in the orthogonal group.
For $\sigma \in S_{n+1}$,
\[ \Bru_\sigma^{\SO} = \{ Q \in \SO_{n+1} \;|\;
\exists U_0, U_1 \in \Up_{n+1}, Q = U_0 P_\sigma U_1 \}; \]
compare with Equation \eqref{equation:BL}.
Let $\Bru_\sigma = \Pi^{-1}[\Bru_\sigma^{\SO}] \subset \Spin_{n+1}$
(where $\Pi: \Spin_{n+1} \to \SO_{n+1}$ is the double cover).
The set $\Bru_\sigma$ has $2^{n+1}$ connected components,
each one containing an element $z$ of $\acute\sigma \Quat_{n+1}$.
For $z \in \tilde B_{n+1}^{+}$,
let $\Bru_z$ be the connected component of $\Bru_{\Pi(z)}$
containing $z$.
The set $\Bru_z$ is a smooth contractible submanifold of $\Spin_{n+1}$
of dimension $\ell = \inv(\sigma)$, $\sigma = \Pi(z)$.
We have
\[ \Bru_\sigma = \bigsqcup_{z \in \acute\sigma \Quat_{n+1}} \Bru_z,
\qquad
 \Spin_{n+1} = \bigsqcup_{z \in \tilde B_{n+1}^{+}} \Bru_z. \]
For a reduced word $\sigma = a_{i_1} \cdots a_{i_\ell}$,
an explicit parametrization of $\Bru_{\acute\sigma}$
is given by $\Phi: (0,\pi)^\ell \to \Bru_{\acute\sigma}$
(see Theorem 1, \cite{Goulart-Saldanha0}):
\[ \Phi(\theta_1, \ldots , \theta_\ell) =
\alpha_{i_1}(\theta_1) \cdots \alpha_{i_\ell}(\theta_\ell). \]
The (strong) Bruhat order in $S_{n+1}$ can be defined by
\[ \sigma_0 \le \sigma_1 \quad\iff\quad
\Bru_{\sigma_0} \subseteq \overline \Bru_{\sigma_1}. \]
If $z \in \Bru_{z_0}$,
$\sigma_0 = \Pi(z_0) < \sigma_1 = \sigma_0 a_i$
then $z \alpha_i(\theta) \in \Bru_{z_0 \acute a_i}$ for $\theta \in (0,\pi)$.
The order here is the strong Bruhat order:
if $\inv(\sigma_0) = \ell$ then $\inv(\sigma_1) = \ell + 1$
and if $\sigma_0 = a_{i_1} \cdots a_{i_\ell}$ is a reduced word
then so is
$\sigma_1 = a_{i_1} \cdots a_{i_\ell} a_i$.
On the other hand, if $z \in \Bru_{z_0}$,
$\sigma_1 < \sigma_0 = \Pi(z_0) = \sigma_1 a_i$
then there exists $\theta = \Theta_i(z) \in (0,\pi)$
such that
$z \alpha_i(-\theta) \in \Bru_{z_1}$, $z_0 = z_1 \acute a_i$.
We also have
$z \alpha_i(\tilde\theta) \in \Bru_{z_1}$
for all $\tilde\theta \in (-\theta,\pi-\theta)$.

We define a partial order in $\tilde B_{n+1}^{+}$,
the \textit{lifted Bruhat order}:
\begin{equation}
\label{equation:bruhatB}
z_0 \le z_1 \quad\iff\quad
\Bru_{z_0} \subseteq \overline \Bru_{z_1}. 
\end{equation}
Clearly, $z_0 \le z_1$ implies $\Pi(z_0) \le \Pi(z_1)$
but the converse does not hold:
for instance, $\Pi(z_0) = \Pi(z_1)$ and $z_0 \le z_1$
together imply $z_0 = z_1$.
We also define a partial order on the set of ancestries
(for a given permutation $\sigma$).
Given two ancestries $\varepsilon$ and $\tilde\varepsilon$,
let $(\varrho_k)$ and $(\tilde\varrho_k)$ be
as in Section \ref{section:ancestry}.
We define
\begin{equation}
\label{equation:poset}
\varepsilon \preceq \tilde\varepsilon \quad\iff\quad
(\forall k,  \varrho_k \le \tilde\varrho_k).
\end{equation}
Notice that $\varepsilon \preceq \tilde\varepsilon$
implies $P(\varepsilon) = P(\tilde\varepsilon)$.
The fact that this is a partial order is straightforward.
A set $U$ of ancestries is an \textit{upper set} 
if $\varepsilon \in U$ and $\varepsilon \preceq \tilde\varepsilon$
imply $\tilde\varepsilon \in U$.
The upper set generated by $\varepsilon$ is
$U_{\varepsilon} =
\{ \tilde\varepsilon \;|\; \varepsilon \preceq \tilde\varepsilon \}$.

\begin{remark}
\label{remark:poset}
If $\varepsilon$ is an ancestry and $\dim(\varepsilon) = 0$
then $\varepsilon$ is $\preceq$-maximal.
Conversely, consider $\varepsilon$ with $\dim(\varepsilon) > 0$.
Define $\tilde\varepsilon$ by
$\tilde\varepsilon(k) = \sign(\varepsilon(k))$:
we have $\varepsilon \prec \tilde\varepsilon$.
In a wiring diagram, the ancestry $\tilde\varepsilon$ is obtained
by replacing each square by a disk of the same color.
\end{remark}

\begin{figure}[ht]
\begin{center}
\includegraphics[scale=0.15]{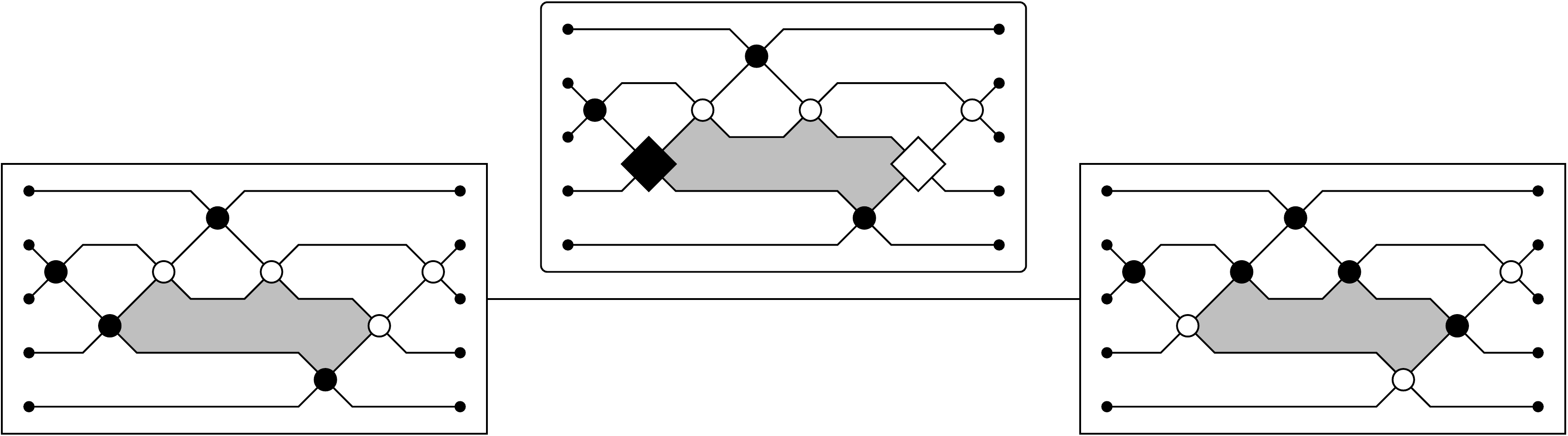}
\end{center}
\caption{An ancestry of dimension one and the upper set generated by it.  }
\label{fig:posetdimone}
\end{figure}

\begin{example}
\label{example:posetdimone}
If $\dim(\varepsilon) = 1$ then
the upper set $U_\varepsilon$ generated by $\varepsilon$
consists of $\varepsilon$ itself and two ancestries of dimension $0$.
One of them is $\tilde\varepsilon = \sign \circ \varepsilon$,
obtained by replacing the two squares by disks of the same color,
as in Remark \ref{remark:poset}.
The second one is obtained from $\tilde\varepsilon$
by changing the colors of all disks on the boundary
of the region corresponding to $\varepsilon$.
Figure \ref{fig:posetdimone} shows an example.
\end{example}

If $\dim(\varepsilon) > 1$ then
the description of the upper set $U$ generated by $\varepsilon$ 
is not so direct.
Recall from Remark \ref{remark:lowdimpreancestry}
that preancestries (and therefore ancestries) of dimension $2$
are classified as type I or type II.

\begin{figure}[ht]
\begin{center}
\includegraphics[scale=0.15]{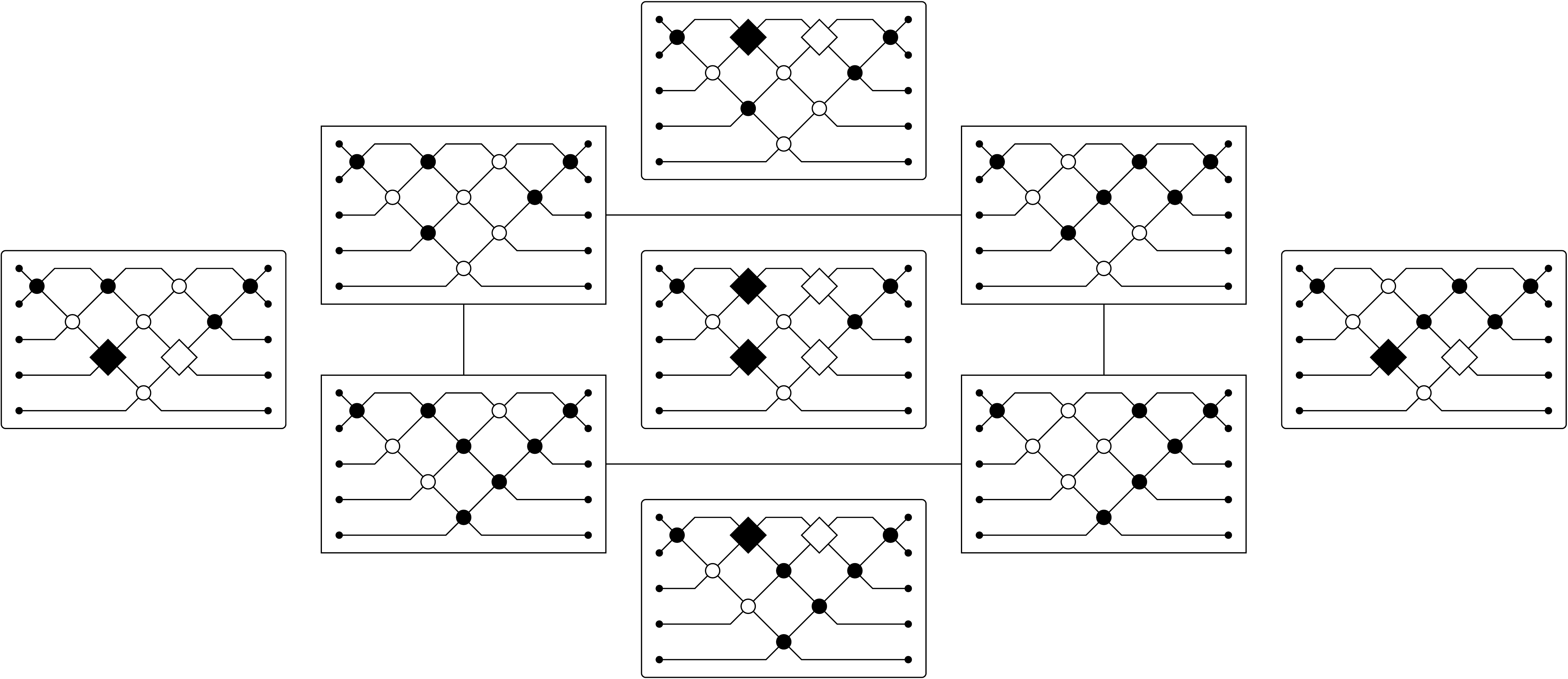}
\end{center}
\caption{The upper set generated by an ancestry of dimension two, type I. }
\label{fig:posetdimtwotypeI}
\end{figure}

\begin{example}
\label{example:posetdimtwo}
Let $\varepsilon$ be an ancestry of dimension $2$, type I.
The set $U_\varepsilon$ contains $4$ elements of dimension $0$,
$4$ elements of dimension $1$ and one element of dimension $2$
(which is $\varepsilon$).
An example is shown in Figure \ref{fig:posetdimtwotypeI}.

\begin{figure}[ht]
\begin{center}
\includegraphics[scale=0.2]{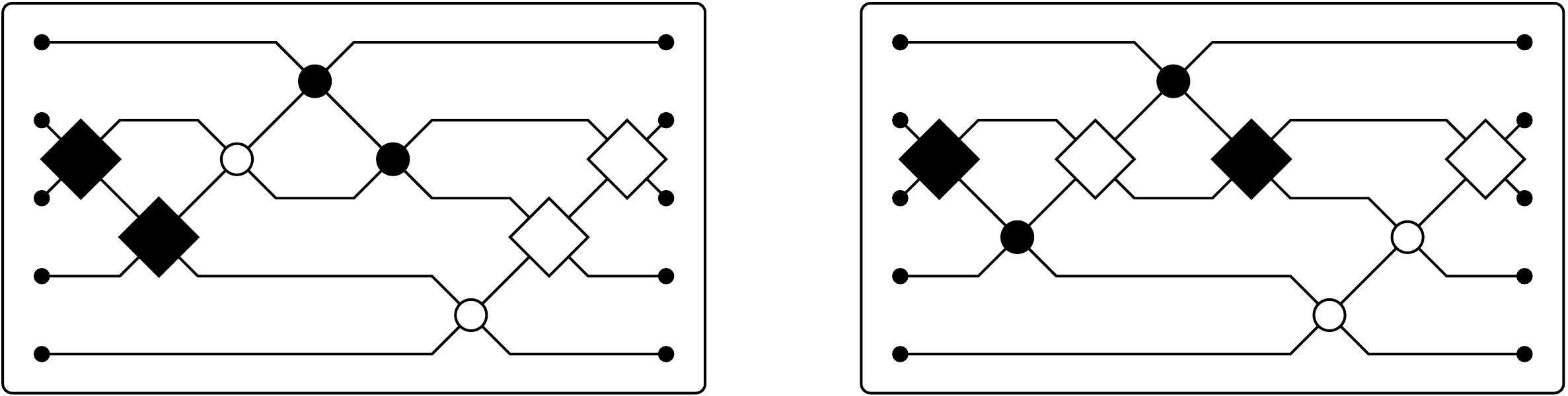}
\end{center}
\caption{Two ancestries of dimension two, types II and I. }
\label{fig:posetdimtwotypeII}
\end{figure}

If $\varepsilon$ is an ancestry of dimension $2$, type II,
the structure of $U_\varepsilon$ is somewhat more complicated.
Indeed, let $\varepsilon$ and $\tilde{\varepsilon}$
be the first and second ancestries in Figure \ref{fig:posetdimtwotypeII},
respectively.
We have
$\varepsilon \preceq \tilde\varepsilon$,
$\tilde\varepsilon \not\preceq \varepsilon$
and $\dim(\varepsilon) = \dim(\tilde\varepsilon) = 2$.
\end{example}


\section{The sets $\BL_z$}
\label{section:BLz}

For a matrix $L \in \Lo_{n+1}^1$, perform the usual $QR$ factorization:
\[ L = QR, \quad Q \in \SO_{n+1}, R \in \Up_{n+1}^{+}; \]
here $\Up_{n+1}^{+} \subset \Up_{n+1}$
is the group of upper triangular matrices
with positive diagonal entries.
This defines a smooth map:
\[ \bQ_{\SO}: \Lo_{n+1}^1 \to \SO_{n+1}, \qquad \bQ_{\SO}(L) = Q. \]
Lift this map to define $\bQ: \Lo_{n+1}^1 \to \Spin_{n+1}$
with $\bQ(I) = 1$.
The set $\cU_1 = \bQ[\Lo_{n+1}^1] \subset \Spin_{n+1}$
is an open contractible neighborhood of $1 \in \Spin_{n+1}^1$.
We have $\cU_1 = \grave\eta \Bru_{\acute\eta}$.
In other words, $\cU_1$ is a top-dimensional Bruhat cell
for the basis described by $\grave\eta$,
which is, up to signs, $e_{n+1}, e_n, \ldots, e_2,e_1$.
The inverse map $\bL = \bQ^{-1}: \cU_1 \to \Lo_{n+1}^1$
is also a smooth diffeomorphism:
it corresponds to the $LU$ factorization.

We are now ready to define the sets $\BL_z$,
already mentioned in Equation \eqref{equation:BLz} in the Introduction:
for $z \in \tilde B_{n+1}^{+}$, take 
\begin{equation}
\label{equation:BLz2}
\BL_z = \bQ^{-1}[\Bru_z] = \bQ^{-1}[\Bru_z \cap \, \grave\eta \Bru_{\acute\eta}]
\subseteq \Lo_{n+1}^1.
\end{equation}
The set $\BL_z$ is diffeomorphic to $\Bru_z \cap \, \grave\eta \Bru_{\acute\eta}$,
the intersection of two Bruhat cells for different bases in $\Spin_{n+1}$.
Indeed, the diffeomorphisms are given by the restrictions
of $\bL$ and $\bQ$. 

\begin{example}
\label{example:aba}
As in Example \ref{example:abaintro},
take $n = 2$ and
$\eta = a_1a_2a_1 = 321 = (13) \in S_3$.
We have
\[ \acute\eta = \frac{\hat a_1 + \hat a_2}{\sqrt{2}}, \qquad
\acute\eta\Quat_3 = \left\{
\frac{\pm 1 \pm \hat a_1 \hat a_2}{\sqrt{2}},
\frac{\pm \hat a_1 \pm \hat a_2}{\sqrt{2}}
\right\}; \]
sign are chosen freely, so that $|\acute\eta\Quat_3| = 8$.
Notice that the value of $\Re(z)$ matches Lemma \ref{lemma:realpartsigma}.
We have 
\begin{equation}
\label{equation:Lxyz}
\Lo_3^1 = \left\{ L = \begin{pmatrix}
1 & 0 & 0 \\ x & 1 & 0 \\ z & y & 1 \end{pmatrix}; x, y, z \in \RR
\right\}, \quad
\BL_\eta = \{ L \;|\; z \ne 0, z \ne xy\} \subset \Lo_3^1. 
\end{equation}
A computation yields
\begin{gather*}
\BL_{\frac{1-\hat a_1\hat a_2}{\sqrt{2}}} =
\{ L \;|\; z > \max\,\{0, xy\} \}, \quad
\BL_{\frac{1+\hat a_1\hat a_2}{\sqrt{2}}} =
\{ L \;|\; z < \min\,\{0, xy\} \}, \\
\BL_{\frac{\hat a_1 + \hat a_2}{\sqrt{2}}} =
\{ L \;|\; x > 0, 0 < z < xy \}, \quad
\BL_{\frac{\hat a_1 - \hat a_2}{\sqrt{2}}} =
\{ L \;|\; x > 0, xy < z < 0\}, \\
\BL_{\frac{-\hat a_1 - \hat a_2}{\sqrt{2}}} =
\{ L \;|\; x < 0, 0 < z < xy \}, \quad
\BL_{\frac{-\hat a_1 + \hat a_2}{\sqrt{2}}} =
\{ L \;|\; x < 0, xy < z < 0\}, \\
\BL_{\frac{-1+\hat a_1\hat a_2}{\sqrt{2}}} =
\BL_{\frac{-1-\hat a_1\hat a_2}{\sqrt{2}}} = \emptyset.
\end{gather*}
Notice that the six non empty sets are contractible.
\end{example}


Let $\lo_{n+1}^1$ be the Lie algebra of $\Lo_{n+1}^1$,
i.e., the set of lower strictly triangular matrices.
For $i \in \{1, \ldots, n\}$, let $\fl_i \in \lo_{n+1}^1$
be the matrix whose only non zero entry is $(\fl_i)_{(i+1,i)} = 1$.
Let $\lambda_i: \RR \to \Lo_{n+1}^1$ be the corresponding
one-parameter subgroup:
\begin{equation}
\label{equation:lambda}
\lambda_i(t) = \exp(t\fl_i) = I + t \fl_i \in \Lo_{n+1}^1.
\end{equation}
We saw in Equation \ref{equation:cE} how
the group $\cE_n = \{\pm 1\}^{\nmesmo}$
acts by automorphisms on $\Spin_{n+1}$.
The same group acts by automorphisms on $\Lo_{n+1}^{1}$:
\begin{equation}
\label{equation:cElambda}
(\lambda_i(t))^E = \lambda_i(E_i t).
\end{equation}
We therefore have $\bQ(L^E) = (\bQ(L))^E$
for all $L \in \Lo_{n+1}^1$ and $E \in \cE_n$.
Thus, the set $\cU_1$ is invariant:
$\cU_1^E = \cU_1$ for all $E \in \cE_n$.
We also have $\bL(z^E) = (\bL(z))^E$
for all $z \in \cU_1$ and $E \in \cE_n$.
For $z \in \tilde B_{n+1}^{+}$ and $E \in \cE_n$, we have
$(\BL_z)^E = \BL_{z^E}$:
in particular, the sets $\BL_z$ and $\BL_{z^E}$
are diffeomorphic through $L \mapsto L^E$.
Thus, in order to determine the homotopy type of $\BL_\sigma$
we decompose $\acute\sigma \Quat_{n+1}$ into $\cE_n$-orbits.
For each orbit, we take a representative $z$
and determine the homotopy type of $\BL_z$.

\begin{example}
\label{example:abacE}
For $\eta \in S_3$, the set $\acute\eta\Quat_{3}$
partitions into three $\cE_n$-orbits defined by $\Re(z)$.
From Example \ref{example:aba} it is clear that, indeed,
$\BL_{(1+\hat a_1\hat a_2)/\sqrt{2}}$ and
$\BL_{(1-\hat a_1\hat a_2)/\sqrt{2}}$ are diffeormorphic
by an action of $E$, i.e., by changes of signs.
Similarly, the four sets
$\BL_{(\pm \hat a_1 \pm \hat a_2)/\sqrt{2}}$
are diffeomorphic.
Finally, the two sets $\BL_{(-1\pm \hat a_1\hat a_2)/\sqrt{2}}$
are empty.
\end{example}

\smallskip

Let $\eta = a_{i_1}\cdots a_{i_m}$ be a reduced word
for the top permutation $\eta$.
A matrix $L \in \Lo_{n+1}^{1}$ is \textit{totally positive}
if and only if there exist positive numbers
$t_1, \ldots, t_m$ such that
\[ L = \lambda_{i_1}(t_1) \cdots \lambda_{i_m}(t_m). \]
The set $\Pos_\eta$ of totally positive matrices
is an open semigroup and a contractible connected component of $\BL_\eta$.
Furthermore, the set does not depend on the choice of the reduced word.
For $\sigma \in S_{n+1}$, set
$\Pos_\sigma = \BL_\sigma \cap \overline{\Pos_\eta}$.
This set can be equivalently defined as follows.
Let $\sigma = a_{i_1} \cdots a_{i_\ell}$ be a reduced word for $\sigma$:
we have $L \in \Pos_\sigma$
if and only if there exist positive numbers
$t_1, \ldots, t_\ell$ such that
\[ L = \lambda_{i_1}(t_1) \cdots \lambda_{i_\ell}(t_\ell). \]
The set $\Pos_\sigma$ is also
a contractible connected component of $\BL_\sigma$.

\begin{example}
\label{example:abapos}
We return to Example \ref{example:abaintro}.
For $n = 2$, it follows easily from Example \ref{example:aba} that
\[
\Pos_\eta = \BL_{\acute\eta},
\qquad
\acute\eta = \frac{\hat a_1 + \hat a_2}{\sqrt{2}}. \]
As we shall see,
for all $n$ and for all $\sigma \in S_{n+1}$,
we have $\Pos_\sigma \subseteq \BL_{\acute\sigma}$.
However, for larger values of $n$ and for most permutations $\sigma$,
$\Pos_\sigma$ is a small connected component
of a much larger set $\BL_{\acute\sigma}$.
\end{example}


\section{First examples of strata $\BLS_\varepsilon$}
\label{section:firstexamples}

Consider a permutation $\sigma$ and a fixed reduced word
$\sigma = a_{i_1}\cdots a_{i_\ell}$.
We define the open subset $\BLS_\varepsilon \subset \BL_\sigma$
if $\varepsilon$ is an ancestry of dimension $0$:
the general definition will be given later. 

There is exactly one preancestry of dimension $0$,
the empty preancestry,
as the first example in Figure \ref{fig:preancestry}. 
An ancestry $\varepsilon$ of dimension $0$
is a function (or sequence)
$\varepsilon: \llbracket \ell \rrbracket \to \{\pm 1\}$.
As a diagram, we mark intersections with black and white disks
(no squares), as in the first example in Figure \ref{fig:ancestry}
or in the boxed ancestries shown in Figure \ref{fig:563412}.
Recall that a black ball denotes $-1$ and a white ball denotes $+1$.

Given an ancestry $\varepsilon$ of dimension $0$, we have 
\begin{equation}
\label{equation:BLSdimzero}
\BLS_\varepsilon =
\{ \lambda_{i_1}(t_1) \cdots \lambda_{i_\ell}(t_\ell) \;|\;
t_k \in \RR \smallsetminus \{0\}, \; \sign(t_k) = \varepsilon_k \}
\subset B_\sigma.
\end{equation}
These subsets have been extensively studied
(see Section \ref{section:history} for references).
We quote the relevant results without proof:
not only are they in the literature
but also they will follow from our results for arbitrary ancestries.
It follows from Corollary 6.5 in \cite{Goulart-Saldanha0} that:
\begin{equation}
\label{equation:Coro65}
\BLS_\varepsilon \subseteq \BL_z, \qquad
z = P(\varepsilon) = (\acute a_{i_1})^{\varepsilon(1)} \cdots 
(\acute a_{i_\ell})^{\varepsilon(\ell)} \in \acute\sigma\Quat_{n+1}.
\end{equation}
The subsets $\BLS_\varepsilon \subset \BL_\sigma$ are open,
and the union over all ancestries of dimension $0$ is open and dense.

We shall define and study the sets $\BLS_\varepsilon$ for any ancestry
in Sections \ref{section:stratification} and \ref{section:goodstrata}.
The sets $\BLS_\varepsilon$ will be disjoint,
and $\BLS_\varepsilon$ will be a contractible
smooth submanifold
$\BLS_{\varepsilon} \subseteq \BL_{P(\varepsilon)}$
of codimension $d = \dim(\varepsilon)$.
Before presenting the general definition,
we discuss some examples and special cases.

Recall that an ancestry $\varepsilon$ (of dimension $0$) is \textit{thin}
if $i_{k_0} = i_{k_1}$ implies $\varepsilon_{k_0} = \varepsilon_{k_1}$.
If $\varepsilon$ is a thin ancestry,
the subset $\BLS_\varepsilon \subset \BL_\sigma$
is also called \textit{thin}.
Notice that $\varepsilon = (+,+,\ldots, +,+)$ is thin,
with $P(\varepsilon) = \acute\sigma$ and
$\BLS_\varepsilon = \Pos_\sigma \subseteq \BL_{\acute\sigma}$,
a contractible connected component.
More generally,
if $\varepsilon$ is a thin ancestry then
there exists $E \in \cE_n$ such that
$\varepsilon(k) = (\acute a_{i_k})^E$ (for all $k$).
We therefore have
$P(\varepsilon) = \acute\sigma^E$ and
$\BLS_\varepsilon = (\Pos_\sigma)^E$.
It follows that
$\BLS_\varepsilon \subseteq \BL_{\acute\sigma^E}$
is a contractible connected component.
The set
\begin{equation}
\label{equation:thickpart}
\BL_{z,\thick} = \BL_z \smallsetminus 
\bigcup_{\varepsilon \textrm{ thin}} \BLS_\varepsilon 
\end{equation}
is called the \textit{thick part} of $\BL_z$.
As we shall see, in general $\BL_{z,\thick}$ can be empty or disconnected.

\begin{example}
\label{example:abazero}
As in Examples \ref{example:abaintro} and \ref{example:aba},
take $n = 2$ and $\eta \in S_3$.
There are two reduced words:
$\eta = a_1a_2a_1 = a_2a_1a_2$,
shown as diagrams in Figure \ref{fig:aba-bab}.
For each reduced word, there are two preancestries,
with dimensions $0$ and $1$,
also shows in Figure \ref{fig:aba-bab}.

\begin{figure}[ht]
\begin{center}
\includegraphics[scale=0.25]{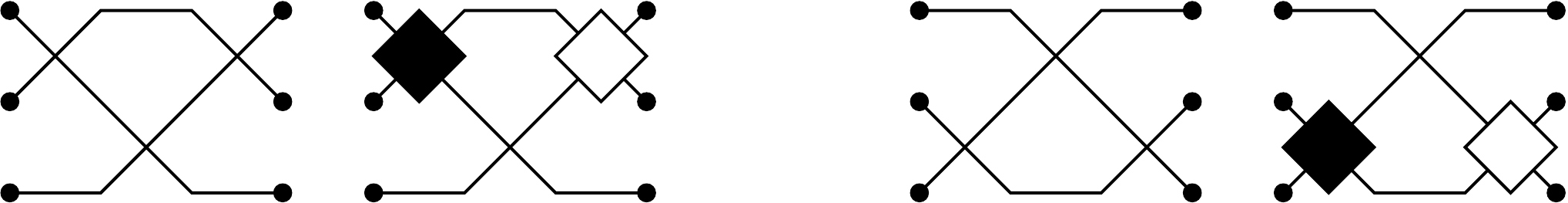}
\end{center}
\caption{Two reduced words for 
$\eta = a_1a_2a_1 = a_2a_1a_2 \in S_3$.
For each word, there are two preancestries.}
\label{fig:aba-bab}
\end{figure}

We first consider the reduced word $\eta = a_1a_2a_1$.
Write $L$ as in Equation \eqref{equation:Lxyz}.
The eight ancestries of dimension $0$ are $(\pm 1, \pm 1, \pm 1)$;
the two ancestries of dimension $1$ are $(-2, \pm 1, +2)$.
A simple computation gives
\[
\lambda_1(t_1)\lambda_2(t_2)\lambda_1(t_3) = \begin{pmatrix}
1 & 0 & 0 \\ t_1+t_3 & 1 & 0 \\ t_2t_3 & t_2 & 1 \end{pmatrix}
\]
and therefore
\begin{align*}
\BL_{(-1,+1,+1)} &= \{ L \;|\; z > \max\{0,xy\}, y > 0\}, \\
\BL_{(+1,-1,-1)} &= \{ L \;|\; z > \max\{0,xy\}, y < 0\}.
\end{align*}
Take $z_0 = (1-\hat a_1\hat a_2)/\sqrt{2}$.
Notice that $P(-1,+1,+1) = P(+1,-1,-1) = P(-2,+1,+2) = z_0$.
These are the only ancestries $\varepsilon$ with $P(\varepsilon) = z_0$,
as can be verified either by computing $P$ for the other ancestries
or from Lemmas \ref{lemma:preminus} and \ref{lemma:preplus}.
As we shall see  later,
\begin{align*}
\BL_{(-2,+1,+2)} &= \{L \;|\; y = 0, z > 0\}, \\
\BL_{z_0} &=
\BL_{(-1,+1,+1)} \sqcup \BL_{(-2,+1,+2)} \sqcup \BL_{(+1,-1,-1)}, 
\end{align*}
consistently with Equation \eqref{equation:Coro65}.
The subset $\BL_{(-2,+1,+2)} \subset \BL_{z_0}$
is a contractible submanifold of codimension $d = 1$.
The CW complex in Figure \ref{fig:aba-CW} is homotopically equivalent
to $\BL_{z_0}$.

A similar decomposition holds for $\BL_{(1+\hat a_1\hat a_2)/\sqrt{2}}$.
The four connected components $\BL_{(\pm \hat a_1 \pm \hat a_2)/\sqrt{2}}$
are thin,  consistently with Example \ref{example:abapos}.

Recall that the permutation $\eta$
also admits the reduced word $\eta = a_2a_1a_2$. 
A similar decomposition exists for the other reduced word,
but the strata are different.
For this other word, the set
$\BL_{(1-\hat a_1\hat a_2)/\sqrt{2}}$
contains two open strata: 
\begin{align*}
\BL_{(+1,+1,-1)} &= \{ L \;|\; z > \max\{0,xy\}, x > 0\}, \\
\BL_{(-1,-1,+1)} &= \{ L \;|\; z > \max\{0,xy\}, x < 0\}
\end{align*}
and a third stratum of codimension $1$:
\[ \BL_{(-2,-1,+2)} = \{ L \;|\; x = 0, z > 0\}. \]
Clearly, the meaning of the ancestry depends
on the reduced word being used.
\end{example}

\begin{example}
\label{example:bacb}
Take $n = 3$ and $\sigma = a_2a_1a_3a_2 = (3,4,1,2) \in S_4$.
Up to transposing adjacent commuting generators,
as in $a_2a_1a_3a_2 = a_2a_3a_1a_2$,
the permutation $\sigma$ admits only one reduced word,
shown in Figure \ref{fig:bacb}.

\begin{figure}[ht]
\begin{center}
\includegraphics[scale=0.25]{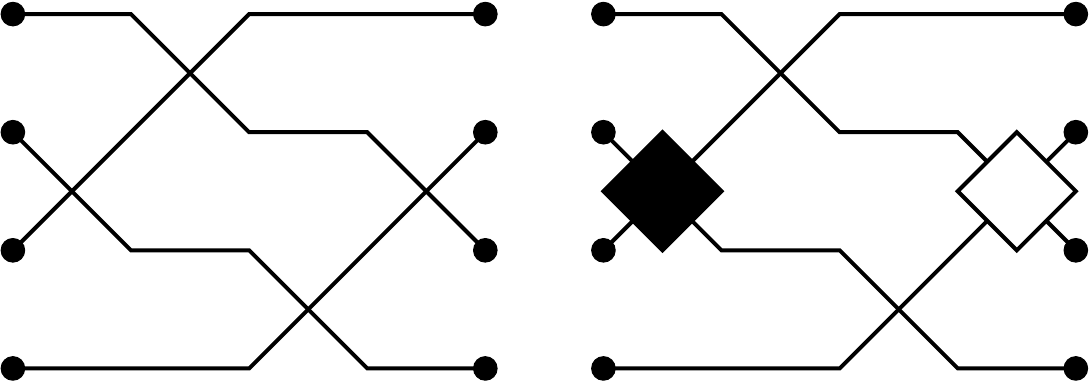}
\end{center}
\caption{The permutation $\sigma = a_2a_1a_3a_2 \in S_4$
and its two preancestries.}
\label{fig:bacb}
\end{figure}

There are two preancestries and $20$ ancestries:
\[ (\pm 1, \pm 1, \pm 1, \pm 1), \quad d = 0; \qquad
(-2, \pm 1, \pm 1, +2), \quad d = 1. \]
For the preancestry of dimension $d = 1$, 
$(\rho_k)$ (as in Equations \eqref{equation:preancestry0} and
\eqref{equation:preancestry})
satisfies $\rho_1 = \rho_2 = \rho_3 = a_1a_2a_3a_2a_1$.
Take $z_0 = -\acute\sigma\hat a_3 =
(1-\hat a_1\hat a_2 - \hat a_1\hat a_3 - \hat a_2\hat a_3)/2$.
The three ancestries with $P(\varepsilon) = z_0$
are $(+1,+1,-1,-1)$, $(-2,-1,+1,+2)$ and $(-1,-1,+1,+1)$.
At this point, the reader may guess that $\BL_{z_0}$
is homotopically equivalent to the CW complex in Figure \ref{fig:bacb-CW},
which is indeed correct. 

\begin{figure}[ht]
\begin{center}
\includegraphics[scale=0.25]{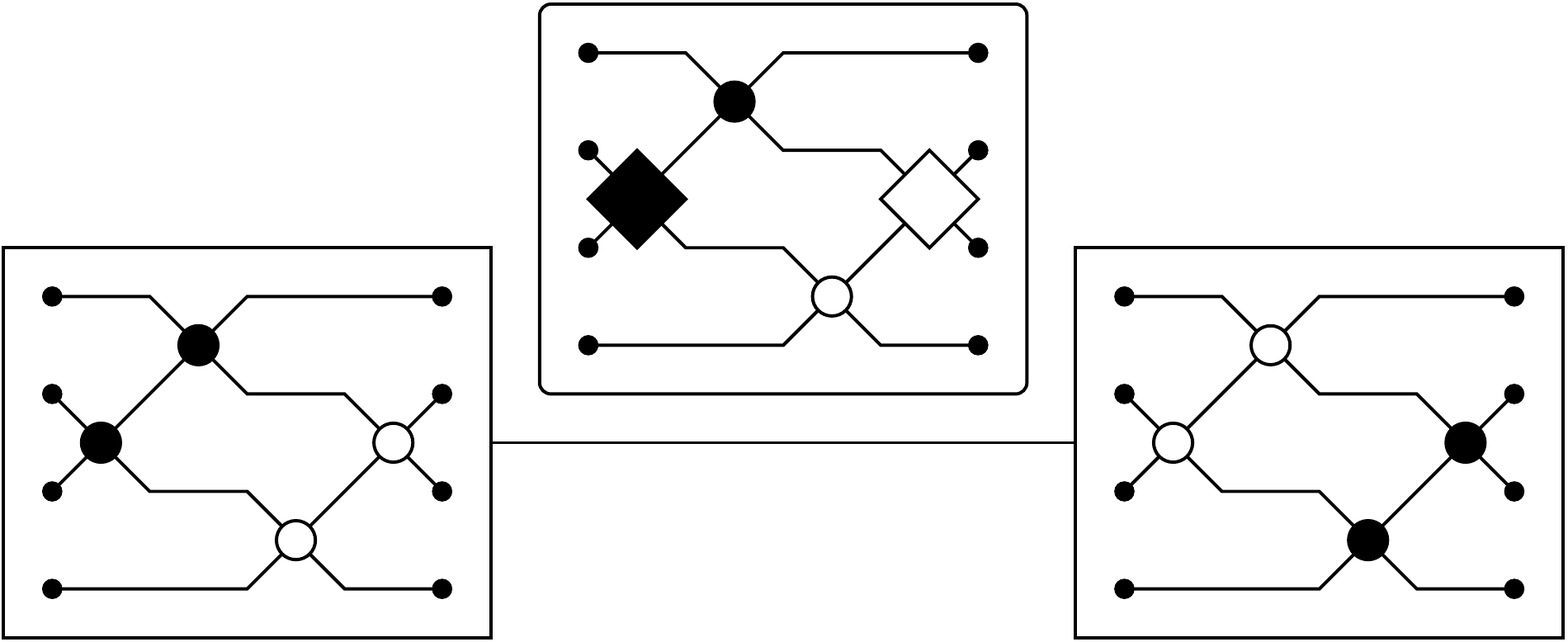}
\end{center}
\caption{A CW complex homotopically equivalent to $\BL_{z_0}$.}
\label{fig:bacb-CW}
\end{figure}


Write $L \in \Lo_4^1$ as
\[ \Lo_4^1 = \left\{ L = \begin{pmatrix}
1 & 0 & 0 & 0 \\
x & 1 & 0 & 0 \\
u & y & 1 & 0 \\
w & v & z & 1 \end{pmatrix}, u, v, w, x, y, z \in \RR \right\}. \]
By applying the definition, the set $\BL_{\sigma}$ is
\[ \BL_{\sigma} = \{ L \;|\; w = 0, u \ne 0, v \ne 0, xyz = xv + zu \}. \]
If $L = \lambda_2(t_1)\lambda_1(t_2)\lambda_3(t_3)\lambda_2(t_4)$ then
\[ u = t_1t_2, \quad v = t_3t_4, \quad w = 0, \quad
x = t_2, \quad y = t_1+t_4, \quad z = t_3. \]
Let $U \subset (0,+\infty)^2 \times \RR^2$
(with coordinates $(u,v,x,y)$) be
the contractible open set defined by $xy < u$.
Consider the map $\Phi: U \to \BL_{\sigma}$ taking $(u,v,x,y)$
to the matrix $L \in \BL_{\sigma}$ with the prescribed values
of $u,v,x,y$ and $z = xv/(xy-u)$.

If $x > 0$ we have $L = \Phi(u,v,x,y) \in \BL_{(+1,+1,-1,-1)}$;
if $x < 0$ we have $L = \Phi(u,v,x,y) \in \BL_{(-1,-1,+1,+1)}$.
We will see that
\begin{align*}
\BL_{(-2,-1,+1,+2)} &= \Phi[(0,+\infty)^2 \times \{0\} \times \RR], \\
\BL_{z_0} &=
\BL_{(-1,-1,+1,+1)} \sqcup \BL_{(-2,-1,+1,+2)} \sqcup \BL_{(-1,-1,+1,+1)}.
\end{align*}
It follows that $\BL_{z_0}$ is contractible.
Similar computations verify that $\BL_{\sigma}$ has $12$ connected components, $8$ thin and $4$ thick, all contractible. 
\end{example}

We conclude this section with an effective enumeration
of the connected components of $\BL_\eta$ for $n \ge 5$.
Recall that it is known that there are $3 \cdot 2^n$ connected components
(see Section \ref{section:history} and \cite{SSV2}).

\begin{prop}
\label{prop:effective}
For $n \ge 5$, the $3 \cdot 2^n$ connected components of $\BL_\eta$ are
\[ \Pos_\eta^E, E \in \cE_n, \qquad
\BL_{z,\thick}, z \in \acute\eta\Quat_{n+1}. \]
The first list are the $2^n$ thin connected components;
the second are the $2^{n+1}$ thick connected components.
\end{prop}

\begin{proof}
We must prove that each $\BL_{z,\thick}$ is connected and non empty.
It follows from Remark \ref{remark:thicketa} that
$\BL_{z,\thick} \ne \emptyset$ for each $z \in \acute\eta\Quat_{n+1}$.
If one of them were to be disconnected, we would have more than
$3\cdot 2^n$ connected components,
contradicting the number of connected components.
\end{proof}


\section{The stratification $\BLS_\varepsilon$}
\label{section:stratification}

Consider as usual a fixed reduced word
$\sigma = a_{i_1}\cdots a_{i_\ell}$.
For $L \in \BL_\sigma$,
we show how to determine the ancestry of $L$.

Let $\tilde z_\ell = \bQ(L)$.
Take $q_\ell \in \Quat_{n+1}^{1}$ such that
$z_\ell = \tilde z_\ell q_\ell \in \Bru_{\acute\sigma}$.
Define recursively $\sigma_0 = 1$, $\sigma_1 = a_{i_1}$,
$\sigma_k = \sigma_{k-1} a_{i_k} = a_{i_1}\cdots a_{i_k}$
so that $\sigma = \sigma_{\ell}$.
From Theorem 1 from \cite{Goulart-Saldanha0},
we have well-defined sequences
$(\theta_k)_{0 < k \le \ell}$ and
$(z_k)_{0 \le k \le \ell}$
such that $z_0 = 1 \in \Spin_{n+1}$ and
\begin{equation}
\label{equation:zk}
z_k = z_{k-1} \alpha_{i_k}(\theta_k) \in \Bru_{\acute\sigma_k},
\qquad
\theta_k \in (0,\pi). 
\end{equation}
Take $\varrho_k \in \tilde B_{n+1}^{+}$ such that
$z_k \in \grave\eta \Bru_{\varrho_k}$:
the sequence $(\varrho_k)$ is the desired ancestry.
The corresponding preancestry is $(\rho_k)$,
$\rho_k = \Pi_{\tilde B_{n+1}^{+}, S_{n+1}}(\varrho_k)$
so that $z_k \in \grave\eta \Bru_{\rho_k}$.
Given an ancestry $\varepsilon$,
let $\BLS_{\varepsilon} \subset \BL_\sigma$
be the set of matrices $L$ with ancestry $\varepsilon$.

We verify that the sequences $(\rho_k)$ and $(\varrho_k)$
are indeed a preancestry and an ancestry
in the sense defined in
Sections \ref{section:preancestry} and \ref{section:ancestry}.
Indeed, we have $\rho_0 = \eta$ and $\varrho_0 = \acute\eta$ since
$z_0 = 1 \in \grave\eta \Bru_{\acute\eta} \subset \grave\eta \Bru_{\eta}$.
If $z_{k-1} \in \grave\eta \Bru_{\varrho_{k-1}}$
and $\rho_{k-1} < \rho_{k-1} a_{i_k}$
then $z_{k-1} \alpha_{i_k}(\theta) \in
\grave\eta \Bru_{\varrho_{k-1} \acute a_{i_k}}$
for all $\theta \in (0,\pi)$:
this implies $\varrho_k = \varrho_{k-1} \acute a_{i_k}$
and $\rho_k = \rho_{k-1} a_{i_k}$.
If $z_{k-1} \in \grave\eta \Bru_{\varrho_{k-1}}$
and $\rho_{k-1} > \rho_{k-1} a_{i_k}$
then $z_{k-1} \alpha_{i_k}(\theta)$ (for $\theta \in (0,\pi)$)
belongs to one of the following three sets:
$\grave\eta \Bru_{\varrho_{k-1}}$,
$\grave\eta \Bru_{\varrho_{k-1} \acute a_{i_k}}$,
$\grave\eta \Bru_{\varrho_{k-1} \hat a_{i_k}}$.
This implies $\varrho_k$ to be one among
$\varrho_{k-1}$, ${\varrho_{k-1} \acute a_{i_k}}$,
${\varrho_{k-1} \hat a_{i_k}}$.
Finally, $\tilde z_k \in \cU_1$ implies
$z_k \in \Bru_{\eta}$ and $\varrho_\ell \in \acute\eta \Quat_{n+1}$.
We thus have
\begin{equation}
\label{equation:stratification}
\BL_\sigma = \bigsqcup_{\varepsilon} \BLS_{\varepsilon},
\qquad
\BL_z = \bigsqcup_{P(\varepsilon) = z} \BLS_{\varepsilon},
\end{equation}
where $\varepsilon$ varies over ancestries
(for a fixed reduced word for $\sigma$).
We are ready to complete the proof of Theorem \ref{theo:two}.

\begin{proof}[Proof of Theorem \ref{theo:two}]
It follows from the definition of $\BLS_{\varepsilon}$ above and
from Equation \eqref{equation:stratification}
that $\BLS_{\varepsilon} \subseteq \BL_{P(\varepsilon)}$.
For any preancestry $\varepsilon_0$ and any $z \in \acute\sigma\Quat_{n+1}$,
we have $\NL_{\varepsilon_0}(z) = N_{\varepsilon_0}(z)$;
see Remark \ref{remark:NLN}.
Theorem \ref{theo:two} now follows directly from
Lemmas \ref{lemma:preminus} and \ref{lemma:preplus}.
\end{proof}

We present the interpretation of $\xi(k)$.
Consider $z_{k-1} \in \grave\eta \Bru_{\varrho_{k-1}}$.
If $\rho_{k-1} < \rho_{k-1} a_{i_k}$ then
$z_{k-1} \alpha_{i_k}(\theta) \in
\grave\eta\Bru_{\varrho_{k-1} \acute a_{i_k}}$
for all $\theta \in (0,\pi)$:
in this case we take $\xi(k) = 1$.
If $\rho_{k-1} > \rho_{k-1} a_{i_k}$ then
there exists a unique $\theta_{\bullet} \in (0,\pi)$ such that 
$z_{k-1} \alpha_{i_k}(\theta_{\bullet}) \in
\grave\eta\Bru_{\varrho_{k-1} \acute a_{i_k}}$.
If $\theta_k < \theta_{\bullet}$ we have
$z_k \in \grave\eta\Bru_{\varrho_{k-1}}$,
$\varrho_k = \varrho_{k-1}$ and $\xi(k) = 0$.
If $\theta_k > \theta_{\bullet}$ we have
$z_k \in \grave\eta\Bru_{\varrho_{k-1} \hat a_{i_k}}$,
$\varrho_k = \varrho_{k-1} \hat a_{i_k}$ and $\xi(k) = 2$.
Finally, 
if $\theta_k = \theta_{\bullet}$ we have
$z_k \in \grave\eta\Bru_{\varrho_{k-1} \acute a_{i_k}}$,
$\varrho_k = \varrho_{k-1} \acute a_{i_k}$ and $\xi(k) = 1$.
Thus, $\xi(k)$ tells us about the size of $\theta_k$:
$\xi(k) = 0$ means that $\theta_k$ is small,
$\xi(k) = 2$ means that $\theta_k$ is large
and
$\xi(k) = 1$ means that $\theta_k$ is just right.

Let us introduce some more notation.  Let
\[ \cU^{\diamond}_1 =
\bigsqcup_{\sigma \in S_{n+1}} \grave\eta \Bru_{\acute\sigma}, \qquad
\cU_1 \subset \cU^{\diamond}_1 \subset \overline{\cU_1} \subset \Spin_{n+1}. \]
The set $\cU^{\diamond}_1$ is a fundamental domain
for the action of $\Quat_{n+1}$ on $\Spin_{n+1}$:
given $z \in \Spin_{n+1}$ there exists a unique $q \in \Quat_{n+1}$
such that $zq \in \cU^{\diamond}_1$ (see \cite{Goulart-Saldanha0}).
For each $k$, write $z_k = \tilde z_k q_k$,
$\tilde z_k \in \cU^{\diamond}_1$, $q_k \in \Quat_{n+1}$.
We clearly have $\tilde z_k \in \grave\eta \Bru_{\acute\rho_k}$.

\begin{lemma}
\label{lemma:tildetheta}
There exist unique $\tilde\theta_k \in (-\pi,0)\cup(0,\pi)$
such that $\tilde z_k = \tilde z_{k-1} \alpha_{i_k}(\tilde\theta_k)$.
Furthermore, for $s = [a_{i_k},q_{k-1}] \in \{\pm 1\}$
we have either $\tilde\theta_k = s\theta_k$ or 
$\tilde\theta_k = s(\theta_k-\pi)$. 
In the first case we have $q_k = q_{k-1}$;
in the second case, $q_k = q_{k-1} \hat a_{i_k}$.
\end{lemma}

\begin{proof}
If $z \in \cU^{\diamond}_1$ and $\theta \in (0,\pi)$
then either $z \alpha_i(\theta) \in \cU^{\diamond}_1$
or $z \alpha_i(\theta-\pi) \in \cU^{\diamond}_1$ (but not both).
Similarly,
either $z \alpha_i(-\theta) \in \cU^{\diamond}_1$
or $z \alpha_i(-\theta+\pi) \in \cU^{\diamond}_1$ (but not both).

We have $z_k = z_{k-1} \alpha_{i_k}(\theta_k)$ and
$z_{k-1} = \tilde z_{k-1} q_{k-1}$.
We therefore have
$z_k = \tilde z_{k-1} q_{k-1} \alpha_{i_k}(\theta_k)$
and therefore
$z_k = \tilde z_{k-1} \alpha_{i_k}(s\theta_k) q_{k-1}$.
By definition, $q_k$ is such that
$z_k = \tilde z_k q_k$, $\tilde z_k \in \cU^{\diamond}_1$.
If $\tilde z_{k-1} \alpha_{i_k}(s\theta_k) \in \cU^{\diamond}_1$
we take $\tilde\theta_k = s\theta_k$,
$\tilde z_k = \tilde z_{k-1} \alpha_{i_k}(s\theta_k)$
and $q_k = q_{k-1}$.
If $s = +1$ and
$\tilde z_{k-1} \alpha_{i_k}(\theta_k) \notin \cU^{\diamond}_1$
we take
$\tilde\theta_k = \theta_k - \pi$,
$\tilde z_k = \tilde z_{k-1} \alpha_{i_k}(\theta_k-\pi) \in \cU^{\diamond}_1$
and $q_k = q_{k-1} \hat a_{i_k} = \hat a_{i_k} q_{k-1}$.
Finally, 
If $s = -1$ and
$\tilde z_{k-1} \alpha_{i_k}(-\theta_k) \notin \cU^{\diamond}_1$
we take
$\tilde\theta_k = -\theta_k + \pi$,
$\tilde z_k = \tilde z_{k-1} \alpha_{i_k}(-\theta_k+\pi) \in \cU^{\diamond}_1$
and $q_k = q_{k-1} \hat a_{i_k} = -\hat a_{i_k} q_{k-1}$.
\end{proof}

We are finally ready to give the interpretation of $\varepsilon(k)$.

\begin{lemma}
\label{lemma:whoisvarepsilon}
We have $\sign(\varepsilon(k)) = \sign(\tilde\theta_k)$.
Also,
$\varepsilon(k) = -2$ if and only if $\rho_k < \rho_{k-1}$;
$\varepsilon(k) = +2$ if and only if $\rho_k > \rho_{k-1}$.
\end{lemma}

\begin{proof}
Compare the formula for $\tilde\theta_k$ in Lemma \ref{lemma:tildetheta}
with Equation \eqref{equation:xi2varepsilon}.
\end{proof}

The above construction can be thought of as a way to extend
the construction in Section \ref{section:firstexamples},
based on the functions $\lambda_i$,
to infinity in the cases where it is impossible in $\Lo_{n+1}^{1}$.
The group $\Lo_{n+1}^{1}$ is unsuitable:
we work instead in the compact group $\Spin_{n+1}$ (or $\SO_{n+1}$),
using the functions $\alpha_i$ instead of $\lambda_i$
and performing the required adaptations.
The examples should clarify these remarks.

\begin{example}
\label{example:abaL}
As in Examples \ref{example:abaintro} and \ref{example:aba},
set $n = 2$ and $\sigma = \eta = a_1a_2a_1$. Consider
\[ L_0 = \begin{pmatrix} 1 & 0 & 0 \\ 0 & 1 & 0 \\ 1 & 0 & 1 \end{pmatrix},
\qquad \tilde z_3 = \bQ(L_0) = \begin{pmatrix}
\sqrt{2}/2 & 0 & -\sqrt{2}/2 \\ 0 & 1 & 0 \\ \sqrt{2}/2 & 0 & \sqrt{2}/2
\end{pmatrix}. \] 
(More correctly, the matrix shown is $\Pi(\tilde z_3) \in \SO_3$.
Here and in other occasions it is easier to do computations
in $\SO_{n+1}$ instead of $\Spin_{n+1}$.)

We have
$\tilde z_3 =
\alpha_1(-\frac{\pi}{2})\alpha_2(\frac{\pi}{4})\alpha_1(\frac{\pi}{2})$.
We have
$\rho_1 = \rho_2 = a_1a_2$ and therefore
(denoting ancestries by $\varepsilon$)
$L_0 \in \BL_{(-2,1,2)}$.
More generally,
it is not hard to verify that, for $L \in \Lo_3^1$,
we have $L \in \BL_{(-2,1,2)}$
if and only if $y = 0$ and $z > 0$.
This justifies the claims made in Example \ref{example:abazero},
particularly Figure \ref{fig:aba-CW}.
\end{example}

\begin{example}
\label{example:bacbL}
As in Example \ref{example:bacb},
set $n = 3$ and $\sigma = a_2a_1a_3a_2$.
Consider
\[ L_0 = \begin{pmatrix} 1 & 0 & 0 & 0 \\ 0 & 1 & 0 & 0 \\
1 & 0 & 1 & 0 \\ 0 & 1 & 0 & 1 \end{pmatrix},
\qquad \tilde z_4 = \bQ(L_0) = \frac{\sqrt{2}}{2} \begin{pmatrix}
1 & 0 & -1 & 0 \\
0 & 1 & 0 & -1 \\
1 & 0 & 1 & 0 \\
0 & 1 & 0 & 1 
\end{pmatrix}. \] 
We have
$\tilde z_4 = \alpha_2(-\frac{\pi}{2})\alpha_1(-\frac{\pi}{4})
\alpha_3(\frac{\pi}{4})\alpha_2(\frac{\pi}{2})$,
$\rho_1 = \rho_2 = \rho_3 = a_1a_2a_3a_2a_1$
and therefore
(denoting ancestries by $\varepsilon$)
$L_0 \in \BL_{(-2,-1,1,2)}$.
The claims in Example \ref{example:bacb} are now easy to verify.
\end{example}


\section{The strata $\BLS_\varepsilon$}
\label{section:goodstrata}

We now prove that the strata $\BLS_{\varepsilon}$
are reasonably well-behaved.

\begin{lemma}
\label{lemma:submanifold}
Consider a permutation and reduced word
$\sigma = a_{i_1}\cdots a_{i_\ell} \in S_{n+1}$
and an ancestry $\varepsilon$.
The subset $\BLS_{\varepsilon} \subseteq \BL_{\sigma}$ is a smooth submanifold
of codimension $d = \dim(\varepsilon)$.
\end{lemma}

\begin{proof}
Let $(\varrho_k)_{0 \le k \le \ell}$ be the ancestry
encoded by $\varepsilon$.
Consider $L \in \BLS_{\varepsilon}$ and $\tilde z_\ell = \bQ(L)$.
Construct $(z_k)_{0 \le k \le \ell}$ and $(\theta_k)_{0 < k \le \ell}$
as above, so that $z_\ell = \tilde z_\ell q_\ell$,
$z_k = z_{k-1} \alpha_{i_k}(\theta_k) \in \grave\eta \Bru_{\varrho_k}$.
Let $K_0 = \{ k \in \ZZ, 1 \le k \le \ell, \varrho_k < \varrho_{k-1} \}$
so that $|K_0| = d$.
For $0 \le j \le \ell$,
let $V_j = \{x \in \RR^j \;|\; k \in K_0 \to x_k = 0 \} \subset \RR^j$,
a linear subspace;
$V_\ell$ has codimension $d$.
We construct a compact neighborhood $U \subset \RR^\ell$ of $0$
and a smooth local diffeomorphism $\Phi: U \to \Bru_{\acute\sigma}$
such that $\Phi(0) = z_\ell$ and
$\Phi(x) \in \bQ[\BLS_{\varepsilon}] q_\ell$ if and only if $x \in V_\ell$.

The set $U$ has the form $U = U_\ell$ where
\[ U_k = [-\epsilon_1,\epsilon_1] \times \cdots
\times [-\epsilon_k,\epsilon_k]. \]
We recursively define $\epsilon_k > 0$ and 
maps $\Phi_k: U_k \to \Bru_{\acute\sigma_k}$
with $\Phi_k(0) = z_k$.
In every case we shall have
$\Phi_k(x_{k-1},\theta) =
\Phi_{k-1}(x_{k-1}) \alpha_{i_k}(\ast)$,
where $\ast$ stands for a smooth function of
$x_{k-1} \in U_{k-1}$ and $\theta \in [-\epsilon_k,\epsilon_k]$.
The case $k = 0$ is trivial.

If $k \notin K_0$, take $\epsilon_k > 0$ sufficiently small
such that the following two conditions hold.
For all $x_{k-1} \in U_{k-1}$ and
$\theta \in [-\epsilon_k,\epsilon_k]$, we have 
$\Phi_{k-1}(x_{k-1})\alpha_{i_k}(\theta_k + \theta) \in
\Bru_{\acute\sigma_k}$.
For all $x_{k-1} \in U_{k-1} \cap V_{k-1}$ and
$\theta \in [-\epsilon_k,\epsilon_k]$, we have 
$\Phi_{k-1}(x_{k-1})\alpha_{i_k}(\theta_k + \theta) \in
\grave\eta\Bru_{\varrho_k}$.
The existence of such $\epsilon_k > 0$ 
follows from Theorem 1 of \cite{Goulart-Saldanha0}.
We then define $\Phi_k(x_{k-1},\theta) = 
\Phi_{k-1}(x_{k-1})\alpha_{i_k}(\theta_k + \theta)$.

If $k \in K_0$,
there exists a smooth function
$\vartheta: U_{k-1} \cap V_{k-1} \to (0,\pi)$
such that, for all $x_{k-1} \in U_{k-1} \cap V_{k-1}$, we have
$\acute\eta \Phi_{k-1}(x_{k-1}) \alpha_{i_k}(\theta_k + \vartheta(x_{k-1}))
\in \Bru_{\varrho_k}$
(here we again use Theorem 1 of \cite{Goulart-Saldanha0}).
Notice that $\vartheta(0) = 0$.
Let $\Pi: U_{k-1} \to U_{k-1} \cap V_{k-1}$ be the orthogonal projection.
Extend $\vartheta$ to $U_{k-1}$ by defining
$\vartheta(x_{k-1}) = \vartheta(\Pi(x_{k-1}))$;
notice that this is a smooth function.
Define
\[ \Phi_k(x_{k-1},\theta) =
\Phi_{k-1}(x_{k-1}) \alpha_{i_k}(\theta_k + \vartheta(x_{k-1}) + \theta). \]
Notice that, for $x_{k-1} \in U_{k-1} \cap V_{k-1}$, we have 
$\Phi_k(x_{k-1},\theta) \in
\grave\eta \Bru_{\varrho_k(\acute a_{i_k})^{\sign(\theta)}}$
(with $\sign(0) = 0$).
Choose sufficiently small $\epsilon_k > 0$ and we are done.
\end{proof}

\begin{lemma}
\label{lemma:contractible}
Consider a permutation $\sigma$,
a reduced word $\sigma = a_{i_1}\cdots a_{i_\ell}$
and an ancestry $\varepsilon$ with $d = \dim(\varepsilon)$.
The smooth submanifold $\BLS_{\varepsilon} \subset \BL_\sigma$
is diffeomorphic to $\RR^{\ell-d}$.
\end{lemma}

\begin{proof}
Let $(\varrho_k)$ be the sequence of elements of $\tilde B_{n+1}^{+}$
defining the ancestry $\varepsilon$; let $\xi$ be as usual.
Let $\Psi_k: (0,\pi)^k \to \Bru_{\acute\sigma_k}$
be the diffeomorphism
\[ \Psi_k(\theta_1,\ldots,\theta_k) =
\alpha_{i_1}(\theta_1)\cdots\alpha_{i_k}(\theta_k). \]
For $0 \le k \le \ell$, let $X_k \subseteq (0,\pi)^k$ be defined by:
\[ (\theta_1,\ldots,\theta_k) \in  X_k \quad\iff\quad
\forall j, (0 \le j \le k) \to
(\Psi_j(\theta_1, \ldots, \theta_j) \in \grave\eta \Bru_{\varrho_j}). \]
As in the proof of Lemma \ref{lemma:submanifold},
let $K_0 = \{ k \in \ZZ, 1 \le k \le \ell, \varrho_k < \varrho_{k-1} \}$.
We prove by induction on $k$ that each set $X_k$ is diffeomorphic to
$\RR^{k - |K_0 \cap [1,k]|}$.
By definition, the restriction
$\Psi_\ell: X_\ell \to \BLS_{\varepsilon}$ is a bijection;
it follows from Lemma \ref{lemma:submanifold}
that it is a diffeomorphism.

The set $X_0$ has a single element, the empty sequence.
We have $X_k \subseteq X_{k-1} \times (0,\pi)$.
For $(\theta_1,\ldots,\theta_k) \in  X_{k-1} \times (0,\pi)$, 
we have $(\theta_1,\ldots,\theta_k) \in X_k$ if and only if
$\Psi_k(\theta_1,\ldots,\theta_k) \in \grave\eta\Bru_{\varrho_k}$.
We again divide our discussion into cases.
If $\xi(k) = 1$ and
$\varrho_{k-1} < \varrho_k = \varrho_{k-1} \acute a_{i_k}$ then
$X_k = X_{k-1} \times (0,\pi)$
and we are done.

Otherwise, we have $\rho_{k-1} a_{i_k} < \rho_{k-1}$.
There exists a smooth function
$\vartheta: X_{k-1} \to (0,\pi)$ such that,
for all $\bar\theta = (\theta_1, \ldots, \theta_{k-1}) \in X_{k-1}$, we have 
\[ \Psi_{k-1}(\bar\theta)\alpha_{i_k}(\vartheta(\bar\theta)) \in
\Bru_{\varrho_{k-1}\acute a_{i_k}}. \]
Indeed, in the notation of Remark 6.6 of \cite{Goulart-Saldanha0},
we have $\vartheta(\bar\theta) =
\pi - \Theta_{i_k}(\Psi_{k-1}(\bar\theta))$.
If $\xi(0) = 0$ we have
\[ X_k = \{(\bar\theta,\theta_k) \in X_{k-1}\times (0,\pi) \;|\;
\theta_k < \vartheta(\bar\theta) \}, \]
which is diffeomorphic to $X_{k-1} \times (0,1)$
and therefore to an open ball.
If $\xi(0) = 2$ we have
\[ X_k = \{(\bar\theta,\theta_k) \in X_{k-1}\times (0,\pi) \;|\;
\theta_k > \vartheta(\bar\theta) \}, \]
and we are done.
Finally, if $\xi(0) = 1$ we have
\[ X_k = \{(\bar\theta,\theta_k) \in X_{k-1}\times (0,\pi) \;|\;
\theta_k = \vartheta(\bar\theta) \}, \]
which is diffeomorphic to $X_{k-1}$.
This completes the proof.
\end{proof}

Recall that in Equation \ref{equation:poset} we define
a partial order on the set of ancestries 
(for a fixed reduced word).
The following lemma is the reason why we defined a partial order
among ancestries.

\begin{lemma}
\label{lemma:poset}
Let $\varepsilon, \tilde\varepsilon$ be ancestries.
If $\BLS_{\varepsilon} \cap \overline{\BLS_{\tilde\varepsilon}} \ne \emptyset$
then $\varepsilon \preceq \tilde\varepsilon$.
\end{lemma}

\begin{proof}
Assume that 
$\BLS_{\varepsilon} \cap \overline{\BLS_{\tilde\varepsilon}} \ne \emptyset$.
Thus, there exists a sequence $(L_j)$ of elements of $\BLS_{\tilde\varepsilon}$
converging to an element $L_\infty$ of $\BLS_{\varepsilon}$.
This implies that $P(\tilde\varepsilon) = P(\varepsilon)$.
If $\tilde z_{j,\ell} = \bQ(L_j)$ and $\tilde z_{\infty,\ell} = \bQ(L_\infty)$,
we have $\lim_j \tilde z_{j,\ell} = \tilde z_{\infty,\ell}$.
This implies that $P(\tilde\varepsilon) = P(\varepsilon)$.
Thus, for $z_{j,\ell}, z_{\infty,\ell} \in \Bru_{\acute\sigma}$
defined as usual we have
$\lim_j z_{j,\ell} = z_{\infty,\ell}$.
Since $z_k$ is a smooth function of $z_\ell$, we have
$\lim_j z_{j,k} = z_{\infty,k}$.
We have $z_{j,k} \in \grave\eta\Bru_{\tilde\varrho_k}$
and $z_{\infty,k} \in \grave\eta\Bru_{\varrho_k}$
and therefore $\Bru_{\varepsilon} \cap 
\overline{\Bru_{\tilde\varepsilon}} \ne \emptyset$.
By definition of the lifted Bruhat order in $\tilde B_{n+1}^{+}$ we have
$\varrho_k \le \tilde\varrho_k$.
Since this holds for all $k$, we have
$\varepsilon \preceq \tilde\varepsilon$, as desired.
\end{proof}

Notice that in Lemma \ref{lemma:poset} we neither claim equivalence, 
nor that either of the above conditions imply
$\BLS_{\varepsilon} \subseteq \overline{\BLS_{\tilde\varepsilon}}$.
The next example shows that Whitney's condition does not always hold.

\begin{example}
\label{example:bcbabdcb}
Consider $\sigma = a_2a_3a_2a_1a_2a_4a_3a_2$
and let $\varepsilon$, $\tilde{\varepsilon}$ 
be the two ancentries shown 
in Figure \ref{fig:posetdimtwotypeII}, in order:
\begin{gather*}
\varepsilon= (-2,-2,+1,-1,-1,+1,+2,+2), \quad
\tilde{\varepsilon} = (-2,-1,+2,-1,-2,+1,+1,+2), \\
\xi = (1,1,0,0,2,0,1,1), \quad \tilde{\xi} = (1,0,1,0,1,0,0,1).
\end{gather*}
We claim that 
\[
\varepsilon \preceq \tilde{\varepsilon}, \qquad
\overline{\BL_{\varepsilon}} \cap \BL_{\tilde{\varepsilon}} \ne \emptyset, \qquad
\BL_{\varepsilon} \not\subseteq \overline{\BL_{\tilde{\varepsilon}}}. \]
The fact that inclusion does not hold follows from the fact
that $\BL_{\varepsilon}$ and $\BL_{\tilde{\varepsilon}}$
are disjoint smooth submanifolds
of the same dimension.
Consider
\[ z_8 = \alpha_2(\theta_1) \alpha_3(\theta_2)
\alpha_2(\theta_3) \alpha_1(\theta_4)
\alpha_2(\theta_5) \alpha_4(\theta_6)
\alpha_3(\theta_7) \alpha_2(\theta_8) \in \Bru_{\acute\sigma}
\]
and corresponding $L$.
Take $\theta_1 = \pi/2$,
$\theta_3 = \theta_4 = \theta_6 = \theta_7 = \theta_8 = \pi/4$.
A computation shows that
for $\theta_2 = \pi/2$ and
$\theta_5 = \pi - \arctan(\sqrt{2})$,
we have $L \in \BL_{\varepsilon}$.
For $\theta_2 \in (\pi/3,\pi/2)$ and
$\theta_5 = \pi - \arctan(\sqrt{2})$,
we have $L \in \BL_{\tilde\varepsilon}$.
\end{example}


\begin{coro}
\label{coro:poset}
If $U$ is an upper set of ancestries then 
\[ \bigcup_{\varepsilon \in U} \BL_\varepsilon \subseteq \BL_\sigma \]
is an open subset.
\end{coro}

\begin{proof}
This follows from Lemma \ref{lemma:poset}.
\end{proof}




\section{Proof of Theorem \ref{theo:CWcomplex}}
\label{section:CW}

In this section we construct the CW complex $\BLC_\sigma$,
the continuous map $i_\sigma: \BLC_\sigma \to \BL_\sigma$
and prove Theorem \ref{theo:CWcomplex}.
The construction is recursive and for the induction step
we need a topological result, Lemma \ref{lemma:topolemma} below.

The idea is that the CW complex $\BLC_\sigma$
is a dual cell structure to the stratification.
Examples of this construction under nicer conditions
should be familiar from Poincaré duality
(see, for instance, \cite{Hatcher}, Section 3.3,
particularly the figure in page 232).
As we shall see, we have sufficient conditions to apply
a similar construction.
We did not find a reference in the required generality,
hence the need to state and prove Lemma \ref{lemma:topolemma}.
Another example of such a construction,
but involving Hilbert manifolds,
can be found in~\cite{Goulart-Saldanha2} (see Theorems 1 and 4).

Here $\Ss^{k-1}_r$ denotes the sphere of radius $r$,
$\BB^k_r$ denotes the open ball,
$\DD^k_r$ denotes the compact disk
and $\KK^k_{r_0,r_1}$ denotes the corona:
\begin{gather*}
\Ss^{k-1}_r = \{v \in \RR^k \;|\; |v| = r \}, \qquad
\BB^k_r = \{ v \in \RR^k \;|\; |v| < r \}, \\
\DD^k_r = \{ v \in \RR^k \;|\; |v| \le r \}, \qquad
\KK^k_{r_0,r_1} = \{ v \in \RR^k \;|\; r_0 \le |v| \le r_1 \}; 
\end{gather*}
also, $\DD^k = \DD^k_1$.
For a CW complex $X$, let $X^{[j]} \subseteq X$ denote
the skeleton of dimension $j$, that is,
the union of cells of dimension at most $j$.

\begin{lemma}
\label{lemma:topolemma}
Let $M_0 \subset M_1$ be smooth manifolds of dimension $\ell$.
Assume that $N_1 = M_1 \smallsetminus M_0 \subset M_1$
is a smooth submanifold of codimension $k$, $0 < k \le \ell$,
and that $N_1$ is diffeomorphic to $\RR^{\ell-k}$.
Assume that $X_0$ is a finite CW complex and that
$i_0: X_0 \to M_0$ is a homotopy equivalence.

There exists a map $\beta: \Ss^{k-1} \to X_0^{[k-1]}$
with the following properties.
Let $X_1$ be obtained from $X_0$ by attaching a cell $C_1$
of dimension $k$ with glueing map $\beta$.
There exists a map $i_1: X_1 \to M_1$ with $i_1|_{X_0} = i_0$ 
such that 
$i_1: X_1 \to M_1$ is a homotopy equivalence.
\end{lemma}

\begin{remark}
\label{remark:topolemma}
Notice that since $M_0 \subset M_1$ is a submanifold
of codimension $0$ it follows that $M_0$ is an open subset of $M_1$.
The subset $N_1 \subset M_1$ is therefore closed.
If $k < \ell$ it follows that $M_1$ is not compact.

The maps $i_0$ and $i_1$ can be taken to be inclusions in many examples
but are not required to be so.
The first paragraph of the proof provides us with a construction of
the map $\beta$ and therefore of the CW complex $X_1$.
The map $i_1: X_1 \to M_1$ is constructed slightly later. 
\end{remark}

\begin{proof}[Proof of Lemma \ref{lemma:topolemma}]
Consider a small transversal section to $N_1$,
the map $\alpha_1: \DD^k_{\frac12} \to M_1$ 
with $\alpha_1(0) = z_1 \in N_1$.
Consider the restriction $\beta_1 = \alpha_1|_{\Ss^{k-1}_{\frac12}}$
such that, ignoring the radius, $\beta_1: \Ss^{k-1} \to M_0$. 
Take $\beta: \Ss^{k-1} \to X_0^{[k-1]}$ such that
$\beta_1$ and $i_0 \circ \beta$ are homotopic in $M_0$.
This is the desired glueing map $\beta: \Ss^{k-1} \to X_0^{[k-1]}$,
completing the construction of the CW complex $X_1$. 
We must now construct $i_1: X_1 \to M_1$
and prove that $i_1$ is a homotopy equivalence.

We first go again through the construction of the map $\beta$,
but more carefully, collecting more notation and information.
By our assumptions, the map $i_0: X_0 \to M_0$ is a homotopy equivalence. 
Thus, there exist a continuous map $p_0: M_0 \to X_0$
and two homotopies $H_0: [0,1]\times M_0 \to M_0$ and
$\tilde H_0: [0,1]\times X_0 \to X_0$
with 
\[ 
H_0(0,z) = z, \quad H_0(1,z) = i_0(p_0(z)), \quad
\tilde H_0(0,x) = x, \quad \tilde H_0(1,x) = p_0(i_0(x)) \]
for all $z \in M_0$ and $x \in X_0$.
Consider a tubular neighborhood of $N_1$ disjoint
from the compact set $i_0[X_0]$.
Since $N_1$ is diffeomorphic to $\RR^{\ell - k}$,
the tubular neighborhood may be assumed to be a smooth injective map
$\Phi: \DD_{\frac12}^k \times \RR^{\ell - k} \to M_1$ with
$\Phi[\{0\} \times \RR^{\ell-k}] = N_1$.
Let $\alpha_1: \DD^k_{\frac12} \to M_1$,
$\alpha_1(x) = \Phi(x,0)$, $z_1 = \alpha_1(0)$
and consider the restriction $\beta_1$, as above.


Define $\beta_2 = p_0 \circ \beta_1: \Ss^{k-1} \to X_0$:
notice that $\beta_1$ and $i_0 \circ \beta_2$ are homotopic in $M_0$,
with homotopy $H_0(\cdot,\beta_1(\cdot))$.
Also, there exists $\beta: \Ss^{k-1} \to X_0^{[k-1]}$
such that $\beta_2$ and $\beta$ are homotopic in $X_0$;
let $H_X: [0,1] \times \Ss^{k-1} \to X_0$ be such a homotopy.
Thus, $\beta_1$ and $i_0 \circ \beta$ are homotopic in $M_0$.

We proceed to construct $i_1: X_1 \to M_1$.
Let $\alpha_2: \KK^k_{\frac12,1} \to M_0$ satisfy
$\alpha_2|_{\Ss^{k-1}_{\frac12}} = \beta_1$,
$\alpha_2|_{\Ss^{k-1}_{\frac34}} = i_0 \circ \beta_2$ and
$\alpha_2|_{\Ss^{k-1}_{1}} = i_0 \circ \beta$.
More precisely, for $r \in [\frac12,1]$ and $u \in \Ss^{k-1}$, set
\[ \alpha_2(ru) = \begin{cases}
H_0(4r-2,\beta_1(u)), & r \in [\frac12,\frac34], \\
i_0(H_X(4r-3,u)), & r \in [\frac34,1].
\end{cases}
\]
Define $\alpha: \DD^k \to M_1$ by
$\alpha|_{\DD^k_{\frac12}} = \alpha_1$ and
$\alpha|_{\KK^k_{\frac12,1}} = \alpha_2$.
Define $i_1$ by $i_1|_{X_0} = i_0$ 
and by $i_1(x) = \alpha(x)$ for $x \in C_1 = \DD^k$.

We need to prove that $i_1$ is a homotopy equivalence.
We could at this point construct $p_1: M_1 \to X_1$
and homotopies $H_1$ and $\tilde H_1$.
Since that construction is rather cumbersome,
we prefer to proceed in a slightly different way:
we first prove that for any $j \ge 0$ the map
$\pi_j(i_1): \pi_j(X_1) \to \pi_j(M_1)$ is a bijection.

We first prove the surjectivity of $\pi_j(i_1)$.
Let $\gamma_{M}: \Ss^j \to M_1$:
we want to prove that there exists  $\gamma_{X}: \Ss^j \to X_1$
such that $\gamma_M$ and $i_1 \circ \gamma_X$ are homotopic.
We may assume that $\gamma_M$ is smooth and transversal to $N_1$:
let $N_{S} = \gamma_M^{-1}[N_1] \subset \Ss^j$,
a smooth submanifold of codimension $k$.
By transversality, there exist $\epsilon \in (0,\frac12)$
and a tubular neighborhood
$\Psi: \DD^k_\epsilon \times N_S \to \Ss^j$
with $\Psi(0,s) = s$ (for all $s \in N_S$).
We may furthermore assume that there exists
a smooth function $f_0: \DD^k_\epsilon \times N_S \to \RR^{\ell - k}$
such that $\gamma_M(\Psi(v,s)) = \Phi(v,f_0(v,s))$
(here $\Phi$ is the tubular neighborhood of $N_1$ described above).

Multiplication of $f_0$ by a bump function takes us from $\gamma_M$
to a homotopic function $\gamma_{M,1}$ such that
$\gamma_M$ and $\gamma_{M,1}$ coincide in
$\Ss^j \smallsetminus \Psi[\BB^k_{\epsilon} \times N_S]$  and
$\gamma_{M,1}(\Psi(v,s)) = \alpha_1(v)$
for all $s \in N_S$ and $v \in \DD^k_{\epsilon/2}$.
Another homotopy takes us to $\gamma_{M,2}$ such that
$\gamma_{M,1}$ and $\gamma_{M_2}$ coincide in 
$\Ss^j \smallsetminus \Psi[\BB^k_{\epsilon/2} \times N_S]$  and,
if $s \in N_S$ and $v \in \DD^k_{\epsilon/2}$ then
\[ \gamma_{M,2}(\Psi(v,s)) = \begin{cases}
\alpha_1(2v/\epsilon), & |v| \le \epsilon/4; \\
\alpha_1(v/(2|v|)), & \epsilon/4 \le |v| \le 3\epsilon/8; \\
\alpha_1(( 2 - (3\epsilon/2) - (4(1-\epsilon)/\epsilon)|v|)v), & 
3\epsilon/8 \le |v| \le \epsilon/2. \end{cases} \]
We now compose with $H_0$ to obtain $\gamma_{M,3}$ as follows:
\[ \gamma_{M,3}(s) = \begin{cases}
(i_0 \circ p_0 \circ \gamma_{M,2})(s), &
s \notin \Psi[\BB^k_{3\epsilon/8} \times N_S], \\
H_0(8(|v|-(\epsilon/4))/\epsilon,\gamma_{M,2}(s)), &
s = \Psi(v,s_0), \, |v| \in [\epsilon/4,3\epsilon/8], \\
\gamma_{M,2}(s), &
s \in \Psi[\BB^k_{\epsilon/4} \times N_S].
\end{cases} \]
Notice that for $s \in N_S$ and $|v| \le 3\epsilon/8$ we have 
$\gamma_{M,3}(\Psi(v,s)) = \alpha(2v/\epsilon)$.
For $s \in \Ss^j \smallsetminus \Psi[\BB^k_{3\epsilon/8} \times N_S]$,
we have $\gamma_{M,3}(s) = (i_0 \circ p_0 \circ \gamma_{M,2})(s)$.
A small adjustment using $H_X$ takes us to
$\gamma_{M,4} = i_1 \circ \gamma_{X,4}$,
completing the proof of surjectivity.

The proof of injectivity of $\pi_j(i_1)$ is similar,
and will be presented in less detail.
Consider $\gamma_X: \Ss^{j} \to X_1$
and assume that $\pi_j(i_1)(\gamma_X) = 0 \in \pi_j(M_1)$.
We may assume that $\gamma_X$ is smooth in the interior of the new cell $C_1$
and that the center $0_{C_1}  \in C_1$ is a regular value.
By hypothesis, there exists $\Gamma_M: \DD^{j+1} \to M_1$
such that $\Gamma_M|_{\Ss^j} = i_1 \circ \gamma_X$.
Again, we may assume that $\Gamma_M$ is smooth
in a neighborhood of $\Gamma_M^{-1}[0_{C_1}]$
and that $0_{C_1}$ is a regular value.
Let $N_D = \Gamma_M^{-1}[0_{C_1}] \subset \DD^{j+1}$:
this is a smooth submanifold of codimension $k$
with boundary $\partial N_D \subset \Ss^{j}$,
also a smooth submanifold of codimension $k$.
As above, pulling back $\Phi$ gives us a tubular neighborhood $\Psi$.
Again as above, we construct $\Gamma_{M,\ast}$ of the form
$\Gamma_{M,\ast} = i_1 \circ \Gamma_{X,\ast}$
where $\Gamma_{X,\ast}: \DD^{j+1} \to X_1$ satisfies
$\Gamma_{X,\ast}|_{\Ss^j} = \gamma_X$.
We thus have $[\gamma_X] = 0 \in \pi_j(X_1)$,
completing the proof of injectivity.

At this point we know that $i_1: X_1 \to M_1$
is such that $\pi_j(i_1)$ is bijective for all $j$.
In other words, $i_1$ is a weak homotopy equivalence.
The set $X_1$ is a CW complex and $M_1$ is a manifold
and therefore homeomorphic to a CW complex.
By Whitehead's Theorem,
the map $i_1$ is a homotopy equivalence, as desired.
\end{proof}

\begin{proof}[Proof of Theorem \ref{theo:CWcomplex}]
Ancestries of dimension $0$ are maximal elements under 
the partial order $\succeq$.
Let $\BL_{\sigma;0} \subseteq \BL_{\sigma}$
be the union of the open, disjoint, contractible sets
$\BL_\varepsilon$ for $\varepsilon$ an ancestry of dimension $0$.
The set $\BL_{\sigma;0}$ is homotopically equivalent
to a finite set with one vertex per ancestry, 
which is of course a CW complex of dimension $0$.
This is the basis of a recursive construction.

We can list the set of ancestries of positive dimension as
$(\varepsilon_i)_{1 \le i \le N_\varepsilon}$
in such a way that $\varepsilon_i \succeq \varepsilon_j$ implies $i \le j$.
Define recursively the subsets 
$\BL_{\sigma;i} = \BL_{\sigma;i-1} \cup \BLS_{\varepsilon_i}
\subseteq \BL_{\sigma}$.
The family of sets $\BL_{\sigma;i}$ defines a filtration:
\begin{equation*}
\label{equation:filtration}
\BL_{\sigma;0} \subset \BL_{\sigma;1} \subset \cdots
\subset \BL_{\sigma;N_\xi-1} \subset \BL_{\sigma;N_\xi} = \BL_{\sigma}.
\end{equation*}
The partial order $\succeq$ and Lemma \ref{lemma:poset} guarantee
that $\BL_{\sigma;i-1} \subset \BL_{\sigma;i}$ is an open subset.
Lemma \ref{lemma:submanifold} tells us that
$\BLS_{\varepsilon_i} = \BL_{\sigma;i} \smallsetminus \BL_{\sigma;i-1}
\subset \BL_{\sigma;i}$ is a smooth submanifold
of codimension $d = \dim(\varepsilon_i)$
and Lemma \ref{lemma:contractible} tells us that
$\BL_{\varepsilon_i}$ is diffeomorphic to $\RR^{\ell-d}$.
Notice that $\BLS_{\varepsilon_i} \subset \BL_{\sigma;i}$
is a closed subset (see Remark \ref{remark:topolemma}).
We may therefore apply Lemma \ref{lemma:topolemma}
to the pair $M_0 = \BL_{\sigma;i-1} \subset \BL_{\sigma;i} = M_1$,
completing the recursive construction and the proof.
\end{proof}


The proofs of Lemma \ref{lemma:topolemma} 
and of Theorem \ref{theo:CWcomplex}
give us instructions for the actual construction of the CW complex
$\BLC_\sigma$ and of the map $i_\sigma$.
This construction of the CW complex and of the glueing maps
are not as direct as might perhaps be desired.
This is the subject of Section \ref{section:glue}.


\section{Euler characteristic}
\label{section:euler}

We begin with an easy formula
for the Euler characteristic of $\BL_z$.

\begin{lemma}
\label{lemma:euler}
For $\sigma \in S_{n+1}$ and $z \in \acute\sigma\Quat_{n+1}$,
we have
\[ \chi(\BL_z) = \sum_{\varepsilon_0}
(-1)^{\dim{\varepsilon_0}} N_{\varepsilon_0}(z). \]
The summation is taken over all preancestries $\varepsilon_0$.
\end{lemma}

\begin{proof}
This follows directly from Theorem \ref{theo:CWcomplex}
and the definition of $N_{\varepsilon_0}(z)$.
\end{proof}

\begin{lemma}
\label{lemma:eveneuler}
Let $z_0 \in \acute\eta\Quat_{n+1}$ be such that $\Re(z_0) > 0$.
We have that $\chi(\BL_{z_0})$ is odd and $\chi(\BL_{-z_0})$ is even.
\end{lemma}

\begin{proof}
We have $\ell = \inv(\eta) = n(n+1)/2$ and
$c = \nc(\eta) = 1 + \lfloor n/2 \rfloor$.
From Lemma \ref{lemma:realpartsigma},
$\Re(z_0) = 2^{-(n+1-c)/2} = 2^{-(\lceil n/2 \rceil)/2}$. 

From Lemma \ref{lemma:eta},
the highest possible dimension $d_{\max} = \lfloor n^2/4 \rfloor$
is achieved by a unique preancestry $\varepsilon_{\max}$ 
(see Figure \ref{fig:preancestryeta}).
This preancestry shall be addressed separately.
We have $\ell - 2d_{\max} = \lceil n/2 \rceil$.
Thus, from
the first equation in Theorem \ref{theo:two} we have
$N_{\varepsilon_{\max}}(z_0) - N_{\varepsilon_{\max}}(-z_0) = 1$.
We have $H_{\varepsilon_{\max}} = H_{\eta}$ 
and $|H_{\eta}| = 2^{n+2-c} = 2^{1+\lceil n/2 \rceil}$.
The second equation in Theorem \ref{theo:two}
thus gives us 
$N_{\varepsilon_{\max}}(z_0) + N_{\varepsilon_{\max}}(-z_0) = 1$.
We thus have
$N_{\varepsilon_{\max}}(z_0) = 1$ and
$N_{\varepsilon_{\max}}(-z_0) = 0$.

Consider a preancestry $\varepsilon_0$
with $d = \dim(\varepsilon_0) = d_{\max} - 1$.
The first equation in Theorem \ref{theo:two}
tells us that
$N_{\varepsilon_0}(z_0) - N_{\varepsilon_0}(-z_0) = 2$.
The second equation tells us that
$N_{\varepsilon_0}(z_0) + N_{\varepsilon_0}(-z_0) =
4|H_{\sigma}|/|H_{\varepsilon_0}|$, a power of $2$.
By definition of $H_{\varepsilon_0}$, however,
we have $|H_{\varepsilon_0}| > |H_{\sigma}|$
and therefore 
$N_{\varepsilon_0}(z_0) + N_{\varepsilon_0}(-z_0) \le 2$.
We thus have
$N_{\varepsilon_{0}}(z_0) = 2$ and
$N_{\varepsilon_{0}}(-z_0) = 0$.

Consider finally a preancestry $\varepsilon_0$
with $d = \dim(\varepsilon_0) < d_{\max} - 1$.
The first equation in Theorem \ref{theo:two}
tells us that
$N_{\varepsilon_0}(z_0) - N_{\varepsilon_0}(-z_0)$
is a power of $2$ and a multiple of $4$.
The second equation tells us that
$N_{\varepsilon_0}(z_0) + N_{\varepsilon_0}(-z_0)$
is also a power of $2$ and therefore also a multiple of $4$.
Thus, $N_{\varepsilon_{0}}(z_0)$ and
$N_{\varepsilon_{0}}(-z_0)$ are both even.

Summing up, $N_{\varepsilon_{0}}(-z_0)$ is even
for every preancestry.
Similarly, $N_{\varepsilon_{0}}(z_0)$ is even
for every preancestry except $\varepsilon_{\max}$.
The result follows from Lemma \ref{lemma:euler}.
\end{proof}

\begin{coro}
\label{coro:eveneuler}
Consider $n \ge 5$ and
$z_0 \in \acute\eta\Quat_{n+1}$ with $\Re(z_0) > 0$.
Then $\BL_{-z_0,\thick}$ is non empty, connected
and its Euler characteristic $\chi(\BL_{-z_0,\thick})$ is even.
\end{coro}

\begin{proof}
From Proposition \ref{prop:effective},
$\BL_{-z_0,\thick}$ is non empty and connected.
From Lemma \ref{lemma:eveneuler},
$\chi(\BL_{-z_0})$ is even.
From Lemma \ref{lemma:Nthin},
$N_{\thin}(-z_0)$ is even.
We have $\chi(\BL_{-z_0,\thick}) = \chi(\BL_{-z_0}) - N_{\thin}(-z_0)$:
the result follows.
\end{proof}

\begin{remark}
\label{remark:eveneuler}
Corollary \ref{coro:eveneuler} implies
the third item in Theorem \ref{theo:one}.
For $n = 5$,
$\BLC_{-z_0}$ is connected and has:
$480$ vertices,
$1120$ cells of dimension $1$,
$864$ cells of dimension $2$,
$228$ cells of dimension $3$,
$6$ cells of dimension $4$
and no cells of higher codimension.
It follows that 
$\chi(\BL_{-z_0}) = 480-1120+864-228+6 = 2$:
in particular, $\BL_{-z_0}$ is not contractible.
It would be very interesting to compute its homotopy type.
\end{remark}


\section{The glueing maps}
\label{section:glue}

The following results gives us a little more information
about the CW complex $\BLC_\sigma$,
particularly about the glueing maps.

Recall that $U_\varepsilon$ is the upper set of ancestries
$\tilde\varepsilon$ with $\varepsilon \preceq \tilde\varepsilon$.
Let $U_\varepsilon^\ast = U_\varepsilon \smallsetminus \{\varepsilon\}$,
also an upper set.
In general, for an upper set $U$ of ancestries, define
\[ \BLS_U = \bigcup_{\varepsilon \in U} \BLS_\varepsilon \subseteq \BL_\sigma,
\qquad
\BLC_U = \bigcup_{\varepsilon \in U} \BLC_\varepsilon \subseteq \BLC_\sigma. \]
We know from Corollary \ref{coro:poset} that $\BLS_U \subseteq \BL_\sigma$
is an open subset.

\begin{lemma}
\label{lemma:upperset}
Let $U$ be an upper set of ancestries.
The subset $\BLC_U \subseteq \BLC_\sigma$ is closed and a CW complex.
The restriction $i_\sigma|_{\BLC_U}: \BLC_U \to \BLS_U$
is a homotopy equivalence.
\end{lemma}

\begin{proof}
We use induction on $|U|$, the cases $|U| \le 1$ being trivial. 
Let $\varepsilon \in U$ be a minimal element,
so that $U^\ast = U \smallsetminus \{\varepsilon\}$
is also an upper set.
By induction hypothesis, the result holds for $U^\ast$.
It follows from Lemma \ref{lemma:topolemma}
and the proof of Theorem \ref{theo:CWcomplex}
that the image of the glueing map for the cell $\BLC_\varepsilon$
is contained in $\BLC_{U^\ast}$,
proving the first claim for $U$.
The second claim follows from Lemma \ref{lemma:topolemma}.
\end{proof}

\begin{coro}
\label{coro:boundary}
The image of the glueing map for $\BLC_\varepsilon$
in contained in $\BLC_{U_\varepsilon^\ast}$.
\end{coro}

\begin{proof}
This follows directly from Lemma \ref{lemma:upperset}.
\end{proof}

Consider an ancestry $\varepsilon$ of positive dimension
$d = \dim(\varepsilon) > 0$.
We define two non empty subsets
$U^{\pm}_\varepsilon \subset U^{\ast}_\varepsilon$.
Let $k_\bullet$ be the largest index $k$
for which $\varepsilon(k) = -2$.
We have $\varrho_{k_\bullet} = \varrho_{k_\bullet - 1} \acute a_{i_k}$:
define $\varrho_{k_\bullet}^{-} = \varrho_{k_\bullet - 1}$
and $\varrho_{k_\bullet}^{+} = \varrho_{k_\bullet - 1} \hat a_{i_k}$.
For $\tilde\varepsilon \in U^{\ast}_\varepsilon$,
let $(\tilde\varrho_k)_{0 \le k \le \ell}$ be defined as usual.
We have
\begin{equation}
\label{equation:Upm}
\tilde\varepsilon \in U^{\pm}_{\varepsilon} \quad\iff\quad
((\tilde\varrho_{k_{\bullet}} = \varrho_{k_\bullet}^{\pm}) \land
(\forall k, (0 \le k < k_\bullet) \to (\tilde\varrho_k = \varrho_k))).
\end{equation}
The sets $U^{\pm}_{\varepsilon}$ are clearly disjoint.

\begin{example}
\label{example:Upmdimone}
We already discussed $U_\varepsilon$ for $\dim(\varepsilon) = 1$
in Example \ref{example:posetdimone}.
In this case, the sets $U^{\pm}_{\varepsilon}$
have one element each.

If $\varepsilon$ is an ancestry of dimension two, type I
(see Remark \ref{remark:lowdimpreancestry} and
Example \ref{example:posetdimtwo})
then the sets  $U^{\pm}_{\varepsilon}$ also
have one element each.
In Figure \ref{fig:posetdimtwotypeI},
these are the edges at the top and bottom of the figure.

If $\varepsilon$ is an ancestry of dimension two, type II
then the sets  $U^{\pm}_{\varepsilon}$ can be more complicated.
In particular, if the ancestries in Figure \ref{fig:posetdimtwotypeII}
are $\varepsilon_0$ and $\varepsilon_1$,
we have $\varepsilon_1 \in U^{-}_{\varepsilon_0}$
(see also Example \ref{example:bcbabdcb}).
\end{example}

Define
\begin{equation}
\label{equation:BLSpm}
\BLS_{\varepsilon}^{\pm} =
\bigcup_{\tilde\varepsilon \in U_\varepsilon^{\pm}} \BLS_{\tilde\varepsilon}. 
\end{equation}
The following result describes these sets near $\BLS_\varepsilon$.

\begin{lemma}
\label{lemma:BLSpm}
Let $\varepsilon$ be an ancestry of dimension $d = \dim(\varepsilon) > 0$.
If $W$ is a sufficiently thin open tubular neighborhood of $\BLS_{\varepsilon}$
then
$(\BLS_{\varepsilon} \cup \BLS^{\pm}_{\varepsilon}) \cap W \subset W$
are smooth submanifolds with boundary.
Both manifolds have codimension $d-1$
and boundary equal to $\BLS_{\varepsilon}$.

Let $W^\ast = W \smallsetminus \BLS_{\varepsilon}$.
There exists a diffeomorphism 
$\Phi: \Ss^{d-1} \times (0,r) \times \RR^{\ell-d} \to W^\ast$
such that
\[
\Phi^{-1}[\BLS^{+}_{\varepsilon}] = \{\bN\} \times (0,r) \times \RR^{\ell-d},
\quad
\Phi^{-1}[\BLS^{-}_{\varepsilon}] = \{\bS\} \times (0,r) \times \RR^{\ell-d},
\]
where $\bN, \bS \in \Ss^{d-1}$ are the north and south poles.
\end{lemma}

\begin{proof}
Let $(z_k)_{0 \le k \le \ell}$ be sequences
as in Equation \eqref{equation:zk}.
For $W \supset \BLS_{\varepsilon}$ as above and $L \in W$,
we identify $L \in W$ with such a sequence
with $z_\ell \in \grave\eta\Bru_{\varrho_\ell}$.
For $L \in W$, we have $L \in \BLS^{\pm}_{\varepsilon}$ if and only if:
\[ (z_{k_{\bullet}} \in \grave\eta\Bru_{\varrho_{k_\bullet}^{\pm}}) \land
(\forall k, \; (0 \le k < k_\bullet) \to
(z_k \in \grave\eta\Bru_{\varrho_k})). \]
The condition above includes $d-1$ transversal equations,
corresponding to $k < k_{\bullet}$, $\varepsilon(k) = -2$:
the remaining are open conditions.
This completes the proof of the first paragraph (the first three sentences in the statement).
The other claims are then easy. 
\end{proof}

Let $M$ be a smooth manifold and
$N \subset M$ be a transversally oriented submanifold of codimension $k$ 
which is also a closed set.
Intersection with $N$ defines an element of $H^k(M;\ZZ)$.
In $W^\ast$ (as defined in Lemma \ref{lemma:BLSpm}),
intersection with either $\BLS^{\pm}_{\varepsilon}$
defines a generator of $H^{d-1}(W^\ast;\ZZ) \approx \ZZ$.

An ancestry $\varepsilon$ with $d = \dim(\varepsilon) > 0$
is \textit{tame} if the following conditions hold:
\begin{enumerate}
\item{The manifold $\BLS_{U^{\ast}_{\varepsilon}}$
is homotopically equivalent to $\Ss^{d-1}$.}
\item{Intersection with $\BLS^{\pm}_{\varepsilon}$
defines generators
of $H^{d-1}(\BLS_{U^{\ast}_{\varepsilon}};\ZZ) \approx \ZZ$.}
\end{enumerate}
Otherwise, $\varepsilon$ is \textit{wild}.

Recall that $\BLS_{U^{\ast}_{\varepsilon}}$
is homotopically equivalent to the CW complex $\BLC_{U^{\ast}_{\varepsilon}}$.
Let us translate the definition above
in terms of $\BLC_{U^{\ast}_{\varepsilon}}$.
The first condition says of course that $\BLC_{U^{\ast}_{\varepsilon}}$
is homotopically equivalent to $\Ss^{d-1}$.
The second condition says that we can build cocycles
$\omega_{\BLC}^{\pm} \in Z^{d-1}(\BLC_{U^{\ast}_{\varepsilon}};\ZZ)$
by taking the elements of $U^{\pm}_{\varepsilon}$ of dimension $d-1$
and interpreting them as cells of $\BLC_{U^{\ast}_{\varepsilon}}$:
the elements $\omega_{\BLC}^{\pm}$ are then generators of $H^{d-1} \approx \ZZ$.

\begin{example}
\label{example:tame}
It follows from Example \ref{example:posetdimone}
that $\dim(\varepsilon) = 1$ implies that $\varepsilon$ is tame. 
Similarly, Example \ref{example:posetdimtwo}
that if $\varepsilon$ is of dimension two type I
then $\varepsilon$ is tame.
\end{example}

The following lemma tells us how to obtain the glueing map
in tame cases.

\begin{lemma}
\label{lemma:tame}
If $\varepsilon$ is tame then the glueing map
$\beta: \Ss^{d-1} \to \BLC_{U^{\ast}_{\varepsilon}}$
is a homotopy equivalence.
\end{lemma}

\begin{proof}
Let $W^\ast$ be as in Lemma \ref{lemma:BLSpm}.
Let $\omega_{W^\ast}^{\pm} \in H^{d-1}(W^\ast;\ZZ)$
be defined by intersection with $\BLS^{\pm}_{\varepsilon}$;
by definition of tameness, either one is a generator.
Let $\beta_1: \Ss^{d-1} \to W^\ast$ be as in the first lines of 
the proof of Lemma \ref{lemma:topolemma}
(with $M_0 = \BLS_{U^{\ast}_{\varepsilon}}$,
$M_1 = \BLS_{U_{\varepsilon}}$ and $k = d$).
We have a pairing
$H^{d-1}(W^\ast;\ZZ) \times \pi_{d-1}(W^\ast) \to \ZZ$.
By Lemma \ref{lemma:BLSpm},
$|\omega_{W^\ast}^{\pm} \beta_1| = 1$.

Let $i: W^\ast \to \BLS_{U^{\ast}_{\varepsilon}}$
be the inclusion.
Let $\omega_{\BLS}^{\pm} \in H^{d-1}(\BLS_{U^{\ast}_{\varepsilon}};\ZZ)$
be defined by intersection with $\BLS^{\pm}_{\varepsilon}$,
as in the definition of tameness.
Let 
\[ i^\ast = H^{d-1}(i): H^{d-1}(\BLS_{U^{\ast}_{\varepsilon}};\ZZ)
\to H^{d-1}(W^\ast;\ZZ); \]
we have $i^\ast(\omega_{\BLS}^{\pm}) = \omega_{W^\ast}^{\pm}$.
We thus have
$\omega_{\BLS}^{\pm} (i \circ \beta_1) = \omega_{W^\ast}^{\pm} \beta_1$
and $i \circ \beta_1$ is therefore a generator of 
$\pi_{d-1}(\BLS_{U^{\ast}_{\varepsilon}})$.
From the proof of Lemma \ref{lemma:topolemma},
so is the glueing map $\beta$.
The result follows.
\end{proof}

\begin{example}
\label{example:tameglue}
It follows from Lemma \ref{lemma:tame}
together with Examples \ref{example:posetdimone} and \ref{example:tame}
that cells $\BLC_\varepsilon$ of dimension $1$ in $\BLC_\sigma$
are edges joining the two vertices corresponding
to the elements of dimension $0$ in $U_\varepsilon$.
This is illustrated in Figure \ref{fig:posetdimone}
(among many others).

Similarly,
Examples \ref{example:posetdimtwo} and \ref{example:tame}
imply that if $\varepsilon$ has dimension two type I
then $\BLC_\varepsilon$ fills in a square hole,
as in Figure \ref{fig:posetdimtwotypeI}.
\end{example}

In the next sections we shall see several other examples of tame ancestries.
So far we have not seen \textit{any} example of a wild ancestry, 
which raises the following doubt:

\begin{question}
\label{question:wild}
Are there \emph{any} wild ancestries?
\end{question}

We do not know the answer:
it seems entirely plausible that in higher dimension wild ancestries exist.




\section{Examples and proof of
Theorems \ref{theo:one} and \ref{theo:collapse}}
\label{section:moreexamples}

In this section we apply the previous results,
particularly Theorem \ref{theo:two},
Lemmas \ref{lemma:tame} and \ref{lemma:euler}
to compute the homotopy type of $\BL_\sigma$
and of its subsets $\BL_z$ for several examples.
Let $N(z) = N_{\varepsilon_0}(z)$
(where $\varepsilon_0$ is the empty preancestry)
be the number of vertices in the CW complex $\BLC_z$.

\begin{example}
\label{example:4231}
Consider $n = 3$ and $\sigma = 4231 = a_1a_2a_3a_2a_1$.
A simple computation gives
$\acute\sigma = (\hat a_2 + \hat a_1\hat a_3)/\sqrt{2}$ and
\[ \acute\sigma\Quat_4 =
\left\{
\frac{\pm 1 \pm \hat a_1\hat a_2\hat a_3}{\sqrt 2},
\frac{\pm \hat a_1 \pm \hat a_2\hat a_3}{\sqrt 2},
\frac{\pm \hat a_2 \pm \hat a_1\hat a_3}{\sqrt 2},
\frac{\pm \hat a_3 \pm \hat a_1\hat a_2}{\sqrt 2} \right\}. \]
The action of $\cE_3$ on $\acute\sigma\Quat_4$ has $5$ orbits.
The orbit
$\cO_{\acute\sigma} = \{ (\pm \hat a_2 \pm \hat a_1\hat a_3)/\sqrt{2} \}$
has size $4$, and for each $z$ we have $N(z) = N_{\thin}(z) = 2$
so that $\BL_z$ has two thin (and therefore contractible)
connected components.
The orbits
$\cO_{\hat a_1\acute\sigma}
= \{ (\pm \hat a_3 \pm \hat a_1\hat a_2)/\sqrt{2} \}$
and
$\cO_{\hat a_3\acute\sigma}
= \{ (\pm \hat a_1 \pm \hat a_2\hat a_3)/\sqrt{2} \}$
both have size $4$.
For each $z$ in one of these two orbits, we have
$N(z) = 2$ and $N_{\thin}(z) = 0$.
Figure \ref{fig:4231x} shows the stratification of $\BL_z$
for one representative of each orbit:
it follows that such sets $\BL_z$ are contractible.
The orbit
$\cO_{\hat a_2\acute\sigma}
= \{ ( - 1 \pm  \hat a_1 \hat a_2\hat a_3)/\sqrt{2} \}$
has size $2$.
For $z$ in this orbit, we have $N(z) = 0$, 
so that the corresponding sets $\BL_z$ are empty.

\begin{figure}[ht]
\begin{center}
\includegraphics[scale=0.35]{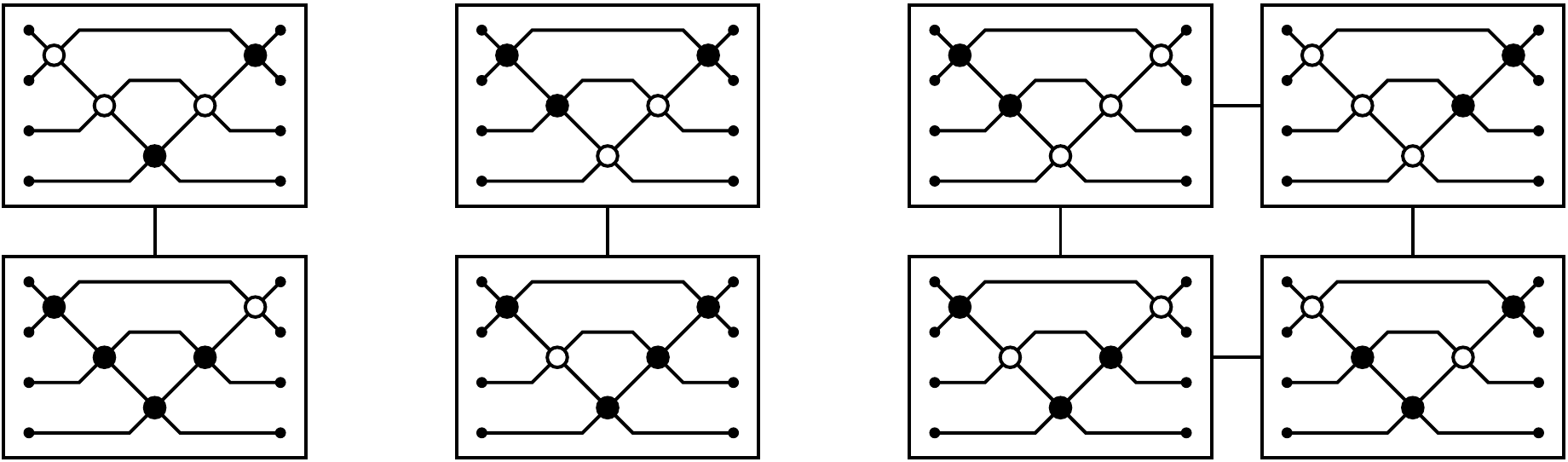}
\end{center}
\caption{The stratifications of
$\BL_{\hat a_1\acute\sigma}$,
$\BL_{\hat a_3\acute\sigma}$,
and
$\BL_{-\hat a_2\acute\sigma}$ for
$\sigma = a_1a_2a_3a_2a_1$.  }
\label{fig:4231x}
\end{figure}

Finally, the orbit
$\cO_{-\hat a_2\acute\sigma}
= \{ ( 1 \pm  \hat a_1 \hat a_2\hat a_3)/\sqrt{2} \}$
also has size $2$:
consider 
$z = -\hat a_2\acute\sigma =
\acute a_1\acute a_2\acute a_3\grave a_2\grave a_1 =
( 1 + \hat a_1 \hat a_2\hat a_3)/\sqrt{2}$.
Figure \ref{fig:4231x} also shows the stratification of $\BL_z$.
A straightforward computation verifies that
the set $\BL_z$ consists of matrices of the form:
\begin{equation}
\label{equation:4231}
L = \begin{pmatrix}
1 & & & \\ l_{21} & 1 & & \\
l_{31} & l_{32} & 1 & \\ l_{41} & l_{42} & l_{43} & 1 
\end{pmatrix}, \qquad
\begin{aligned}
l_{41} &> \max\{0, l_{21}l_{42}, l_{31}l_{43} \}, \\
l_{32} &= \frac{l_{31}l_{42}}{l_{41}}.
\end{aligned}
\end{equation}
The above description makes it clear that $\BL_z$ is contractible,
but we want to explore its decomposition into strata. 
The set $\BL_z$ contains $4$ open strata,
$4$ strata of codimension $1$
and one stratum of codimension $2$,
with ancestry $\varepsilon = (-2,-2,+1,+2,+2)$.
A computation shows that $\varepsilon$ is tame.
The set $\BL_{\varepsilon}$ is the set of matrices of the above form  
with $l_{32} = l_{42} = l_{43} = 0$.
The open strata are characterized by the signs of 
$l_{42}$ and $l_{43}$.
In our CW complex, the ancestry $\varepsilon = (-2,-2,+1,+2,+2)$
corresponds to a square, glued along the four edges
(corresponding to ancestries of codimension $1$)
in the obvious way.
Thus, $\BL_z$ is also contractible.
Summing up, $\BL_{\sigma}$ has $18$ connected components.
Each connected component of $\BLC_{\sigma}$ collapses to a point.
Each connected component of $\BL_{\sigma}$ is therefore contractible.
\end{example}

\begin{example}
\label{example:4321}
We take $n = 3$ and $\sigma = \eta = a_1a_2a_1a_3a_2a_1$,
the top permutation (with $\ell = 6$);
see Figure \ref{fig:4321-a12}.
We shall verify that the set $\BL_\eta$ has $20$ connected components,
all contractible.
The total number of connected components of $\BL_\eta$
was first calculated by Vl.~Kostov, B.~Shapiro and M.~Shapiro
using  ad hoc methods back in 1987 (unpublished);
see also \cite{GeShVa, SSV1}.

\begin{figure}[ht]
\begin{center}
\includegraphics[scale=0.3]{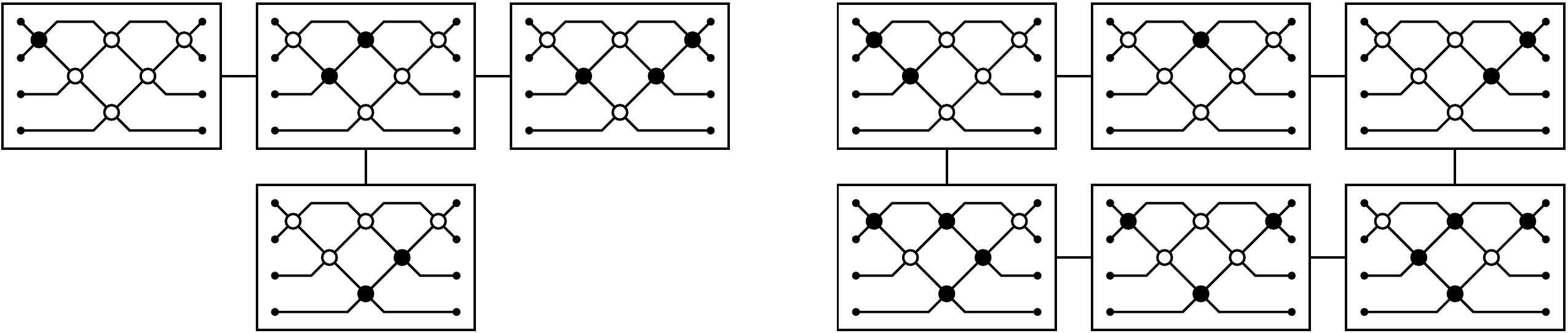}
\end{center}
\caption{The CW complexes
$\BLC_{-\hat a_1\acute\eta}$ and $\BLC_{-\hat a_2\acute\eta}$.  }
\label{fig:4321-a12}
\end{figure}

We have
\begin{gather*}
\acute\eta =
\frac{-1+\hat a_2 + \hat a_1\hat a_3 - \hat a_1\hat a_2\hat a_3}{2},
\qquad
\grave\eta =
\frac{-1-\hat a_2 + \hat a_1\hat a_3 + \hat a_1\hat a_2\hat a_3}{2}, \\
\acute\eta\Quat_4 = \left\{
\frac{\pm 1 \pm \hat a_2 \pm \hat a_1\hat a_3 \pm \hat a_1\hat a_2\hat a_3}{2},
\frac{\pm \hat a_1 \pm \hat a_1\hat a_2 \pm \hat a_3 \pm \hat a_2\hat a_3}{2}
\right\}
\end{gather*}
where we must take an even number of `$-$' signs
(so that the above set has $16$ elements).
There are three $\cE_3$-orbits, determined by real part
(two of size $4$, one of size $8$).
If $\Re(z) = -\frac12$, the set $\BL_z$ has two thin connected components
(and no thick ones).

In order to study the orbit $\Re(z) = 0$, 
take
\[ z_1 =
\grave a_1 \acute a_2 \acute a_1 \acute a_3 \acute a_2 \acute a_1 =
(-\hat a_1)\acute\eta = \acute\eta(-\hat a_3) = 
\frac{\hat a_1 - \hat a_1\hat a_2 + \hat a_3 - \hat a_2\hat a_3}{2}.  \]
A case-by-case verification shows that $\BL_{z_1}$ 
has $4$ ancestries of dimension $0$,
$3$ ancestries of dimension $1$ and
no ancestries of dimension higher that $1$. 
These numbers can also be obtained from Theorem \ref{theo:two}.
In any case,
$\BLC_{z_1}$ is
the first graph in Figure \ref{fig:4321-a12}:
$\BLC_{z_1}$ therefore collapses to a point and
$\BL_{z_1}$ is contractible.



Figure \ref{fig:4321-a12} also shows $\BLC_{z_0}$
where
\[ z_0 =
\grave a_1 \grave a_2 \acute a_1 \acute a_3 \acute a_2 \acute a_1 =
-\hat a_2\acute\eta = -\acute\eta\hat a_2 =
\frac{1+\hat a_2 + \hat a_1\hat a_3 + \hat a_1\hat a_2\hat a_3}{2}. \]
In $\BLC_{z_0}$ there are $6$ ancestries of dimension $0$,
$6$ ancestries of dimension $1$ and
exactly one ancestry of dimension $2$:
$\varepsilon = (-2,-2,+1,+1,+2,+2)$.
A computation shows that $\varepsilon$ is tame.
Thus, in the CW complex shown in Figure \ref{fig:4321-a12},
the cell of dimension $2$ glues in the obvious way,
so that $\BLC_{z_0}$ collapses to a point and is contractible.
This completes this example.
\end{example}

\begin{example}
\label{example:563412}
Set $n = 5$ and $\sigma = 563412 = a_2a_1a_3a_2a_4a_3a_5a_4a_2a_1a_3a_2$.
We have $\ell = 12$, $b = |\Block(\sigma)| = 0$, $c = \nc(\sigma) = 4$ and
\[
\acute\sigma = -\grave\sigma =
\frac{-\hat a_1-\hat a_2\hat a_3\hat a_4-\hat a_5+
\hat a_1\hat a_2\hat a_3\hat a_4\hat a_5}{2}.
\]
Consider $z_0 = \hat a_1 \acute\sigma$:
we have $\Re(z_0) = \frac12$ and $N(-z_0) = 48$.
A computation shows that $\BLC_{-z_0}$ has $56$ cells of dimension $1$,
$8$ cells of dimension $2$ and no cells of higher dimension.
This implies that the Euler characteristic is $\chi(\BL_{-z_0}) = 0$.
It turns out that $\BL_{-z_0}$ has two connected components.
We construct the CW complex corresponding to one of the components
in Figure \ref{fig:563412};
the other connected component is similar,
changing all signs for strata of codimension $0$.
The ancestries of dimension $2$ can be verified to be tame.
Thus, the $2$-cells in the CW complex fill in the hexagons and octagons
in the figure, including the hexagon which joins to top and bottom
of the figure, as in Figure \ref{fig:563412CW}.
Thus, each connected component of $\BL_{-z_0}$
is homotopically equivalent to the circle $\Ss^1$.
This completes the proof of the second item in Theorem \ref{theo:one}.
\end{example}

\begin{example}
\label{example:54321}
Set  $\sigma = \eta = a_1a_2a_1a_3a_2a_1a_4a_3a_2a_1$.
We have $\ell = 10$, $b = |\Block(\eta)| = 0$, $c = \nc(\eta) = 3$ and
\[
\acute\eta = - \grave\eta = 
\frac{-\hat a_1-\hat a_1\hat a_2\hat a_3-\hat a_4-\hat a_2\hat a_3\hat a_4}{2},
\]
It follows from Theorem \ref{theo:two} that
$N(z) = 32 + 16 \Re(z)$.
In this example, it turns out that $\acute\eta\Quat_5$
contains $4$ elements with $\Re(z) = \frac12$,
$4$ elements with $\Re(z) = -\frac12$ and
$24$ elements with $\Re(z) = 0$.
The set  $\acute\eta\Quat_5$ has $5$ orbits under $\cE_4$
of sizes $8, 4, 4, 8, 8$, shown below.
\begin{align*}
\cO_{\acute\eta} &= \left\{
\frac{\pm \hat a_1 \pm \hat a_1\hat a_2\hat a_3
\pm \hat a_4 \pm \hat a_2\hat a_3\hat a_4}{2} 
\right\},
\quad N(z) = 32,
\quad N_{\thin}(z) = 2,  \\
\cO_{\hat a_1\acute\eta} &= \left\{
\frac{1 \pm \hat a_2\hat a_3
\pm \hat a_1\hat a_4 \pm \hat a_1\hat a_2\hat a_3\hat a_4}{2}
\right\}, \quad N(z) = 40, 
\quad N_{\thin}(z) = 0,  \\
\cO_{-\hat a_1\acute\eta} &= \left\{
\frac{-1 \pm \hat a_2\hat a_3
\pm \hat a_1\hat a_4 \pm \hat a_1\hat a_2\hat a_3\hat a_4}{2}
\right\}, \quad N(z) = 24,
\quad N_{\thin}(z) = 0,  \\
\cO_{\hat a_2 \acute\eta} &= \left\{
\frac{\pm \hat a_1 \hat a_2 \pm \hat a_1 \hat a_3
\pm \hat a_2\hat a_4 \pm \hat a_1\hat a_4}{2}
\right\}, \quad N(z) = 32, 
\quad N_{\thin}(z) = 0,  \\
\cO_{\hat a_1\hat a_2\acute\eta} &= \left\{
\frac{\pm \hat a_2 \pm \hat a_3
\pm \hat a_1\hat a_2\hat a_4 \pm \hat a_1\hat a_3\hat a_4}{2}
\right\}, \quad N(z) = 32,
\quad N_{\thin}(z) = 0,  
\end{align*}
In the expressions in the Clifford algebra notation,
we must always have an even number of `$-$' signs.

\begin{figure}[ht]
\begin{center}
\includegraphics[scale=0.35]{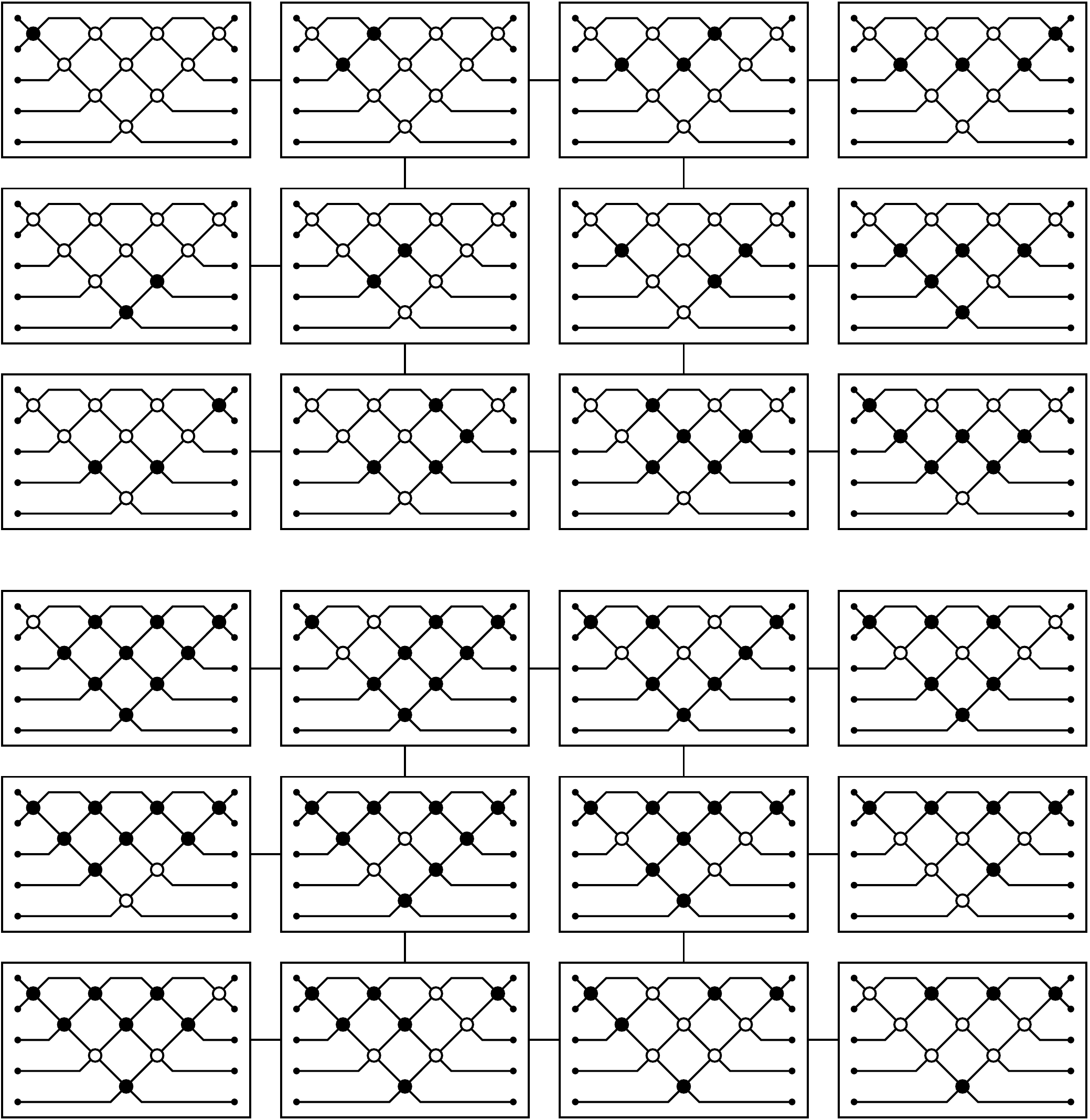}
\end{center}
\caption{The CW complex $\BLC_{-\hat a_1\acute\eta}$;
see Example \ref{example:54321}.}
\label{fig:54321x}
\end{figure}

In order to count connected components
and obtain further information about the topology
of the sets $\BL_z$, $z \in \acute\eta\Quat_{n+1}$,
we can pick one representative from each orbit
and draw the CW complex $\BLC_z$.
As a sample, we do this in Figure \ref{fig:54321x}
for $z = {-\hat a_1\acute\eta} = -\acute\eta\hat a_4$.
In this case, there are exactly two ancestries of dimension $2$,
both tame:
\begin{gather*}
(+1,-2,-2,+1,+1,+2,+1,+1,+2,+1), \\
(-1,-2,-2,+1,-1,+2,-1,+1,+2,-1);
\end{gather*}
there are no ancestries of higher dimension.
It follows that $\BL_{-\hat a_1\acute\eta}$ is homotopically equivalent
to the disjoint union of two points.
In other words, each of the two connected components 
of $\BLC_{-\hat a_1\acute\eta} = \BLC_{-\hat a_1\acute\eta,\thick}$
collapses to a point.

The other cases are more laborious
but the CW complexes can be drawn.
In all other cases, $\BL_{z,\thick}$ is non empty,
connected and contractible.
This confirms that $\BL_\eta$ has $52$ connected components
and that they are all contractible.
\end{example}

\begin{figure}[ht]
\begin{center}
\includegraphics[scale=0.35]{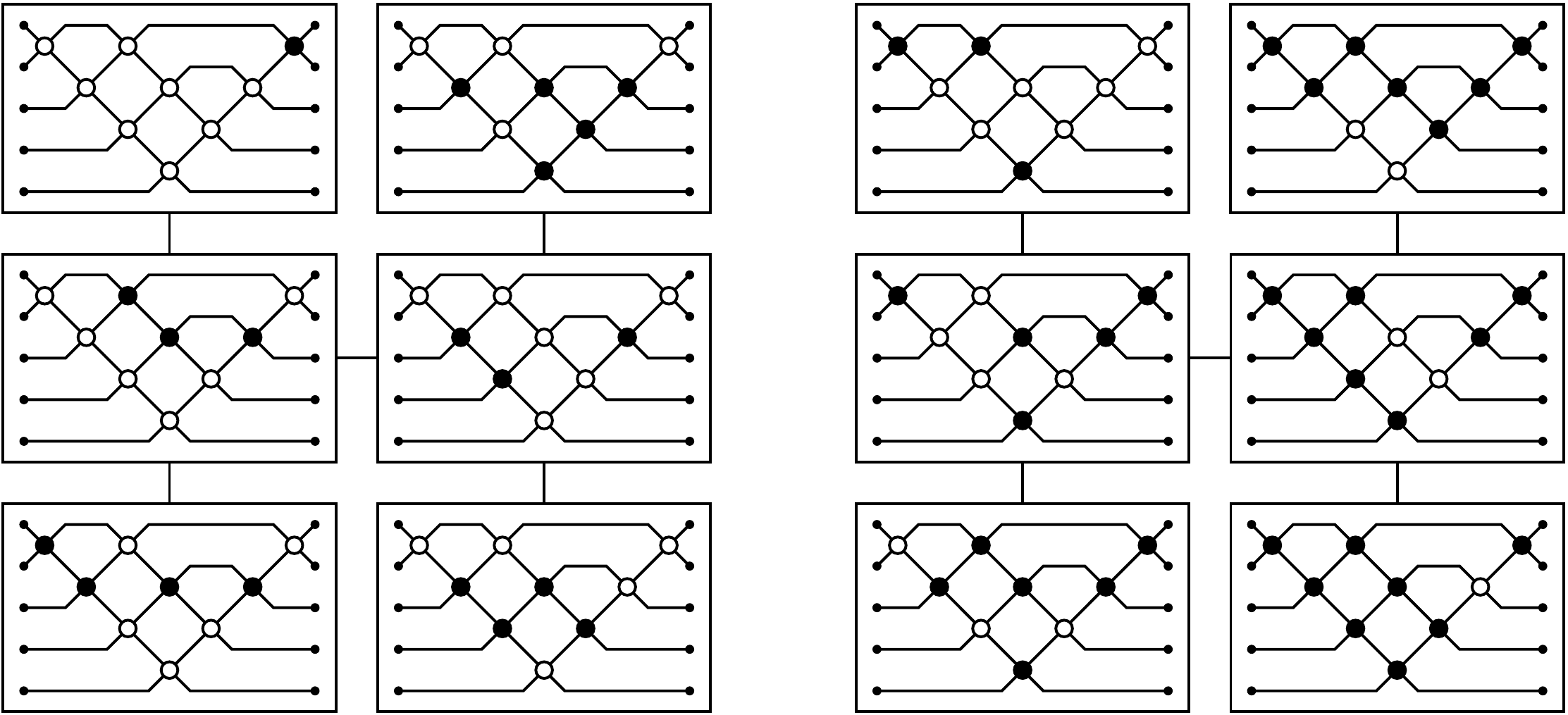}
\end{center}
\caption{The CW complex $\BLC_{-\acute\sigma\hat a_1}$;
see Example \ref{example:54231}.}
\label{fig:54231x}
\end{figure}

\begin{example}
\label{example:54231}
Consider now $\sigma = 54231 = a_1a_2a_1a_3a_2a_1a_4a_3a_2a_1$;
Figure \ref{fig:54231} shows this reduced word as a diagram.
In the cycle notation, $\sigma = (15)(243)$; 
we therefore have $n = 4$, $\ell = 9$, $c = 2$ and $b = 0$.
Theorem \ref{theo:two} tells us that, for $z \in \acute\sigma \Quat_5$,
we have $N(z) = 16 + 8\sqrt{2} \Re(z)$.
We have
\[ \acute\sigma = 
\frac{-\hat a_1 + \hat a_1\hat a_2
+ \hat a_1\hat a_3 - \hat a_1\hat a_2\hat a_3
- \hat a_4 + \hat a_2\hat a_4
+ \hat a_3\hat a_4 - \hat a_2\hat a_3\hat a_4}{2\sqrt{2}} \]
and $\Re(\pm\hat a_1\acute\sigma) = \pm\sqrt{2}/4$.

\begin{figure}[ht]
\begin{center}
\includegraphics[scale=0.25]{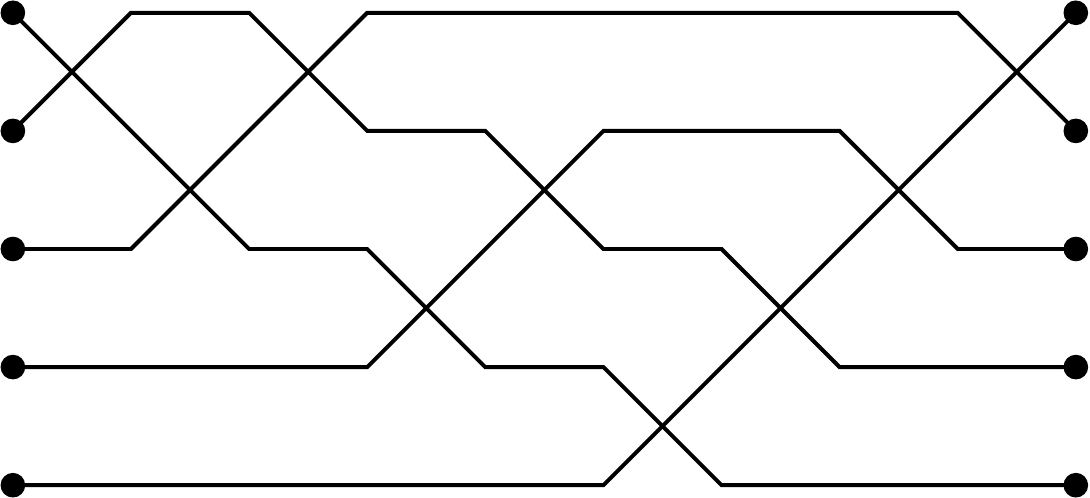}
\end{center}
\caption{The permutation $\sigma \in S_5$; see Example \ref{example:54231}.}
\label{fig:54231}
\end{figure}

It turns out that $\Re(z) = 0$ implies $c_{\anti} = 1$:
the set $\acute\sigma \Quat_5$ thus has $3$ orbits under $\cE$,
of sizes $16$, $8$ and $8$:
\begin{gather*}
\cO_{\acute\sigma},  \quad \Re(z) = 0, \quad
N(z) = 16, \quad N_{\thin}(z) = 1,  \\
\cO_{\hat a_1\acute\sigma}, \quad \Re(z) = \frac{\sqrt{2}}{4}, \quad
N(z) = 20, \quad N_{\thin}(z) = 0, \\
\cO_{-\hat a_1\acute\sigma}, \quad \Re(z) = -\frac{\sqrt{2}}{4}, \quad
N(z) = 12, \quad N_{\thin}(z) = 0.
\end{gather*}
The thick part of $\BL_{\acute\sigma}$ is connected,
so that $\BL_{\acute\sigma}$ has two connected components.
The set $\BL_{\hat a_1\acute\sigma}$ is also connected,
but $\BL_{-\hat a_1\acute\sigma}$ has two connected components.
In Figure \ref{fig:54231x} we show the two connected components
of $\BL_z$ for $z = -\acute\sigma\hat a_1$.
There are no ancestries of codimension $2$ or higher 
and therefore these connected components are contractible.
Notice that an involution takes one component of $\BL_z$ to the other.
The total number of connected components of $\BL_{\sigma}$ is therefore equal to $56$.
\end{example}

\begin{proof}[Proof of Theorem \ref{theo:collapse}]
This follows from a long computation:
samples are given in
Examples \ref{example:4231}, \ref{example:4321},
\ref{example:54321} and \ref{example:54231}.
More examples are presented in greater detail in
\cite{Alves-Saldanha-2}.
In particular, for $n = 4$ and $\sigma = \eta$,
we have a larger example than anything presented in this paper,
including cells of dimension up to $4$.
In \cite{Alves-Saldanha-2} we present an explicit sequence of collapses
taking to a CW complex homeomorphic to a $2$-disk.
Such computations rely of course on such results
as Lemma \ref{lemma:tame}.
\end{proof}

\begin{proof}[Proof of Theorem \ref{theo:one}]
The first item follows from
Theorems \ref{theo:CWcomplex} and \ref{theo:collapse}.
The second item follows from Example \ref{example:563412}.
The third item follows from Corollary \ref{coro:eveneuler} and
Remark \ref{remark:eveneuler}.
\end{proof}

\end{document}